\documentclass{amsbook}%
\usepackage{amsfonts}
\usepackage{amsmath}
\usepackage{amssymb}
\usepackage{graphicx}%
\theoremstyle{plain}

\newtheorem{case}{Case}

\newtheorem{corollary}{Corollary}

\newtheorem{lemma}{Lemma}

\newtheorem{proposition}{Proposition}
\newtheorem{remark}{Remark}

\newtheorem{theorem}{Theorem}
\numberwithin{equation}{section}

\numberwithin{figure}{part}
\renewcommand{\thefigure}{\Roman{part}-\arabic{figure}}
\begin{document}
\frontmatter
\title[Suspensions of homology spheres]{Suspensions of homology spheres}
\author[R.D. Edwards]{Robert D. Edwards}
\address{Department of Mathematics, UCLA, Los Angeles, CA 90095-1555}
\email{rde@math.ucla.edu}
\thanks{}
\date{}
\dedicatory{\vspace{5.25in}%
\_\_\_\_\_\_\_\_\_\_\_\_\_\_\_\_\_\_\_\_\_\_\_\_\_\_\_\_\_\_\_\_\_\_\_\medskip\\This research was supported in part under NSF Grant No. MCS 76-06903.
Preparation of the electronic manuscript was supported by NSF Grant
DMS-0407583. Final editing was carried out by Fredric Ancel, Craig Guilbault
and Gerard Venema.}
\begin{abstract}
This article is one of three highly influential articles on the topology of
manifolds written by Robert D. Edwards in the 1970's but never published. This
article \textquotedblleft Suspensions of homology spheres\textquotedblright%
\ presents the initial solutions of the fabled Double Suspension Conjecture.
The article \textquotedblleft Approximating certain cell-like maps by
homeomorphisms\textquotedblright\ presents the definitive theorem on the
recognition of manifolds among resolvable generalized manifolds. (This work
garnered Edwards an invitation to give a one-hour plenary address to the 1978
International Congress of Mathematicians.) The third article \textquotedblleft
Topological regular neighborhoods\textquotedblright\ develops a comprehensive
theory of regular neighborhoods of locally flatly embedded topological
manifolds in high dimensional topological manifolds. The manuscripts of these
three articles have circulated privately since their creation. The organizers
of the Workshops in Geometric Topology
(http://www.math.oregonstate.edu/\symbol{126}topology/workshop.htm) with the
support of the National Science Foundation have facilitated the preparation of
electronic versions of these articles to make them publicly available.\medskip

This article contains four major theorems:

\begin{enumerate}
\item[I.] The double suspension of Mazur's homology 3-sphere is a 5-sphere,

\item[II.] The double suspension of any homology n-sphere that bounds a
contractible (n+1)-manifold is an (n+2)-sphere,

\item[III.] The double suspension of any homology 3-sphere is the cell-like
image of a 5-sphere.

\item[IV.] The triple suspension of any homology 3-sphere is a
6-sphere.\medskip

\end{enumerate}

Edwards' proof of I was the first evidence that the suspension process could
transform a non-simply connected manifold into a sphere, thereby answering a
question that had puzzled topologists since the mid-1950's if not earlier.
Results II, III and IV represent significant advances toward resolving the
general double suspension conjecture: the double suspension of every homology
n-sphere is an (n+2)-sphere. [The general conjecture was subsequently proved
by J. W. Cannon (Annals of Math. 110 (1979), 83-112).]

\end{abstract}

\maketitle
\tableofcontents

\centerline{\large{\bf Introduction}}\bigskip

In the study of triangulations of manifolds, there arises naturally the
following.\bigskip

\textsc{Suspension Question.} \emph{Is there any manifold }$H^{n}$\emph{,
other than the }$n$\emph{-sphere}$\footnote{Throughout this paper,
superscripts on topological spaces always denote dimension (as is customary).
Their occasional appearance or disappearance in later sections has no
notational significance, e.g., $K^{2}$ and $K$ denote the same space.}$\emph{,
such that for some }$k\geq1$\emph{, the }$k$\emph{-fold suspension }$%
\Sigma^kH^n$\emph{ is homeomorphic to the }$(n+k)$\emph{-sphere? }(see
Prologue, Section I; definitions are below).\bigskip

By Alexander duality applied to the $(k-1)$-dimensional \textquotedblleft
suspension sphere\textquotedblright\ $S^{k-1}$ in $\Sigma^{k}H^{n}$, coupled
with the fact that $\Sigma^{k}H^{n}-S^{k-l}$ is homeomorphic to $H^{n}\times
H^{k}$, such a manifold $H^{n}$ satisfying the suspension question must
necessarily be a \textit{homology} n-\textit{sphere}, which means a closed
topological $n$-manifold having the integral homology groups of the $n$-sphere
(this explains the customary letter $H$). This paper examines the sufficiency
of this condition. The potential singular set (i.e. nonmanifold set) of
$\Sigma^{k}H^{n}$ is along the suspension $(k-1)$-sphere. There one readily
sees that $\Sigma^{k}H^{n}$ is locally homeomorphic to $\mathring{c}%
H^{n}\times\mathbb{R}^{k-1}$, where $\mathring{c}H^{n}$ denotes the open cone
on $H^{n}$. It is known that the suspension question really amounts to that of
deciding whether $\mathring{c}H^{n}\times\mathbb{R}^{k-1}$ is a topological
manifold. This is because:

\begin{quotation}
\emph{If the suspension of a compact space }$Y$\emph{ is a manifold, then
}$\Sigma Y$\emph{ must be a sphere.}
\end{quotation}

\noindent The quickest proof of this fact is to argue that the complement in
$\Sigma Y$ of a suspension point (which is homeomorphic to $\mathring{c}Y)$
must be homeomorphic to euclidean space, since every compact subset of it lies
in an open ball (details of this are in \cite{Brwn}).

The $k=2$ case of the above suspension question has come to be known as the
Double Suspension Conjecture: \emph{If }$H^{n}$\emph{ is a homology }%
$n$\emph{-sphere, then} $\Sigma^{2}H^{n}\approx S^{n+2}$ (where $\approx$
means \textquotedblleft is homeomorphic to\textquotedblright). The
unrestricted version of this conjecture (i.e. allowing $k>2$) is sometimes
called the Multiple Suspension Conjecture. (I do not know why the Double
Suspension Conjecture is a conjecture, rather than a question or problem; that
is the way most people refer to it now. This conjecture has been raised as a
question in various places, for example \cite{CZ} and \cite{Gl1}-\cite{Gl4};
probably the best known place is Milnor's list of problems \cite[ p. 579]%
{Las}. Earlier it was referred to as the Sphere Problem \cite[p. 16]{Moi}. The
suspension question must have confronted all researchers who since the time of
Brouwer's work \cite{Brou} have attempted to understand triangulations of
manifolds. For example, see \cite{Cai1} or \cite{Cai2} (especially $\S $ 6);
the latter work has a particularly good bibliography. See also \cite{Ku}.

There exist many homology spheres, which are not themselves homeomorphic to
spheres, on which the conjecture can be tested. Perhaps the most famous is
Poincar\'{e}'s binary dodecahedral homology 3-sphere, but unfortunately it
turns out to be somewhat difficult to work with, as we will see. There exist
in profusion more tractable homology n-spheres in all dimensions $n\geq3$ (see
Parts I and II).

All known genuine (i.e. nontrivial) homology n-spheres are
nonsimply-connected, necessarily so when $n\geq5$, for in these dimensions the
generalized topological Poincar\'{e} conjecture \cite{New1} establishes that a
simply-connected homology n-sphere is homeomorphic to a sphere. This explains
why the $k=1$ case of the suspension question is not emphasized, for in that
case an additional necessary condition on $H^{n}$ is that it be
simply-connected (assuming $n\geq2$). Hence that case of the question can have
nontrivial content only possibly when $n=3$ or $4$. If $H^{4}$ is a homotopy
4-sphere, it is known that $\Sigma H^{4}\approx S^{5}$ (see \cite[Assertion p.
83]{Si3} or \cite[Appendix I]{Si6}; compare \cite{Hi2}); if $H^{3}$ is a
homotopy $3$-sphere, it is \textbf{not} known whether $\Sigma H^{3}\approx
S^{4}$, but it is known that $\sigma^{2}H^{3}\approx S^{5}$ (\cite[Theorem
A]{Si3}; compare \cite{Gl1}).

These results represent the first cases of the suspension conjecture ever
settled (possibly vacuously, of course). For nonsimply connected homology
spheres, the double suspension conjecture remained one of those tantalizing
problems that had three possible outcomes: it could be true never, sometimes,
or always. Now, after several stages of development due to myself and J.
Cannon, the conjecture has been completely settled in the affirmative by
the\bigskip

\textsc{Double Suspension Theorem.} \textit{The double suspension $\Sigma
^{2}H^{n}$ of any homology sphere $H^{n}$ is homeomorphic to a sphere.\bigskip
}

The purpose of this paper is to present all of my work related to the
suspension question, which took place during the period 1974-76. Briefly
stated, I established as a general rule that $\Sigma^{k}H^{n}\approx S^{n+k}$
whenever $n+k\geq6$ $(k\geq2)$. Recently, Cannon crystallized and extended my
work to prove a much more general theorem, one consequence of which was the
final case $\Sigma^{2}H^{3}\approx S^{5}$ \cite{Can2}.

This paper is written in four parts, elaborated in the manner I have always
intended, corresponding to the four stages in the development of my
understanding of the problem. The complete contents of the paper are as
follows:\bigskip

\noindent\textsc{Prologue:}\textbf{ }\textit{Consequences of the Double
Suspension Theorem; plus preliminaries, definitions and notation. \medskip}

As motivation for the main results presented in Parts I-IV, a number of
consequences of the double suspension are presented. In addition, some
preliminaries are discussed; specifically, cell-like decompositions and the
Bing Shrinking Theorem. Finally, some definitions and notation are
established.\bigskip

\noindent\textsc{Part I:} \textit{The double suspension of Mazur's homology
3-sphere is} $S^{5}$.\textit{\medskip}

This represents the first case of the conjecture settled for a
nonsimply-connected homology sphere. The proof amounts to a self-contained,
bare-hands construction of a discernible homeomorphism.\bigskip

\noindent\textsc{Part II:} \textit{The double suspension of any homology
$n$-sphere which bounds} \textit{a contractible $(n+1)$-manifold is a
sphere.\medskip}

In particular, this holds for any homology $n$-sphere, $n\geq4$. This part
generalizes the construction in Part I. It makes essential use of an important
construction introduced by \v{S}tan'ko in \cite{St1}. A Postscript to Part II
explains a certain Replacement Principle for Cell-like Compacta in manifolds.
An Appendix to Part II explains how to shrink the decomposition which is at
the heart of Parts I and II, namely, the spun Bing decomposition.\bigskip

\textsc{Part III:} \textit{The double suspension of any homology $3$-sphere is
the image of $S^{5}$ under a cell-like map.\medskip}

Thus, the suspension question for an arbitrary homology $3$-sphere is a
cell-like decomposition space problem. This result was proved independently by
Cannon in \cite{Can3}. This part contains in passing a simplified construction
of Siebenmann's non PL-triangulable topological manifold, extracted from
\cite{Sch}.\bigskip

\noindent\textsc{Part IV:} \textit{The triple suspension of any homology
$3$-sphere is $S^{6}$.\medskip}

This ad hoc argument rests in the end on a clever construction of Bing, which
he used to establish the shrinkability of a countable collection of flat arcs
in a manifold.\bigskip

Although these results are being formally published here for the first time,
various informal lecture notes from my talks, taken by various people, have
been around for some time, e.g. from Orsay (Spring 1975), Cambridge (July
1975), Nantes (Spring 1976), Institute for Advanced Study (October 1976), and
the AMS St. Louis meeting (January 1977). I am grateful in particular to A.
Fathi, L. Guillou and Y. M. Vissetti for their clear write-up of my Orsay talks.

\mainmatter

\part*{Prologue: Consequences of the Double Suspension Theorem; plus
preliminaries, definitions and notation}

\section{Consequences of the Double Suspension Theorem}

This section is a compilation of various corollaries, most of them known to
various people before the theorem was proved. (In the seven month period
before Cannon's work, when only my triple suspension theorem (Part IV) was
known as a general rule, some of the results here were restricted to ambient
dimension $\geq6$, rather than the present $\geq5$,)

\subsection{\textbf{Exotic triangulations}}

The most noteworthy corollary is:

\begin{theorem}
\textit{Not all triangulations of topological manifolds need be piecewise
linear ($PL$; i.e. combinatorial) triangulations}. \textit{In fact, for any
given dimension $q\geq5$, there exists a triangulated topological manifold
$Q^{q}$ which is not even homotopically equivalent to any $PL$ manifold}.
\end{theorem}

\textsc{Discussion:} A $PL$ manifold-without-boundary can be defined to be a
polyhedron $P$ which is piecewise linearly homogeneous, that is, for any two
points $x,y\in P$, there is a $PL$ homeomorphism $h:P\rightarrow P$ such that
$h(x)=y$. In this spirit, a triangulated topological manifold can be defined
to be polyhedron which is \textbf{topologically} homogeneous. The point, then,
of the first part of the assertion above is that there exists a topologically
homogeneous polyhedron which is not piecewise linearly homogeneous. Such an
example is provided by $\Sigma^{2}H^{n}\approx S^{n+2}$ for any nontrivial
polyhedral homology sphere $H^{n}$, e.g. Newman's, Poenaru's or Mazur's (see
Part I).

Still, this example leaves something to be desired, since $\Sigma^{2} H^{n}$
is topologically homeomorphic to a piecewise linearly homogeneous polyhedron.
This leads to the second part of the assertion.

The most easily described example of the polyhedron $Q$ is the familiar one
$Q^{q}=(M^{4}\cup c(\partial M^{4}))\times T^{q-4}$, where $T^{q-4}$ denotes
the $(q-4)$-torus, and where $M^{4}$ is the 4-manifold-with-boundary described
by Milnor as being eight copies of the tangent disc bundle of $S^{2}$ plumbed
together according to the $E_{8}$-diagram (see \cite[Chap 5]{Brwd} or
\cite[\S 8]{HNK}), or (equivalently) described by Hirzebruch as the resolution
of the singularity of the equation $z_{1}^{2}+z_{2}^{3}+z_{3}^{5}=0$,
restricted to the unit ball in $C^{3}$ \cite[$\S $9, esp. Exercise (5.8)]%
{HNK}. The manifold $M^{4}$, being smooth, can be given $PL$ manifold
structure. The boundary $\partial M^{4}$ is Poincar\'{e}'s binary dodecahedral
homology 3-sphere. The significance of the manifold-with-singularity
$M^{4}\cup c(\partial M^{4})$ is that it is parallelizable off of the
conepoint, and it has index 8. This contrasts with the bedrock theorem of
Rokhlin that any closed $PL$ 4-manifold which is parallelizable off of a point
must have index divisible by 16.

If $Q$ were homotopically equivalent to a $PL$ manifold, $W$ say, then one
could homotope the map
\[
W\ \overset{\text{hom.equiv}}{\longrightarrow}\ Q\overset{\text{proj.}%
}{\longrightarrow}\ T^{q-4}%
\]
to be $PL$ and transverse to a point in $T^{q-4}$, producing as point-preimage
a closed $PL$ 4-manifold $N^{4}$. One argues that $N^{4}$ would be
almost-parallelizable, and also would have index 8, contradicting Rokhlin's
theorem. See \cite[Section 5]{Si2} and \cite[IV, App. B]{KS2} for a discussion
of this and related details, and for other references.

The $PL$ singular set of $Q^{q}$ is the subset $c\times T^{q-4}$, along which
$Q^{q}$ is locally homeomorphic to $\mathring{c}(\partial M^{4})\times
\mathbb{R}^{q-4}$. By the Double Suspension Theorem, $Q^{q}$ is locally
euclidean there, and hence is a topological manifold.

The above simple description of the manifold $Q^{q}$ rests of course on the
Double Suspension Theorem and all the work it entails. A topological manifold
homotopically equivalent (now known homeomorphic) to $Q^{q}$, for any $q\geq
5$, was first constructed in 1969 by Siebenmann \cite[Section 5]{Si2}{, as a
counterexample to the $PL$ triangulation conjecture for topological manifolds.
His construction, which rests in the end only on the local contractibility of
the homeomorphism group of a topological manifold, is the starting point of
the work in Part III. }

This discussion closes with some quick remarks related to $Q^{q}$ and to
triangulations of manifolds.\bigskip

\textsc{1. Triangulating topological manifolds.}\quad The question of whether
all topological manifolds are triangulable (i.e., homeomorphic to polyhedra)
remains open. Completing the line of investigation begun in \cite{Si3}, D.
Galewski - R. Stern \cite{GS2} and T. Matumoto \cite{Mat2} have established
the following stimulating result (incorporating now the Double Suspension
Theorem into their work): \textit{All topological manifolds of dimension
$\geq5$ are triangulable if and only if there exists a homology $3$-sphere
$H^{3}$, having Rokhlin invariant} $1\in\mathbb{Z}/2$ (meaning: $H^{3}$ bounds
a parallelizable $PL$ 4-manifold of index 8), \textit{such that $H^{3}\#H^{3}$
bounds a $PL$ acyclic $4$-manifold (i.e., such that $H^{3}$ is $PL$
homology-cobordant to $-H^{3}$)}. In fact, there is a specific, easily
constructed \textquotedblleft test\textquotedblright\ topological
5-manifold-without-boundary $M^{5}$, which is triangulable if and only if such
a homology 3-sphere exists.\bigskip

\textsc{2. Algebraic varieties.}\quad It is commonly known that $\Sigma
^{k}H^{n}$ is a real algebraic variety for many homology spheres $H^{n}$, e.g.
Poincar\'{e}'s $H^{3}$, $k$ arbitrary. Recently S. Akbulut and R. Mandelbaum
remarked that the non-PL-triangulable manifold $Q^{q}$ described above can be
realized as a real algebraic variety (or complex algebraic variety, if $q$ is
even) (see \cite{AK}).\bigskip

\textsc{3. Stratified spaces.}\quad As already shown, the topological space
$Q^{q}$ is an example of a space whose intrinsic topological stratification
has only one stratum, yet for \textbf{any} polyhedral structure on $Q$ (even
up to homotopy), its intrinsic $PL$ stratification must have at least two
strata. D.R. Anderson constructed in \cite{An} an example of such a space $Q$
having no $PL$ stratification as simple as its intrinsic topological
stratification, without using any suspension theorems for homology
spheres.\bigskip

\textsc{4. Topological transversality.}\quad Assuming the Double Suspension
Theorem true for a single homology 3-sphere $H^{3}$ of Rokhlin invariant 1,
e.g. Poincar\'{e}'s $H^{3}$, M. Scharlemann established the following
transversality theorem at dimension 4 (this is the nonrelative version; for
complete statements, see \cite[Thms. B and C]{Sch}): \textit{Given a map
$f:M^{m}\rightarrow E$ from a topological $m$-manifold to the total space $E$
of a euclidean microbundle $\xi:X\hookrightarrow E\rightarrow X$ of fiber
dimension $m-4$ over a space $X$, then there is an arbitrarily small homotopy
of $f$ to a map $f_{\ast}$ such that $f_{\ast}$ is microbundle-transverse to
the core $X$ of $E$, and such that $f_{\ast}^{-1}(X)$ is a $4$%
-manifold-with-isolated-singularities, each singularity being homeomorphic to
$\overset{\circ}{c}H^{3}$.} Recall that if $m-k\neq4\neq m$, where $k=$ fiber
dimension of $E$, then Kirby-Siebenmann \cite[III, Section 1]{KS2} proved the
best possible topological transversality theorem, which concludes that
$f_{\ast}^{-1}(X)$ is an $(m-k)$-manifold, without singularities.

\textsc{5. Collapsibility and shellability.}\quad In \cite{BCC}, the authors
prove the following curious result, using the fact that the double suspension
of some nonsimply connected homology $n$-sphere $H^{n}$ is topologically a
sphere: \textit{For every $n\geq3$, there is a polyhedron $P^{n+2}$,
topologically homeomorphic to the $(n+2)$-ball, such that $P\times I^{n-2}$ is
not collapsible, but $P\times I^{n+3}$ is collapsible.} The example is
$P=\Sigma^{2}H^{n}-\operatorname*{int}B^{n+2}$, where $B^{n+2}$ is any $PL$
cell in $\Sigma^{2}H^{n}$ disjoint from the suspension circle. In a similar
vein, $\Sigma^{2}H^{n}$ provides an example of a nonshellable triangulation of
a sphere (see \cite{DK}).\bigskip

\subsection{\textbf{Polyhedral homology manifolds.}}

Another consequence of the Double Suspension Theorem is the second half of the
following assertion (the first half being already known).

\begin{theorem}
\textit{Suppose $P$ is a connected polyhedron. If $P$ is a topological
manifold, then $P$ is a homotopy manifold (as defined below). Conversely, if
$P$ is a homotopy manifold and $\dim P\neq4\neq\dim\partial P$, then $P$ is a
topological manifold.}
\end{theorem}

The dimension restriction arises because it is not known whether the cone on a
homotopy 3-sphere is a manifold.

Two special cases of the above assertion are as follows:

\begin{corollary}
\textit{Suppose $P_{1},P_{2}$ are two connected polyhedra, with $\dim
P_{i}\neq0\neq\dim\partial P_{i}$ (the latter nonequality meaning that $P_{i}$
contains no open set homeomorphic to $[0,1))$. Then}

\begin{enumerate}
\item \textit{the product $P_{1}\times P_{2}$ is a topological manifold
$\Longleftrightarrow P_{1}$ and $P_{2}$ are homology manifolds with
homotopically collared boundaries, and}

\item \textit{the join $P_{1}\star P_{2}$ is a compact topological manifold
(necessarily a ball or sphere) $\Longleftrightarrow P_{1}$ and $P_{2}$ are
compact homology manifolds with homotopically collared boundaries, such that
each of $P_{1}$ and $P_{2}$ either has the integral homology groups of a
sphere, or else is contractible.}
\end{enumerate}
\end{corollary}

A \textit{polyhedral homology manifold} is a polyhedron $P$ in which the link
of each component of each stratum has the integral homology groups of either a
sphere or a point. Using duality, this property holds if it is known to hold
only for the links of those strata-components which are closed subsets of $P$.
The union of the strata-components whose links are acyclic comprise a
subpolyhedron of $P$, called the \textit{boundary} of $P$ and denoted
$\partial P$; it can be shown to be a homology manifold-without-boundary.

The strata referred to in this discussion are from any $PL$ stratification of
$P$ compatible with the $PL$ structure of the polyhedron. For example, one
could let the $i^{\text{th}}$ stratum be $K^{(i)}-K^{(i-1)}$, where $K$ is a
simplicial complex which is $PL$ homeomorphic to $P$, and where $K^{(j)}$
denotes the $j$-skeleton of $K$. The most natural stratification of a
polyhedron $P$ is the minimal one, i.e., the intrinsic $PL$ stratification
\cite[p. 421]{Ak}. A \textit{vertex} of a stratification is simply a point of
the 0-stratum; this stratum may be empty. If one prefers, one may talk about
simplicial homology manifolds instead of polyhedral homology manifolds, but
the polyhedral setting is the natural one (just as the notion of a PL manifold
is a more natural notion than that of a combinatorial manifold).

Much of this material is explained more fully, from the simplicial standpoint,
in \cite[Chapter V]{Mau}.

If a polyhedral homology manifold $P$ is to have the property that
$P\times\mathbb{R}^{n}$ is a topological manifold for some euclidean space
$\mathbb{R}^{n}$, then $P$ must have the additional property that its acyclic
links are in fact contractible (this condition is vacuous if $\partial
P=\emptyset)$. Let such a $P$ be called a \textit{polyhedral homology manifold
with homotopically collared boundary.} The reason for this name is that it can
be verified that a polyhedral homology manifold $P$ satisfies this additional
property if and only if there is a \textquotedblleft homotopy
collaring\textquotedblright\ map $\psi:\partial P\times\lbrack0,1)\rightarrow
P$ extending the inclusion $\partial P\times0=\partial P\hookrightarrow P$
such that $\psi(\partial P\times(0,1))\cap\partial P=\emptyset$. As an
example, the polyhedron $P=\Sigma^{k}M$, for any acyclic compact $PL$ manifold
$M$ and any $k\geq1$, is a polyhedral homology manifold. However, its boundary
is homotopically collared if and only if $M$ is simply-connected.

If a polyhedral homology manifold $P$ is itself to be a topological manifold,
then in addition it must have the property that the links in both $P$ and
$\partial P$ of vertices (as defined above) must be simply connected whenever
they have dimension $\geq2$. Let a polyhedron $P$ which is a polyhedral
homology manifold with homotopically collared boundary, and which satisfies
this link condition, be called a \textit{polyhedral homotopy manifold} (this
terminology may contrast with some earlier usage, but it seems justified, here
at least).

\begin{proof}
[Proof of Theorem]In the following discussion, all stratifications of
polyhedra will be assumed to be intrinsic stratifications, in order to
minimize hypotheses.

The basic case of the assertion is when $P$ is a homotopy
manifold-without-boundary. The proof in this case proceeds by induction on the
\textit{depth} of $P$, which is defined to be $\dim P-\dim P_{0}$, where
$P_{0}$ is the nonempty stratum of $P$ of lowest dimension (see for example
\cite{Si5}). For induction purposes, then, assume $n\geq0$, and assume it has
been established that given any polyhedral homotopy manifold-without-boundary
$Q$, with depth $Q\leq n$ and $\dim Q\neq4$, then $Q$ is a topological
manifold-without-boundary. Suppose $P$ is a polyhedral homotopy
manifold-without-boundary with depth $P=n+1$ and $\dim P\neq4$. Let $P_{0}$ be
the nonempty stratum of $P$ of lowest dimension. Then depth $(P-P_{0})\leq n$,
and hence by the induction hypothesis $P-P_{0}$ is a topological
manifold-without-boundary. Let $k=\dim P_{0}$, and let $L$ be the link of any
component of $P_{0}$. In order to show that $P$ is a manifold, it suffices to
establish:\medskip

\noindent\textbf{Claim:} \textit{$\overset{\circ}{c}L\times\mathbb{R}^{k}$ is
a topological manifold.}

The proof of this is divided into several cases, depending on the dimension of
$L$.

\begin{case}
$\dim L\leq2$.
\end{case}

Then $L$ is topologically a sphere, hence the claim follows.

\begin{case}
$\dim L=3$.
\end{case}

Then $L$ is a (manifold) homology 3-sphere. Hence the claim follows from the
Double Suspension Theorem (note that $k\geq1$ in this case, since $\dim
P\neq4).$

\begin{case}
$\dim L\geq4$.
\end{case}

We have that depth$(L\times\mathbb{R}^{1})$ (= depth $L)<$ depth $P$, and
hence by the induction hypothesis the polyhedral homotopy
manifold-without-boundary $L\times\mathbb{R}^{1}$ is a topological manifold.
By the Proposition in Part II, $L\times\mathbb{R}^{1}$ is embeddable as a
neighborhood of the end of some open contractible topological manifold $M$. If
$k=0$, then by hypothesis $L$ is homotopically a sphere, hence $M$ is
homeomorphic to $\mathbb{R}^{m}$ (\cite{Si1}, or \cite{St} coupled with the
fact that $M$ has $PL$ manifold structure). Hence $\overset{\circ}{c}L$, which
is homeomorphic to the 1-point semicompactification of the end, is a manifold.
If $k\geq1$, then $\overset{\circ}{c}L\times\mathbb{R}^{k}$, being
homeomorphic to $M/X\times\mathbb{R}^{1}$, is a manifold by Part II, where
$X=M-L\times\mathbb{R}^{1}$. This establishes the basic case of the
Theorem.\newline(Technical aside: When only the triple suspension theorem was
known for homology 3-spheres, and hence the above assertion required the
additional restriction $\dim P\neq5\neq\dim\partial P$, the last part of this
argument, for the case $\dim L=4$ and $k\geq2$, had to be established by a
more complicated argument.)

In the general case, when $\partial P\neq\emptyset$, the preceding argument
can be augmented by using the following collar argument, which avoids some
delicate link analysis at $\partial P$. Let $Q$ be gotten from $P$ by
attaching to $P$ an exterior boundary collar, this being denoted
\[
Q=(P+\partial P\times\lbrack0,1])/\partial P=\partial P\times0.
\]
Hence $\partial Q=\partial P\times1$. By the without-boundary case, both
$\partial Q$ and $Q-\partial Q$ are manifolds, hence $Q$ is a manifold.
Furthermore the given collar $\partial P\times\lbrack0,1]$ for $\partial Q$ in
$Q$ is 1-LCC (1-locally co-connected) in $Q$, that is, small loops in
$Q-\partial P\times\lbrack0,1]$ (sizes measured in the metric of $Q$) are
null-homotopic in $Q-\partial P\times\lbrack0,1]$ by small homotopies. Hence,
a now-standard radial engulfing argument (\cite{Se}; see \cite[\S 3]{Da3})
shows that the interval-fibers of this collar in $Q$ can be shrunk to points
by pseudoisotopy of $Q$, and hence $Q$ is homeomorphic to the quotient space,
which is $P$.
\end{proof}

\begin{remark}
It can be deduced from the recent work of Cannon \cite{Can2} combined with
Cannon \cite{Can3} or Bryant-Lacher \cite{BL2} that the purely topological
analogue of the assertion above is true. That is, it remains true with
\emph{polyhedron} replaced by its topological analogue, namely a
\emph{conelike-stratified (CS) set} as defined by Siebenmann in \cite{Si5}.
Although links in cone-like stratified sets are not intrinsically
well-defined, they become intrinsically well-defined after crossing them with
euclidean space of dimension equal to their codimension. Hence the homotopy
type of \textquotedblleft the\textquotedblright\ link makes sense. The
existence of intrinsic stratifications for conelike-stratified sets was
establishes by M. Handel in \cite{Han}.
\end{remark}

\subsection{\textbf{Wild embeddings with mapping cylinder neighborhoods.}}

Another consequence of the Double Suspension Theorem is:

\begin{theorem}
\textit{For any given dimension $m\geq5$, there exists a topological embedding
of a circle into an $m$-manifold such that the embedding has a manifold
mapping cylinder neighborhood, and yet the embedding is not locally flat.}
\end{theorem}

Such an example cannot exist in dimensions $m\leq4$ (see \cite{BL2}). The
example is provided by the suspension circle embedded in $\Sigma^{2}%
H^{n}\approx S^{n+2}$, for any non-simply-connected homology sphere $H^{n}$.
This particular embedding is homogeneous by ambient isotopy. Note that this
shows, for example, that there is no reasonable notion of an `intrinsically'
good point of an embedding (one might conjecture such a notion after learning
statements such as \textquotedblleft any circle in euclidean space is tame
modulo a 0-dimensional $F_{\sigma}$ subset\textquotedblright).

\subsection{Codimension two embeddings.}

Since the work of Kirby-Siebenmann, it has been known that there exists a
topologically locally flat codimension 2 embedding of one $PL$ manifold into
another, such that the embedding cannot be ambient isotoped to be piecewise
linearly locally flat. For example, one can take the inclusion $S^{3}%
=S^{3}\times0\hookrightarrow(S^{3}\times\mathbb{R}^{2})_{\theta})$, where
$(S^{3}\times\mathbb{R}^{2})_{\theta}$ denotes the $PL$ manifold obtained by
putting the nonstandard $PL$ structure $\theta$ on $S^{3}\times\mathbb{R}^{2}$
\cite{KS1}. The Double (or Multiple) Suspension Theorem provides another
example, for given any homology 3-sphere $H^{3}$, and any $PL$ manifold
$M^{m},m\geq0$, \textit{the homeomorphism $\Sigma^{2}H^{3}\times M^{m}\approx
S^{5}\times M^{m}$ can be chosen to be piecewise linearly locally flat on}
$H^{3}\times M^{m}\Leftrightarrow$ \textit{the Rokhlin invariant of} $H^{3}$
is 0. The proof will not be given here; it uses the same argument that was
used in Section I above, involving transversality and Rokhlin's theorem.

\subsection{\textbf{Exotic group actions.}}

It is well known that the Double Suspension Theorem establishes that
\textit{there is a nonstandard semifree action of $S^{1}$ on $S^{5}$ which is
piecewise linear with respect to some polyhedral structure on $S^{5}$.} For
there is a natural semifree action of $S^{1}$ on $\Sigma^{2}H^{3}=S^{1}\ast
S^{3}\approx S^{5}$ with fixed point set $H^{3}$. Another interesting group
action, pointed out to me by F. Raymond and H. Samelson, is the following:
\textit{there is a topological action of $SO(3)$ on $S^{7}$ with all isotropy
groups discrete.} (Compare \cite{MS} and \cite{Ol}; the latter work in fact
exhibits such an action which is smooth.) This action follows from the fact
that $S^{7}\approx H^{3}\ast H^{3}$, where $H^{3}$ is Poincar\'{e}'s binary
dodecahedral homology 3-sphere (cf. II above). The action of $SO(3)$ is the
diagonal action.

\section[Preliminaries: Cell-like decompositions]{Preliminaries: Cell-like
decompositions and the Bing Shrinking Criterion}

The suspension problem for homology spheres is a problem in the theory of
cell-like upper semicontinuous decompositions of manifolds. This venerable
subject, fathered by R.L. Moore and developed largely by R.H. Bing, studied
the class of proper cell-like maps $\{f:M\rightarrow Q\}$ from manifolds onto
metric spaces (definitions below). The major problem of the subject is
\textit{to decide when the quotient} (i.e. {target) space $Q$ is a manifold,
or at least when $Q\times R^{k}$ is a manifold, for some $k$.}

The link between the suspension question for homology spheres and cell-like
decomposition theory became apparent as a result of the Newman, Poenaru and
Mazur constructions of homology $n$-spheres which bound contractible
$(n+1)$-manifolds (recalled in Part I). For suppose $H^{3}=\partial M^{4}$,
where $M^{4}$ is a contractible 4-manifold. Let $X$ be a spine of $M$, that
is, let $X=M-\partial M\times\lbrack0,1)$, where $\partial M\times\lbrack0,1)$
denotes any open collar for $\partial M$ in $M$. Then $\Sigma^{k}H^{3}$ is a
sphere if and only if $M/X\times\mathbb{R}^{k-1}$ is a manifold (as explained
in the Introduction), where $M/X$ denotes the quotient space of $M$ gotten by
identifying $X$ to a point; clearly $M/X\approx c(\partial M)$. So the
suspension question becomes one of determining whether the target of the
cell-like map $M\rightarrow M/X$ is stably a manifold.

In the subject of cell-like decompositions of manifolds, there are two
landmark results of a general character, the Shrinking Theorem of Bing, and
the Cellular Approximation Theorem of Armentrout and Siebenmann. Both concern
the question of approximating certain maps by homeomorphisms.

The version of Bing's theorem which I prefer is the following if-and-only-if
version. In applications of this theorem, invariably $M$ is a manifold, $f$ is
a cell-like surjection and the problem is to decide whether $Q$ is a manifold.
The significance of Bing's theorem was to turn attention from the target space
$Q$, where it had been focused since R.L. Moore's work \cite{Mo} on cell-like
decompositions of the plane, to the source space $M$, where one had the
obvious advantage of working in a space known to be a manifold. The
justification for the approximation statement in the theorem will become clear
after the next theorem.\bigskip

\textsc{Shrinking Theorem} (Bing 1952, from \cite[Section 3 II, III]{Bi1};
compact version). \textit{A surjection $f:M\rightarrow Q$ of compact metric
spaces is approximable by homeomorphisms if and only if the following
condition holds (now known as the Bing Shrinking Criterion): Given any
$\epsilon>0$, there exists a homeomorphism $h:M\rightarrow M$ such that}

\begin{enumerate}
\item distance $(fh,f)<\epsilon,$ \textit{and}

\item \textit{for each} $y\in Q$, diam$(h(f^{-1}(y)))<\epsilon$.\bigskip
\end{enumerate}

There are many refinements and addenda one can make to this theorem, the most
significant having to do with realizing $f$ by pseudoisotopy versus realizing
$h$ by ambient isotopy. These will be not discussed here. Bing used the
theorem in much of his subsequent work, including \cite[proofs of
Theorems]{Bi2}, \cite[Section 8]{Bi3} and \cite[\S 3 and Thm2]{Bi4}. McAuley
\cite{McA1} was the first person to broaden the theorem to the above
generality; his proof was a straightforward adaptation of Bing's.

The easy half of the proof of the theorem is the implication $\Rightarrow$;
one simply writes $h=g_{0}^{-1}g_{1}$ for two successively chosen
homeomorphisms $g_{0},g_{1}$ approximating $f$. The nontrivial, and
significant half of the theorem is the implication $\Leftarrow$. Bing's idea
here was to construct a surjection $p:M\rightarrow M$, with $fp$ close to $f$,
such that the point-inverse sets of $p$ coincide exactly with those of $f$.
Then $g\equiv fp^{-1}:M\rightarrow Q$ defines the desired homeomorphism
approximating $f$. The map $p$ is constructed as a limit of homeomorphisms
$p=\lim_{i\rightarrow\infty}h_{1}h_{2}\ldots h_{i}$ where the $h_{i}^{\prime
}s$ are provided by the Shrinking Criterion, for $\epsilon_{i}$ values which
go to 0 and $i$ goes to $\infty$. The heart of the proof is showing how to
choose each $\epsilon_{i}$, the main point being that it depends on the
composition $h_{1}\ldots h_{i-1}$. Details are given in many places, for
example, \cite[pp. 45,46]{Ch}.

This proof of the implication $\Leftarrow$ becomes quite transparent when
recast as a Baire category argument. For in the Baire space $C(M,Q)$ of maps
from $M$ to $Q$, with the uniform metric topology, let $E$ be the closure of
the set $\{fh\mid h:M\rightarrow M$ is a homeomorphism$\}$. The Bing Shrinking
Criterion amounts to saying that for any $\epsilon>0$, the open subset of
$\epsilon$-maps in $E$, call it $E_{\epsilon}$, is dense in $E$. Hence
$E_{0}\equiv\cap_{\epsilon>0}E_{\epsilon}$ is dense in $E$, since $E$ is a
Baire space. Since $E_{0}$ consists of homeomorphisms, this show that $f\in E$
is approximable by homeomorphisms.

\indent In applying the above Shrinking Theorem to show that various of his
decomposition spaces were manifolds, possibly after stabilizing (\cite{Bi1}%
\cite{Bi2}\cite{Bi3}), Bing used only the \textquotedblleft
if\textquotedblright\ part of the theorem, establishing that his Criterion
held, i.e. establishing that the decomposition was shrinkable. But to show
that his dogbone decomposition space was \textbf{not} a manifold, Bing used
the \textquotedblleft only if\textquotedblright\ part of the theorem, bridging
the gap by showing that if the decomposition space $Q$ were a manifold, then
the quotient map $f:\mathbb{R}^{3}\rightarrow Q$ would have to be approximable
by homeomorphisms \cite[Section 8]{Bi3}. This result has been generalized to
become the \bigskip

\textsc{Cellular Approximation Theorem.}\textbf{ }(extending Moore \cite{Mo}
for $m=2$, see \cite[\S 11]{McA2}; Armentrout \cite{Ar} for $m=3$; Siebenmann
\cite{Si6} for $m\geq5$; compact without-boundary version). \textit{Suppose
$f:M^{m}\rightarrow Q^{m}$ is a map of closed manifolds, $m\neq4$. Then $f$ is
a cell-like surjection} (read \textit{cellular if $m=3$) if and only if $f$ is
approximable by homeomorphisms.\bigskip}

\indent The previously known ``if'' part is quite easy, and holds for all $m$.
The ``only if'' part requires some highly sophisticated geometrical analysis.

\indent The significance of these two theorems taken together is clear:

\begin{quotation}
\textit{The problem of deciding whether a particular cell-like image of a
manifold-without-boundary is itself a manifold is equivalent to deciding
whether the Bing Shrinking Criterion holds in the source } (dimension
$\mathbf{4}$ excepted).\medskip
\end{quotation}

\indent The subject of cell-like decompositions of manifolds developed in the
1950's into a robust and active theory, thanks to Bing's spectacular series of
theorems. By comparison, progress throughout the 1960's was gradual. From 1960
on, the duncehat decomposition problem for $M^{4}\times\mathbb{R}^{1}$
(explained in Part I) was a natural problem to work on; from about 1970 on, it
was the natural problem to work on.

\indent In this paper, only the \textquotedblleft if\textquotedblright\ part
of Bing's Shrinking Theorem is used. But the psychological comfort provided by
Siebenmann's theorem was vital to this work.

\section{Definitions and notation}

Throughout this paper, \textit{all spaces are locally compact separable
metric,} and all manifolds are topological manifolds, possibly with boundary,
unless states otherwise. All definitions are standard.

\indent The (single) \textit{suspension} of a compact metric space $Y$,
denoted $\Sigma Y$, is the quotient of $Y\times\lbrack-1,1]$ gotten by
identifying to distinct points the two subsets $Y\times\pm1$. The
$k$-\textit{fold suspension} $\Sigma^{k}Y$ is $\Sigma(\ldots(\Sigma Y)\ldots)$
($k$ times); if $k=0$, this is understood to be $Y$. Equivalently, $\Sigma
^{k}Y$ can be defined to the join, $Y\ast S^{k-1}$ of $Y$ with the
$(k-1)$-sphere. It follows that $Y\ast S^{k-1}-S^{k-1}$ is homeomorphic to
$Y\times\mathbb{R}^{k}$, and $Y\ast S^{k-1}-Y$ is homeomorphic to
$\mathring{c}Y\times S^{k-1}$, where $\mathring{c}Y$ denotes the \textit{open
cone} on $Y$, that is, $\mathring{c}Y=Y\times\lbrack0,1)/Y\times0$.

A compact metric space is \textit{cell-like} provided that whenever it is
embedded in an ANR (e.g., in the Hilbert cube), it is null-homotopic inside of
any given neighborhood of itself in the ANR (see \cite{Lac1} and \cite{Lac2}
for an elaboration of this and related facts). A map $f:X\rightarrow Y$ of
spaces is \textit{proper} provided the preimage of any compact subset of $Y$
is compact, and $f$ is \textit{cell-like} provided it is onto and each
point-inverse is cell-like. A \textit{cell-like upper-semicontinuous
decomposition} of a space $X$ is nothing more than a proper cell-like
surjection from $X$ onto some (quotient) space $Y$.

\indent Standard notations are $\operatorname*{int}_{X}A,\operatorname*{cl}%
_{X}A$ and $\operatorname*{fr}_{X}A$ for the \textit{interior},
\textit{closure} and \textit{frontier} of a subset $A$ in a space $X$; the
subscript $X$ is omitted whenever it is clear. In this paper careful
distinction is made between the notions of \textit{boundary} and
\textit{frontier}, the former being used only in the manifold sense, e.g.
$\partial M$, and the latter being used only in the point-set sense, e.g.
$\operatorname*{fr}A.$

\indent This paper is concerned largely with constructing homeomorphisms and
ambient isotopies of manifolds. The ambient manifolds can always be assumed to
be piecewise linear (or smooth, if one prefers), and it will be convenient to
use the well-established concepts of piecewise linear topology, e.g., regular
neighborhoods, expansions and collapses, and transversality. As an example,
the phrase \textquotedblleft$J$ \textit{is transverse to} $M\times0$
\textit{in} $M\times\mathbb{R}^{1}$\textquotedblright, where $J$ is a
polyhedron (homeomorphic to $K\times\mathbb{R}^{1}$, in Part I) and $M$ is a
manifold, will mean that there is a bicollar $\alpha:M\times(-1,1)\rightarrow
M\times\mathbb{R}^{1}$ for $M=M\times0$ in $M\times\mathbb{R}^{1}$ (which,
after ambient isotopy fixing $M\times0$, can be assumed to be the standard
bicollar) such that $\alpha|(J\cap M\times(-1,1))$ is a bicollar for $J\cap
M\times0$ in $J$ (everything $PL$).

\part{The double suspension of Mazur's homology $3$--sphere is $S^{5}$}

\indent The goal of this part is to show that the double suspension of a
certain homology $3$-sphere described by Mazur is homeomorphic to $S^{5}$
(this was announced in \cite{Ed2}). The proof is completely descriptive, so
that in particular one can see the wild suspension circle in $S^{5}$
materialize as a limit of tame circles.

\indent By Lefschetz duality, the boundary of any compact contractible
manifold is a homology sphere. Examples of contractible $(n+1)$-manifolds with
nonsimply connected boundary (i.e. genuine homology $n$-sphere boundary), for
$n\geq4,$ were first constructed by Newman \cite{New1}, by taking the
complement in $S^{n+1}$ of the interior of a regular neighborhood of any
acyclic, nonsimply-connected $2$-complex (some applications of this
construction appear in \cite{Wh2}, \cite{CW} and \cite{CZ}). Constructions
were subsequently given for $n=3$ by Poenaru \cite{Po} and Mazur \cite{Maz}.
In this part we are interested in the simplest possible such construction, of
lowest possible dimension, namely $n=3$.

\indent The Poenaru construction goes as follows. Take the product
\[
P^{4}\equiv(S^{3}-\operatorname*{int}N(K_{0}))\times\lbrack-1,1]
\]
of a compact, nontrivial knot complement with an interval, and attach to it a
4-ball $B^{4}$ by gluing a tubular neighborhood $N(K_{1})$ of a nontrivial
knot $K_{1}$ in $\partial B^{4}$ to a tubular neighborhood of a meridian
$\mu\times0$ in $P^{4}$. This is the desired manifold, which can be written
$M^{4}=(P^{4}\cup B^{4})/[N(\mu\times0)=N(K_{1})]$. The advantage of this
construction of Poenaru is the ease with which one can verify that (1) $M^{4}$
is contractible, since the attached 4-ball kills $\pi_{1}(P^{4})$, and (2)
$\partial M^{4}$ is not simply connected, since it is the union of two knot
complements, and the Loop Theorem implies that for any nontrivial knot $K$ in
$S^{3}$, $\pi_{1}(\partial N(K))\rightarrow\pi_{1}(S^{3}-\operatorname*{int}%
N(K))$ is monic. The disadvantage of this construction (for the purposes of
Part I) is that the spine of $M^{4}$, although it can be chosen $2$%
-dimensional, is not especially simple.

\indent Mazur's construction, on the other hand, starts with the simplest
possible spine, i.e. the simplest contractible but noncollapsible polyhedron,
namely the $2$-dimen\-sion\-al duncehat, and thickens it to be a 4-manifold
with nonsimply connected boundary. Recall that the duncehat can quickly be
described as the space gotten by attaching a $2$-cell $B^{2}$ to a circle
$S^{1}$ by sewing $\partial B^{2}$ to $S^{1}$ by a degree one map which goes
twice around $S^{1}$ in one direction, and once around in the other direction.
Since Mazur's homology 3-sphere is the one we work with, we recall further
details and establish some notation.

\indent The duncehat $K^{2}$ and its thickening $M^{4}$ will be described
together, making use of Figure I-1.

\begin{figure}[th]
\centerline{
\includegraphics{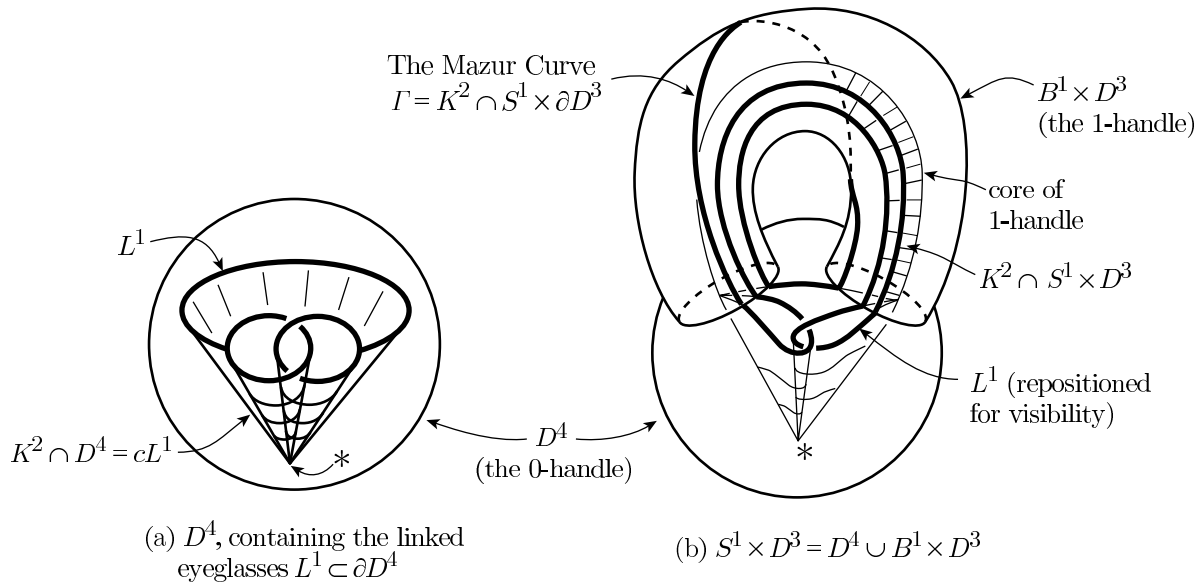}
}\caption{The Mazur 4-manifold $M^{4}=\text{0-handle\ }D^{4}\cup\text{1-handle
}B^{1}\times D^{3}\cup\text{2-handle\ }B^{2}\times D^{2}$, and its duncehat
spine $K^{2}$}%
\end{figure}

It is convenient to express $M^{4}$ as a union of a 0-handle $D^{4}$, a
1-handle $B^{1}\times D^{3}$ and a $2$-handle $B^{2}\times D^{2}$. Let $L^{1}$
be the 1-complex embedded in $\partial D^{4}$ as shown in Figure I-1a.

$L^{1}$ is usually referred to as the \emph{linked eyeglasses}. The
intersection $K^{2}\cap D^{4}$ is to be $cL^{1}$, where the coning is done to
the centerpoint \*of $D^{4}$. The 1-handle $B^{1}\times D^{3}$ is attached to
$D^{4}$ as shown in Figure I-1b, and $K^{2}\cap S^{1}\times\partial D^{3}$ is
defined to be the familiar Mazur curve $\Gamma$ in $S^{1}\times\partial D^{3}%
$. The intersection $K^{2}\cap S^{1}\times D^{3}$ is to be regarded as the
mapping cylinder of an apparent degree one map from $\Gamma$ onto the core
$S^{1}\times0$ of $S^{1}\times D^{3}$. Finally, the $2$-handle $B^{2}\times
D^{2}$ is attached to $S^{1}\times B^{3}$ by identifying $\partial B^{2}\times
D^{2}$ with a neighborhood of $\Gamma$ in $S^{1}\times\partial D^{3}$ in such
a manner that $\partial B^{2}\times0$ is identified to $\Gamma$. Then
$K^{2}\cap B^{2}\times D^{2}$ can be taken to be $B^{2}\times0$. The
\textquotedblleft twisting\textquotedblright\ of the attaching map here (in
the direction perpendicular to $\Gamma$) is not specified, because it is
irrelevant in the upcoming constructions.

The only properties of the pair $(M^{4},K^{2})$ used below are:

\begin{itemize}
\item[(i)] $M^{4}$ is a contractible $PL$ manifold, and $K^{2}$ is a
polyhedral spine of $M^{4}$;

\item[(ii)] $K^{2}-K^{(1)}$ is an open $2$-cell locally flatly embedded in
$\operatorname*{int}M^{4}$, where $K^{(1)}$ denotes the intrinsic 1-skeleton
of $K^{2}$, namely $S^{1}\times0$; and

\item[(iii)] near $\ast\in K^{2}$, the embedding of $K^{2}$ in
$\operatorname*{int}M^{4}$ is as described in Figure I-1.
\end{itemize}

\indent To prove that $\Sigma^{2}(\partial M^{4})\approx S^{5}$, i.e. that
$M^{4}/K^{2}\times\mathbb{R}^{1}$ is a manifold, it suffices to establish that
the Bing Shrinking Criterion holds for the stabilized quotient map $\pi\times
id_{\mathbb{R}^{1}}:M\times\mathbb{R}^{1}\rightarrow M/K\times\mathbb{R}^{1}$.
This amounts to establishing the\bigskip

\textsc{Shrinking Proposition. }\emph{Given any }$\varepsilon>0$\emph{, there
is a homeomorphism }$h:M^{4}\times\mathbb{R}^{1}\rightarrow M^{4}%
\times\mathbb{R}^{1}$\emph{, fixed on }$\partial M^{4}\times R^{1}$\emph{,
such that for each }$t\in R^{1}$\emph{, }

\begin{enumerate}
\item $h\left(  M\times t\right)  \subset M\times\left[  t-\varepsilon
,t+\varepsilon\right]  $\emph{, and}

\item $\operatorname*{diam}\left(  h\left(  K^{2}\right)  \times t\right)
<\varepsilon$\emph{.\bigskip}
\end{enumerate}

\textsc{Technical Notes. }The following brief comments concern the precise
relation of the above statement to the Bing Shrinking Criterion as stated in
the Preliminaries: (a) If one prefers to stay in the world of compacta, one
can replace $\mathbb{R}^{1}$ by $S^{1}$; the reasons that $\mathbb{R}^{1}$ is
used here are (i) tradition, and (ii) the notion of, and notation for,
\textquotedblleft vertical\textquotedblright\ and \textquotedblleft
horizontal\textquotedblright\ motions are more natural in $M\times
\mathbb{R}^{1}$ than in $M\times S^{1}$; (b) tradition has the homeomorphism
$h$ being uniformly continuous, which in fact it is by construction but that
is not really necessary; and (c) in the above statement of the Shrinking
Proposition, the Bing Shrinking Criterion really demands that $\partial M$ be
replaced by $M-\operatorname*{int}N_{\epsilon}(K)$, but this stronger
statement is clearly deducible from the given one by replacing $M$ by a small
regular neighborhood of $K$ in $M$.

\indent The above homeomorphism $h$ will be isotopic to the identity
$rel\,\partial M\times\mathbb{R}^{1}$, by construction. This is useful to keep
in mind when trying to visualize it.

\indent To understand the motivation for the following construction of $h$,
recall the basic principle:

\begin{quotation}
$h$ \textit{is to be thought of as being an arbitrarily close approximation to
a surjection} $p:M\times\mathbb{R}^{1}\rightarrow M\times\mathbb{R}^{1}$,
\textit{whose nontrivial point-inverses are precisely the sets} $\{K\times
t\mid t\in\mathbb{R}^{1}\}$.
\end{quotation}

\noindent This is because in Bing's proof of the implication $\Leftarrow$ of
his theorem (see the Preliminaries), he constructs $p$ as a limit
$p=\lim_{i\rightarrow\infty}h_{1}h_{2}\ldots h_{i}$, where the $h_{i}$'s are
homeomorphisms provided by the Shrinking Proposition, for smaller and smaller
values of $\epsilon$. (Recall that from this map $p$ one gets the desired
homeomorphism $g$ approximating the quotient map $\pi\times id_{\mathbb{R}%
^{1}}:M\times\mathbb{R}^{1}\rightarrow M/K\times\mathbb{R}^{1}$ by defining
$g=(\pi\times id_{\mathbb{R}^{1}})p^{-1}:M\times\mathbb{R}^{1}\overset
{\approx}{\rightarrow}M/K\times\mathbb{R}^{1}$.)

\indent Now, the image $p(K\times\mathbb{R}^{1})$ $(\approx\mathbb{R}^{1})$ in
$M\times\mathbb{R}^{1}$ must have the following remarkable pathological
property (recognized by Glaser in \cite[p. 16, last two paragraphs]{Gl2}):

\begin{quotation}
\textit{for each} $t\in\mathbb{R}^{1}$, \textit{the intersection}
$p(K\times\mathbb{R}^{1})\cap M^{4}\times t$ \textit{must be} $0$%
-\textit{dimensional, yet wild enough so that the inclusion $\partial
M^{4}\times t\hookrightarrow M\times t-p(K\times\mathbb{R}^{1})$ induces a
$\pi_{1}$-monomorphism.}
\end{quotation}

\noindent The second assertion follows because the inclusion
\[
\partial M^{4}\times t\hookrightarrow M^{4}\times\mathbb{R}^{1}-p(K\times
\mathbb{R}^{1})
\]
must induce a $\pi_{1}$-monomorphism, which is clear from the properties of
$p$. The first assertion follows because if $p(K\times\mathbb{R}^{1})\cap
M^{4}\times t$ were not 0-dimension, then there would exist some interval
$(a,b)\subset\mathbb{R}^{1}$ such that $p(K\times(a,b))\subset M^{4}\times t$,
which would then guarantee that the interval $p(K\times(a,b))$ would be
locally flatly embedded in $M\times\mathbb{R}^{1}$ (say by the Klee trick, as
in \cite[Part II]{CW}), and hence $\partial M$ would be simply connected (e.g.
by general positioning off of $p(K\times\mathbb{R}^{1})$ in $M\times
\mathbb{R}^{1}$ a $2$-disc initially mapped into $p(M\times c)$, $c=$ midpoint
of $(a,b)$, and then pushing this $2$-disc out to $\partial M\times
\mathbb{R}^{1})$. Consequently, if the homeomorphism $h$ of the Shrinking
Proposition is to approximate $p$, then given $\epsilon>0$ one must be able to
construct a homeomorphism $h$ having the following property:

\begin{quotation}
\textit{for each $t\in\mathbb{R}^{1}$, each component of the intersection}
$h(K\times\mathbb{R}^{1})\cap M\times t$ \textit{has diameter} $<\epsilon$
(and yet necessarily the inclusion $\partial M\times t\hookrightarrow M\times
t-h(K\times\mathbb{R}^{1})$ must induce a $\pi_{1}$-monomorphism).
\end{quotation}

\noindent This is what motivates the construction below.

\indent The Basic Lemma below (or more accurately, the first 3 steps of it)
shows how to construct such a homeomorphism $h$ having this intersection
property with respect to the single level $M \times0$. This Lemma is the heart
of the entire proof; everything after it amounts to tidying up.

\indent The first, second and fourth steps of the Basic Lemma involve
splitting operations, which are called either meiosis or mitosis, depending on
the context. They are the freshest ingredients of the proof.

\indent The Lemma requires a definition. Given a homeomorphism $h:M\times
\mathbb{R}^{1}\rightarrow M\times\mathbb{R}^{1}$, then a component $C$ of
$h(K\times\mathbb{R}^{1})\cap M\times0$ is \textit{source-isolated} if
$h^{-1}(C)$ is isolated in $K\times\mathbb{R}^{1}$, in the sense that there is
an interval $(a,b)\subset\mathbb{R}^{1}$ such that $K\times(a,b)\cap
h^{-1}(M\times0)=h^{-1}(C)$. In the statement below, $C$ most often will be a
2-sphere lying in $h(K-K^{(1)})\times\mathbb{R}^{1}$, in which case
source-isolation will imply that $h^{-1}(C)$ bounds a 3-cell in $K\times
\mathbb{R}^{1}$ whose interior misses $h^{-1}(M\times0)$.

\indent One should think of the following lemma as starting at $i=0$, with
$h_{0}=$ identity, and producing successively better approximations to the
desired homeomorphism $h$ of the Shrinking Proposition.\bigskip

\textsc{Basic Lemma}$_{i}$ (to be read separately and successively for
$i=1,2,3,4$). \emph{Given }$\delta>0$\emph{, there exists a homeomorphism
}$h_{i}:M\times\mathbb{R}^{1}\rightarrow M\times\mathbb{R}^{1}$\emph{, with
compact support in }$\operatorname*{int}M\times\mathbb{R}^{1}$\emph{, such
that }$h_{i}(K\times\mathbb{R}^{1})$\emph{ is transverse to }$M\times0$\emph{
and the components of }$h_{i}(K\times\mathbb{R}^{1})\cap M\times0$\emph{ are
as follows.}

\begin{itemize}
\item[($i=1$)] \emph{There are two components. One is a pierced duncehat
(described below) having diameter }$<\delta$\emph{, and the other is a
source-isolated }$2$\emph{-sphere }$\Sigma^{2}\times0$\emph{, where }%
$\Sigma^{2}\subset M^{4}$\emph{.}

\emph{Let }$\Sigma^{2}\times D^{2}$\emph{ be a product neighborhood of
}$\Sigma^{2}$\emph{ in }$M^{4}$\emph{.}

\item[($i=2$)] \emph{There are }$1+2^{p}$\emph{ components }$(p=p(\delta
))$\emph{. One is the }$\delta$\emph{-small pierced duncehat from Step 1, and
the remaining }$2^{p}$\emph{ are source-isolated 2-spheres lying in }%
$\Sigma^{2}\times D^{2}\times0$\emph{, such that this collection of }$2^{p}%
$\emph{ 2-spheres in }$\Sigma^{2}\times D^{2}\times0$\emph{ is equivalent to
the }$p^{th}$\emph{ stage in the spun Bing collection of 2-spheres in }%
$S^{2}\times D^{2}$\emph{ (described below).}

\item[($i=3$)] \emph{There are }$1+2^{p}$\emph{ components (same }%
$p=p(\delta))$\emph{. One is the }$\delta$\emph{-small pierced duncehat, and
the remaining }$2^{p}$\emph{ are source-isolated }$\delta$\emph{-small
2-spheres in }$\Sigma^{2}\times D^{2}\times0$\emph{.}

\item[($i=4$)] \emph{There are }$1+2^{p+1}$\emph{ components (same
}$p=p(\delta))$\emph{. One is the }$\delta$\emph{-small pierced duncehat, and
the remaining }$2^{p+1}$\emph{ are }$\delta$\emph{-small embedded duncehats
lying in }$\Sigma^{2}\times D^{2}\times0$\emph{, each being the image under
}$h_{4}$\emph{ of some entire level }$K\times t_{q},1\leq q\leq2^{p+1}%
$\emph{.\bigskip}
\end{itemize}

\textsc{Note. }In using the above Lemma later, all we will really care about
is the intersection of $h_{i}(K\times\mathbb{R}^{1})$ with the 2-handle
$B^{2}\times D^{2}\times0\subset M^{4}\times0$, not with the entire level
$M^{4}\times0$. But we have stated the Lemma in the above elaborated form to
make the overall process a little clearer.

\begin{proof}
It is best to gain first a qualitative understanding of the components of
$h_{i}(K\times\mathbb{R}^{1})\cap M\times0$, for $i=1,\ldots,4$, without
paying attention to their size. For this purpose, it is easiest to work in the
source, understanding there the pre-images $K\times\mathbb{R}^{1}\cap
h_{i}^{-1}(M\times0)$. After this one can come to grips with the images
$h_{i}(K\times\mathbb{R}^{1})\cap M\times0$, and the size of their components.
Following this advice, we present the construction in two rounds, in the first
round just working in the source, describing there the sets $K\times
\mathbb{R}^{1}\cap h_{i}^{-1}(M\times0),1\leq i\leq4$, and their components.
In the second round we pay attention to how these components are embedded in
$M\times0$, and their size.

In the first round of the construction, each set $h_{i}^{-1}(M\times0),1\leq
i\leq4$, will be described as the frontier in $M\times\mathbb{R}^{1}$ of an
arbitrarily small (relative) regular neighborhood in $M\times\mathbb{R}^{1}$
of a subpolyhedron $A_{i}$ of $M\times\mathbb{R}^{1}$, where%
\[
M\times(-\infty,0]=A_{0}\,\includegraphics{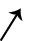}\,A_{1}%
\,\includegraphics{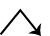}\,\ldots\,\includegraphics{x-udarrow}\,A_{2}%
=A_{3}\,\includegraphics{x-udarrow}\,A_{4}.
\]
Here the arrows indicate expansions and collapses. Each $A_{i}$ will contain
$A_{0}$, and $\operatorname*{cl}(A_{i}-A_{0})$ will be a compact 2-dimensional
polyhedron in $\operatorname*{int}M\times\mathbb{R}^{1}$. The regular
neighborhoods will be relative to $\partial M\times(-\infty,0)]$, and hence
the intersection of their frontiers with $\partial M\times\mathbb{R}^{1}$ will
be $\partial M\times0$. Each regular neighborhood is to be chosen so that its
restriction to $K\times\mathbb{R}^{1}$ is a regular neighborhood in
$K\times\mathbb{R}^{1}$ of the subpolyhedron $B_{i}\equiv A_{i}\cap
K\times\mathbb{R}^{1}$. The successive sets $\{h_{i}^{-1}(M\times
(-\infty,0])\},1\leq i\leq4$, can be thought of as gotten by applying the
uniqueness of regular neighborhoods principle, with each homeomorphism $h_{i}$
being the end of an ambient isotopy. (This regular neighborhood description,
from my Cambridge notes, is not quite the manner in which I originally
perceived the $h_{i}$'s, but it is the quickest way to describe them.)

The $(A_{i},B_{i})$'s are as follows (see Figure I-3).

\textbf{Step i = 0.} $A_{0}=M\times(-\infty,0];B_{0}=K\times(-\infty,0]$.

\textbf{Step i = 1} (first round). Let $a\in L^{1}$ be the midpoint of the
interval part of the linked eyeglasses $L^{1}$ (see Figure I-2).

\begin{figure}[th]
\centerline{
\includegraphics{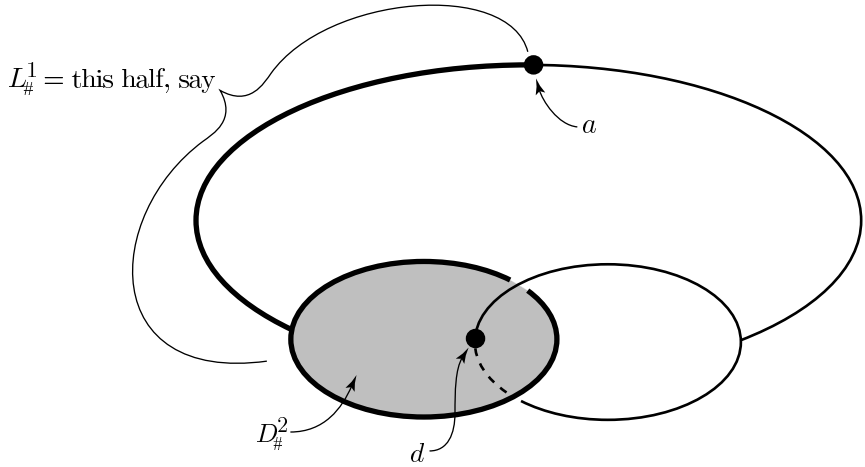}
}\caption{$L^{1}=$ linked eyeglasses $\subset\partial D^{4}$}%
\end{figure}

Let $L_{\#}^{1}$ be the closure of one of the components of $L^{1}-a$ (either
one); hence $L_{\#}^{1}$ is a circle with a feeler attached, with the free end
of the feeler being the point $a$. Let $D_{\#}^{2}$ be the natural 2-disc that
the circle of $L_{\#}^{1}$ bounds in $\partial D^{4}$, so that $D_{\#}^{2}\cap
L^{1}=\partial D_{\#}^{2}\cup d$, where $d$ is an interior point of
$D_{\#}^{2}$. Writing the 0-handle $D^{4}$ of $M^{4}$ as $D^{4}=\partial
D^{4}\times\lbrack0,1)\cup\ast$, then define
\[
A_{1}=A_{0}\cup(L_{\#}^{1}\cup D_{\#}^{2})\times\frac{1}{2}\times1\cup
a\times\frac{1}{2}\times\lbrack0,1]\subset A_{0}\cup D^{4}\times\mathbb{R}%
^{1}\subset M\times\mathbb{R}^{1}%
\]
(see Figure I-3$_{1}$).

\begin{figure}[th]
\renewcommand{\thefigure}{\Roman{part}-3$_1$} \centerline{
\includegraphics{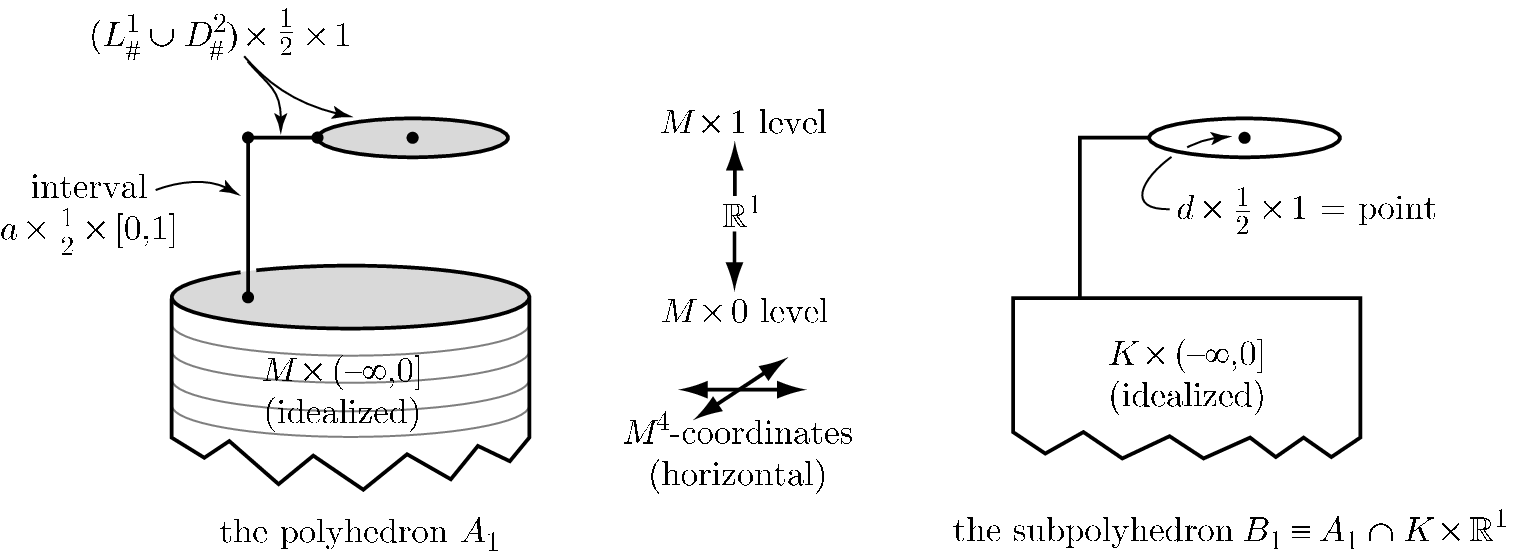}
}\end{figure}

\begin{figure}[th]
\renewcommand{\thefigure}{\Roman{part}-3$_2$, \textnormal{start}}
\centerline{
\includegraphics{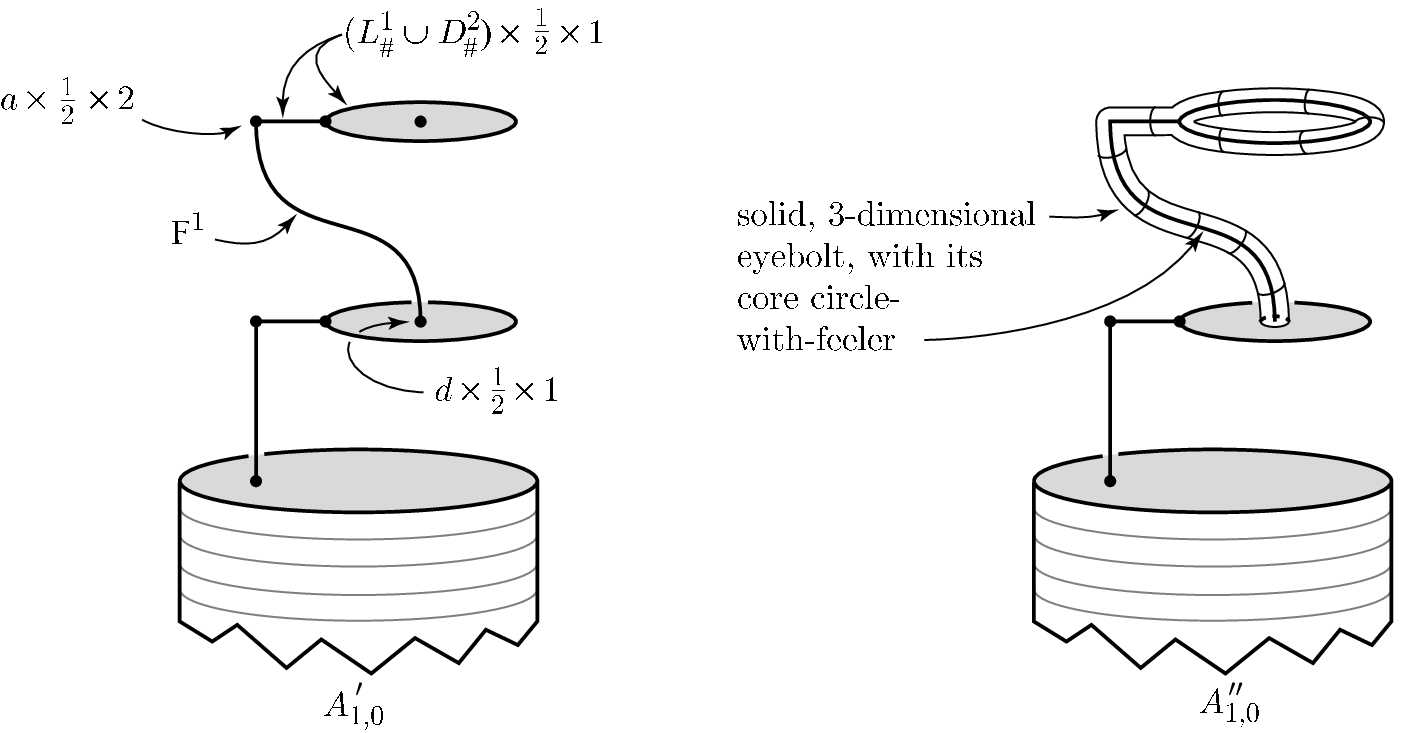}
}\end{figure}

\begin{figure}[th]
\renewcommand{\thefigure}{\Roman{part}-3$_2$, \textnormal{cont.}}
\centerline{
\includegraphics{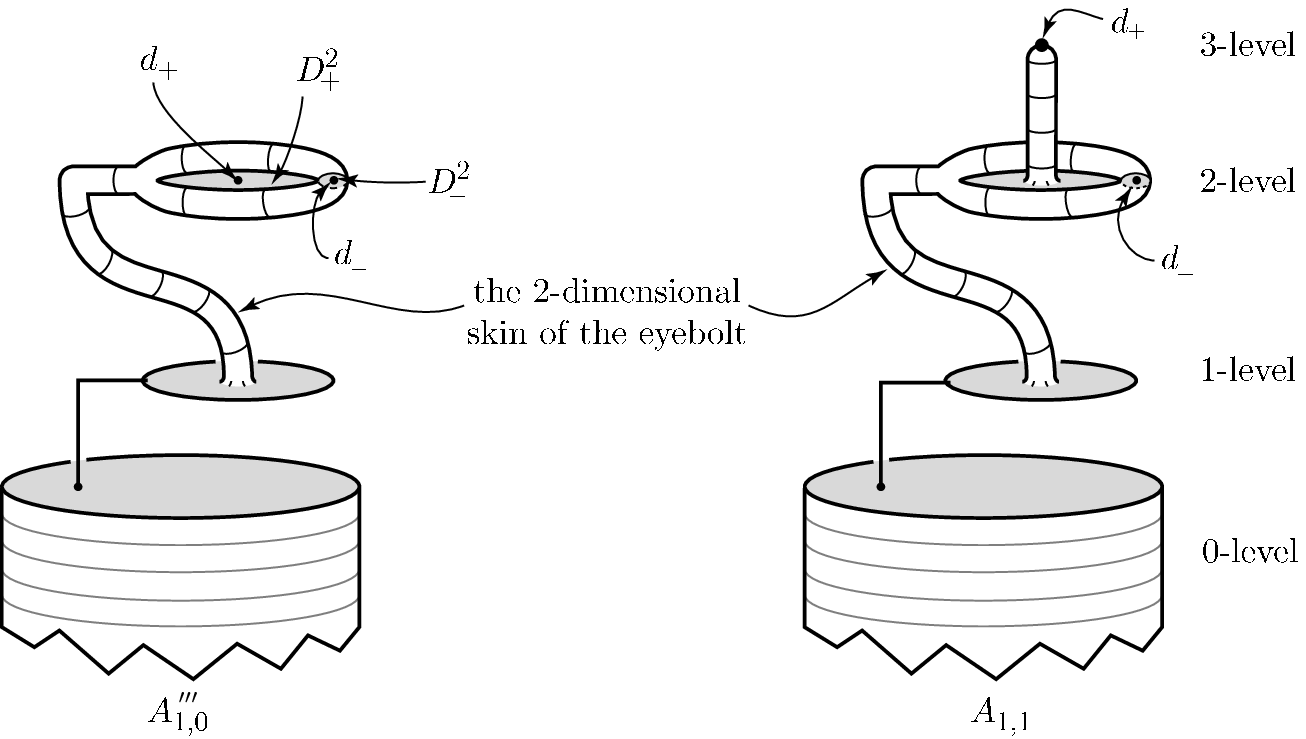}
}\end{figure}

\begin{figure}[th]
\setcounter{figure}{2} \renewcommand{\thefigure}{\Roman{part}-\arabic{figure}}
\centerline{
\includegraphics{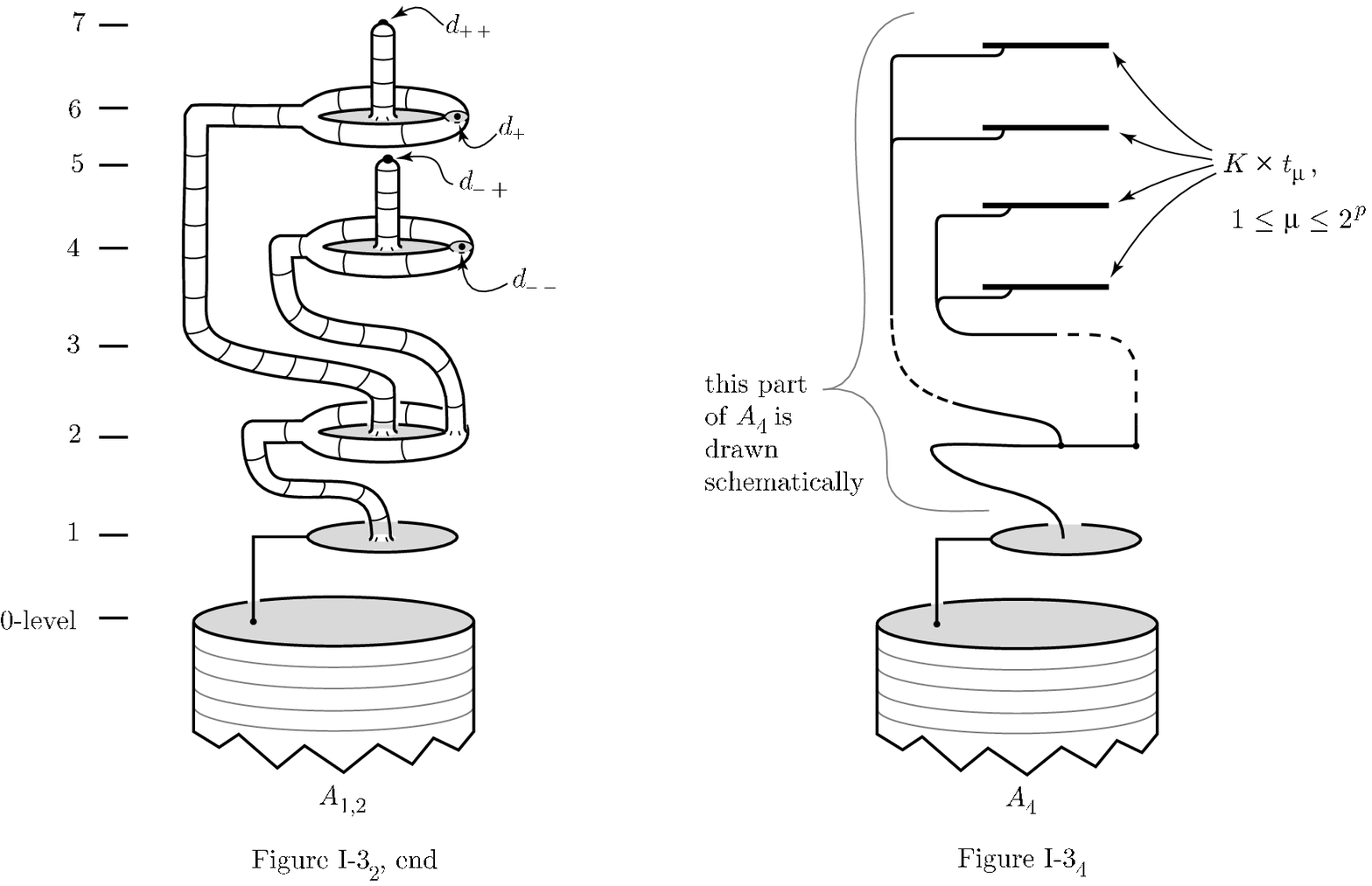}
}\end{figure}

In other words, $A_{1}-A_{0}$ lies in $\partial D^{4}\times\frac{1}{2}%
\times\lbrack0,1]$, where [0,1] is a sub-interval of the vertical coordinate
$\mathbb{R}^{1}$ of $M\times\mathbb{R}^{1}$, and $A_{1}-A_{0}$ consists of a
horizontal 2-disc-with-feeler lying in the $M\times1$ level, together with the
vertical interval joining the free end of the feeler to the $M\times0$ level.
Let $A_{\delta}$ denote $A_{1}$ minus the interior of the 2-disc $D_{\#}%
^{2}\times\frac{1}{2}\times1$, and let $B_{\delta}=A_{\delta}\cap
K\times\mathbb{R}^{1}$, which is $K\times(-\infty,0]$ with a
circle-with-feeler attached. By the construction, $B_{1}=B_{\delta}\cup
d\times\frac{1}{2}\times1$, which has two components, one of them a point.

The set $K\times\mathbb{R}^{1}\cap h_{1}^{-1}(M\times0)$, being the frontier
in $K\times\mathbb{R}^{1}$ of a regular neighborhood of $B_{1}$ in
$K\times\mathbb{R}^{1}$, will consist of two components, one of them the
2-sphere boundary of a 3-ball, and the other a \textit{pierced duncehat} (see
Figure I-4). As an abstract set, a pierced duncehat $K_{\#}^{2}$ can be
described as the union of a duncehat-with-hole and a 2-torus, the boundary of
the hole being identified with an essential curve on the 2-torus. In symbols,
\[
K_{\#}^{2}=(K^{2}-\operatorname*{int}C^{2})\cup_{\partial C^{2}=S^{1}%
\times\text{point}}S^{1}\times S^{1},
\]
where $C^{2}$ is a $2$-cell in $K^{2}-K^{(1)}$. This will be taken up again in
the second round; for the moment, our main concern is to define the sets
$\{h_{i}^{-1}(M\times(-\infty,0])\}$ and (therefore) $\{K\times\mathbb{R}%
^{1}\cap h_{i}^{-1}(M\times(-\infty,0])\}$, leaving to be specified the actual
behavior of the $h_{i}$'s on these sets.

\begin{figure}[ptb]
\centerline{
\begin{tabular}{c}
\includegraphics{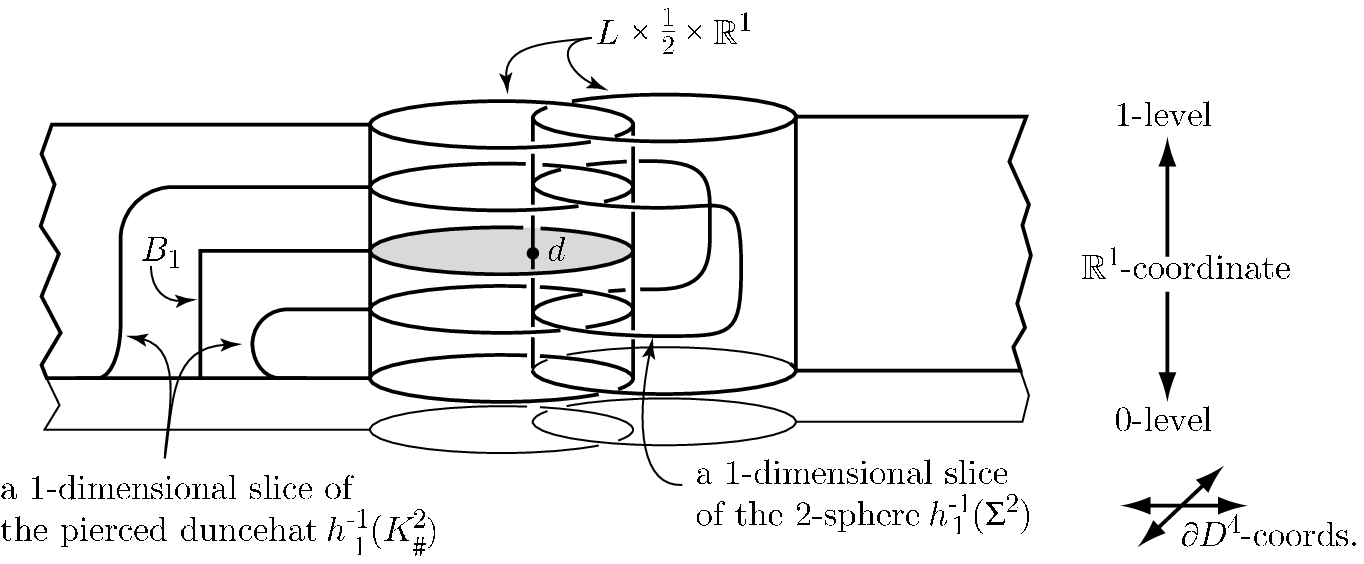}\\
Part (a).  A slice of $K^2\times \mathbb{R}^1 \cap h_1^{-1}(M^4\times 0)$.  The ambient space\\
of the picture is the 4-dimensional slice $\partial D \times \frac12 \times \mathbb{R}^1$ of $M^4\times\mathbb{R}^1$.\\
\bigskip{}\\
\includegraphics{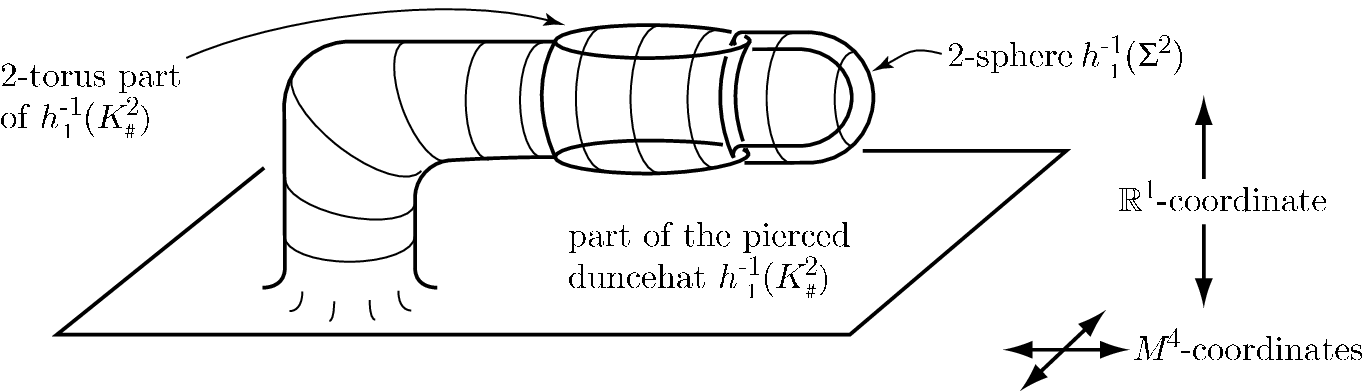}\\
Part (b).  Another view of $K^2\times \mathbb{R}^1 \cap h_1^{-1}(M^4\times 0)$, which consists\\
of the pierced duncehat $h^{-1}_1(K^2_\#)$ and the 2-sphere $h^{-1}_1(\Sigma^2)$.\\
The view in Part (a) is a slice of this.\\
\bigskip{}\\
\includegraphics{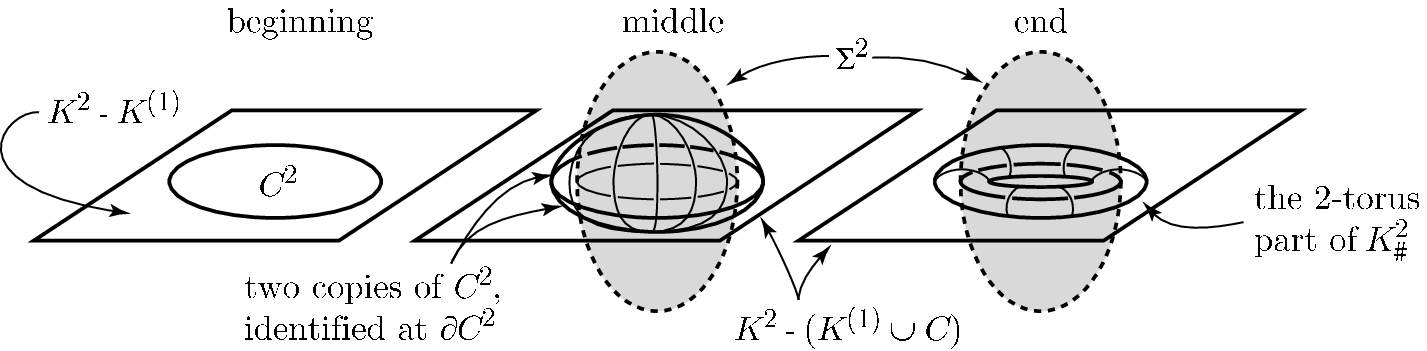}\\
Part (c).  Three stages in the growth of the pierced duncehat $K^2_\#$ and the \\
linking 2-sphere $\Sigma^2$, taking place in $M^4\times 0$.  The 2-sphere $\Sigma^2$ has dotted \\
lines to indicate that the hemispheres of $\Sigma^2$ lie in the 4th dimension; \\
only the equator of $\Sigma^2$ lies in this 3-dimensional slice of $M^4\times 0$.
\end{tabular}
}\caption{Some aspects of Step $i=1$, including the pierced duncehat}%
\end{figure}

\textbf{Step i = 2} (first round). This is the interesting step. The
polyhedron $A_{2}$ will be described in stages, by describing a finite
sequence of polyhedra%
\[
A_{1}=A_{1,0}\includegraphics{x-udarrow}A_{1,1}\includegraphics{x-udarrow}\ldots
\includegraphics{x-udarrow}A_{1,p-1}\includegraphics{x-udarrow}A_{1,p}=A_{2}%
\]
where $p$ is to be specified. The first expansion-collapse is typical of them
all; it is described using two expansions and a collapse, plus a final
repositioning:
\[
A_{1,0}\includegraphics{x-uarrow}A_{1,0}^{^{\prime}}%
\includegraphics{x-uarrow}A_{1,0}^{^{\prime\prime}}%
\includegraphics{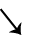}A_{1,0}^{^{\prime\prime\prime}}\leadsto A_{1,1}%
\]
(see Figure I-$3_{2}$).

The polyhedron $A_{1,0}^{\prime}$ is gotten from $A_{1,0}=A_{1}$ by starting
at the point $d\times\frac{1}{2}\times1$ in $A_{1}$, and sending out a feeler
(= interval) $F^{1}$ upwards from this point to the point $a\times\frac{1}%
{2}\times2$, always staying in $(K-K^{(1)})\times\lbrack1,2]$, and then
adjoining the disc-with-feeler $(L_{\#}^{1}\cup D_{\#}^{2})\times\frac{1}%
{2}\times2$ to $A_{1,0}\cup F^{1}$. Note that the feeler $F^{1}$ may have to
travel a long distance to join the two points, but nevertheless it is
possible, since the points $a$ and $d$ are joinable by an arc in the open
2-cell $K^{2}-K^{(1)}$. (Aside: This entire step could be done without any
vertical motion, working always in the $M\times1$ level to change $A_{1}$ to
$A_{2}$, postponing to Step 3 or 4 the matter of source-isolation. This is the
point of view adopted in the more general program in Part II. But for Part I
the present description seems a bit easier.)

Next, $A_{1,0}^{^{\prime\prime}}$ is gotten by thickening part of
$A_{1,0}^{\prime}$, adding to it a small 3-dimensional \textquotedblleft
tubular neighborhood\textquotedblright\ of $F^{1}\cup L_{\#}^{1}\times\frac
{1}{2}\times2$, resembling a 3-dimensional eyebolt where \textquotedblleft
tubular neighborhood\textquotedblright\ of $F^{1}$ means the restriction to
$F^{1}$ of a 2-disc normal bundle of $(K^{2}-K^{(1)})\times\mathbb{R}^{1}$ in
$M^{4}\times\mathbb{R}^{1}$, and where \textquotedblleft tubular
neighborhood\textquotedblright\ of $L_{\#}^{1}\times\frac{1}{2}\times2$ means
a genuine regular neighborhood of it (rel $a\times\frac{1}{2}\times2$) in
$\partial D^{4}\times\frac{1}{2}\times2$. We assume these thickenings match up
nicely with each other, and with what has already been defined, as suggested
by Figure I-$3_{2}.$

To get $A_{1,0}^{^{\prime\prime\prime}}$ from $A_{1,0}^{^{\prime\prime}}$,
most of the newly added 3-dimensional eyebolt is collapsed away, starting at
the $(d\times\frac{1}{2}\times1)$-end and collapsing upwards, so that the only
material which is left behind is the 2-dimensional outer skin of the eyebolt,
plus a spanning 2-disc $D_{-}^{2}$ in the handle at the end of the eyebolt.
Also left undisturbed is the other 2-disc $D_{+}^{2}$ spanning the eye of the
eyebolt. Thus, $A_{1,0}^{^{\prime\prime\prime}}$ is homeomorphic to
\[
A_{\delta}\cup\left(  2-disc\text{ with handle}\right)  \cup D_{-}^{2}\cup
D_{+}^{2}\text{,}%
\]
where the boundaries of $D_{-}^{2}$ and $D_{+}^{2}$ are identified with the
circles in the figure-eight spine of the disc-with-handle (which is a
punctured 2-tours), and where the boundary of the disc-with-handle is
identified with the circle in $A_{\delta}$. In the collapsing of
$A_{1,0}^{\prime\prime}$ to $A_{1,0}^{\prime\prime\prime}$, some of the
intersection of $A_{1,0}^{\prime\prime}$ with $K\times\mathbb{R}^{1}$ is
collapsed away. Originally it consisted of $B_{\delta}$, plus a circle with a
feeler attached, plus a point $d_{+}\in\operatorname*{int}D_{+}^{2}$. When
done with the collapse the circle-with-feeler has been reduced to just a
single point on the circle, namely the center $d_{-}$ of $D_{-}^{2}$, while
$B_{\delta}$ and the point $d_{+}$ remain undisturbed.

The final repositioning of $A_{1,0}^{\prime\prime\prime}$, done to achieve the
source-isolation property of the Lemma, is to move one of the two intersection
points, say $d_{+}$, upwards to say the $M\times3$ level, keeping fixed the
complement of a small neighborhood of this point in $A_{1,0}^{\prime
\prime\prime}$. Let $A_{1,1}$ be this repositioned copy of $A_{1,0}%
^{\prime\prime\prime}$. By construction,
\[
B_{1,1}=(B_{1,0}-d\times\frac{1}{2}\times1)\cup\{d_{-},d_{+}\}=B_{\delta}%
\cup\{d_{-},d_{+}\}\text{,}%
\]
where one of the new points is at the 2-level, and the other is at the 3-level.

Qualitatively, $A_{1,1}$ is gotten from $A_{1,0}$ in a single operation,
simply by \textit{puckering} the $2$-disc $D_{\#}^{2}\times\frac{1}{2}%
\times1\subset A_{1,0}$, where a \textit{puckered} 2-\textit{disc} is a
2-disc-with-handle with two spanning $2$-discs attached to make it
contractible, as remarked in the description of $A_{1,0}^{\prime\prime\prime}$
above. Note that a puckered 2-disc really is obtained by puckering a genuine
$2$-disc; see Figure I-6a.

To get $A_{1,2}$ from $A_{1,1}$, qualitatively one simply puckers each of the
two spanning $2$-discs $D_{-}^{2}$ and $D_{+}^{2}$ in $A_{1,1}$.
Quantitatively, one does two independent, non-interfering expansion-collapse
processes, each process being a copy of the one used above to go from
$A_{1,0}$ to $A_{1,1}$, this time starting at the two new points $d_{-}$ and
$d_{+}$ of $B_{1,1}$. The initial feelers sent out should go up to say the
4-level (from the $2$-level) and to the 6-level (from the 3-level), and then
the bulk of the activity takes place in those levels, except that the final
repositioning operation makes the four new points of $B_{1,2}=B_{\delta}%
\cup\{d_{--},d_{-+},d_{+-},d_{++}\}$ lie in say the $4,5,6$ and $7$-levels.

Now the pattern is established. The set $K\times\mathbb{R}^{1}\cap
h_{i,j}^{-1}(M\times0)$, being the frontier in $K\times\mathbb{R}^{1}$ of a
regular neighborhood in $K\times\mathbb{R}^{1}$ of $B_{1,j}=B_{\delta}%
\cup2^{j}$ points, will consist of a pierced duncehat and $2^{j}$ $2$-spheres,
and these components are isolated in $K\times\mathbb{R}^{1}$, in that they
have no overlap in the $\mathbb{R}^{1}$-coordinate (i.e., their $\mathbb{R}%
^{1}$-projections are disjoint). On the other hand, the images under $h_{1,j}$
of these components are linked in $M^{4}\times0$ in a very interesting manner,
as will be explained in the second round.

\textbf{Step i = 3} (first round). $A_{3}=A_{2}$. Hence $h_{3}$ is
qualitatively the same as $h_{2}$. The quantitative difference between $h_{3}$
and $h_{2}$ will be explained in the second round.

\textbf{Step i = 4} (first round). See Figures I-3$_{4}$ and I-7. Given
$A_{2}$ $(=A_{3})$, let $\{t_{\mu}\mid1\leq\mu\leq2^{p}\}$ denote the
$\mathbb{R}^{1}$-levels of the $2^{p}$ points of $B_{2}$ (in the above
construction, these levels were $t_{\mu}=2^{p}+\mu-1)$. Then define
$A_{4}=A_{2}\cup\bigcup\{K^{2}\times t_{\mu}\mid1\leq\mu\leq2^{p}\}$. The
expansion-collapse $A_{2}\includegraphics{x-udarrow}A_{4}$ follows from the fact
that $K\times I$ is collapsible. Note that $B_{4}=B_{2}\cup\bigcup
\{K^{2}\times t_{\mu}\mid1\leq\mu\leq2^{p}\}=B_{\delta}\cup\bigcup
\{K^{2}\times t_{\mu}\mid1\leq\mu\leq2^{p}\}$. It is clear that the
intersection $K^{2}\times\mathbb{R}^{1}\cap h_{4}^{-1}(M^{4}\times0)$ consists
of $1+2^{p+1}$ components, one being the pierced duncehat, and the others
being say the $2^{p+1}$ duncehats $\{K^{2}\times(t_{\mu}\pm\lambda)\mid
1\leq\mu\leq2^{p}\}$, for some small $\lambda>0$.

Having established the qualitative definitions of the $h_{i}$'s, we seek now
to gain a better understanding of their behavior, paying close attention to
the target intersections $h_{i}(K^{2}\times\mathbb{R}^{1})\cap M^{4}\times0$.

\textbf{Step i = 1} (second round). Analysis of the two-component intersection
$h_{1}(K\times\mathbb{R}^{1})\cap M\times0$ reveals that it can be described
in the following way (this amplifies the earlier remarks; see Figure I-4c).

Starting with $K^{2}$ in $M^{4}$, remove a small $2$-cell $C^{2}$ from
$K^{2}-K^{(1)}$, and replace it with a $2$-torus to produce the pierced
duncehat component $K_{\#}^{2}$. At the same time, add to the picture a
disjoint $2$-sphere $\Sigma^{2}$ lying in $M^{4}$, so that in $M^{4}$ it links
the $2$-torus of $K_{\#}^{2}$, as suggested by Figure I-4c. Note that the
inclusion $\partial M^{4}\hookrightarrow M^{4}-(K_{\#}^{2}\cup\Sigma^{2})$ is
monic on the fundamental group, as it must be (but this need not be verified).

It is interesting to actually watch the original duncehat
\[
K^{2}\times0=K^{2}\times\mathbb{R}^{1}\cap M^{4}\times0
\]
in $M^{4}\times0$ transform into these two components $K_{\#}^{2}\times0$ and
$\Sigma^{2}\times0$, by following the ambient isotopy $h_{t}$ of
$h_{0}=identity$ to $h_{1}$, focusing on the intersection $h_{t}%
(K\times\mathbb{R}^{1})\cap M\times0$ (or more easily, on $K\times
\mathbb{R}^{1}\cap h_{t}^{-1}(M\times0)$, using the earlier description). The
changes in $K^{2}\times0$ take place in an arbitrarily small neighborhood of
$C^{2}\times0$ in $M^{4}\times0$, and are symmetric under rotation of this
4-dimensional neighborhood about the fixed $2$-plane transverse to
$C^{2}\times0$ at its center (the \textquotedblleft spinning\textquotedblright%
\ point of view, explained in the upcoming $i=2$ case, is also useful here).
The changes are (Figure I-4c): first $C^{2}\times0$ divides, or splits, into
two parallel copies of itself, with boundaries remaining joined together at
$\partial C^{2}\times0$; then the $2$-sphere $\Sigma^{2}\times0$ materializes,
growing from say the original center point of $C^{2}\times0$, until its
diameter is almost that of $C^{2}\times0$; finally, the two parallel
boundary-identified copies of $C^{2}\times0$ join together near their centers,
becoming pierced there, so their union forms the $2$-torus which links the
$2$-sphere in $M\times0$. This dividing operation can be regarded as
\textit{duncehat meiosis}.

The pierced duncehat $K_{\#}\times0$ can now be made arbitrarily small in
$M\times0$, by level-preserving ambient isotopy of $M\times\mathbb{R}^{1}$,
because it can be pushed arbitrarily close to the circle $K^{(1)}\times0$,
which itself can be made small homotopically hence isotopically, using the
contractibility of $K^{2}$. This motion will stretch $\Sigma^{2}\times0$
large, but that is allowed, for Steps 2 and 3 will take care of that. We
henceforth will assume these motions have been incorporated into the
homeomorphism $h_{1}$.

\textbf{Step i = 2} (second round). In this step the $2$-sphere component of
intersection $\Sigma^{2}\times0$ from Step 1 is replaced by $2^{p}$ new
$2$-sphere components of intersection, which are embedded in a tubular
neighborhood $\Sigma^{2}\times D^{2}\times0$ of $\Sigma^{2}\times0$ in
$M^{4}\times0$. The goal here is to describe how these new $2$-spheres are
embedded. The pierced duncehat component $K_{\#}^{2}\times0$ from Step 1 is
left untouched (for the remainder of the proof, in fact).

The model situation comes from Bing's foundational 1952 paper \cite{Bi1}.
There Bing described a certain nest of 3-dimensional solid tori, which he used
to define a decomposition of $\mathbb{R}^{3}$. They can be described as
follows (see Figure I-5).

\begin{figure}[th]
\centerline{
\includegraphics{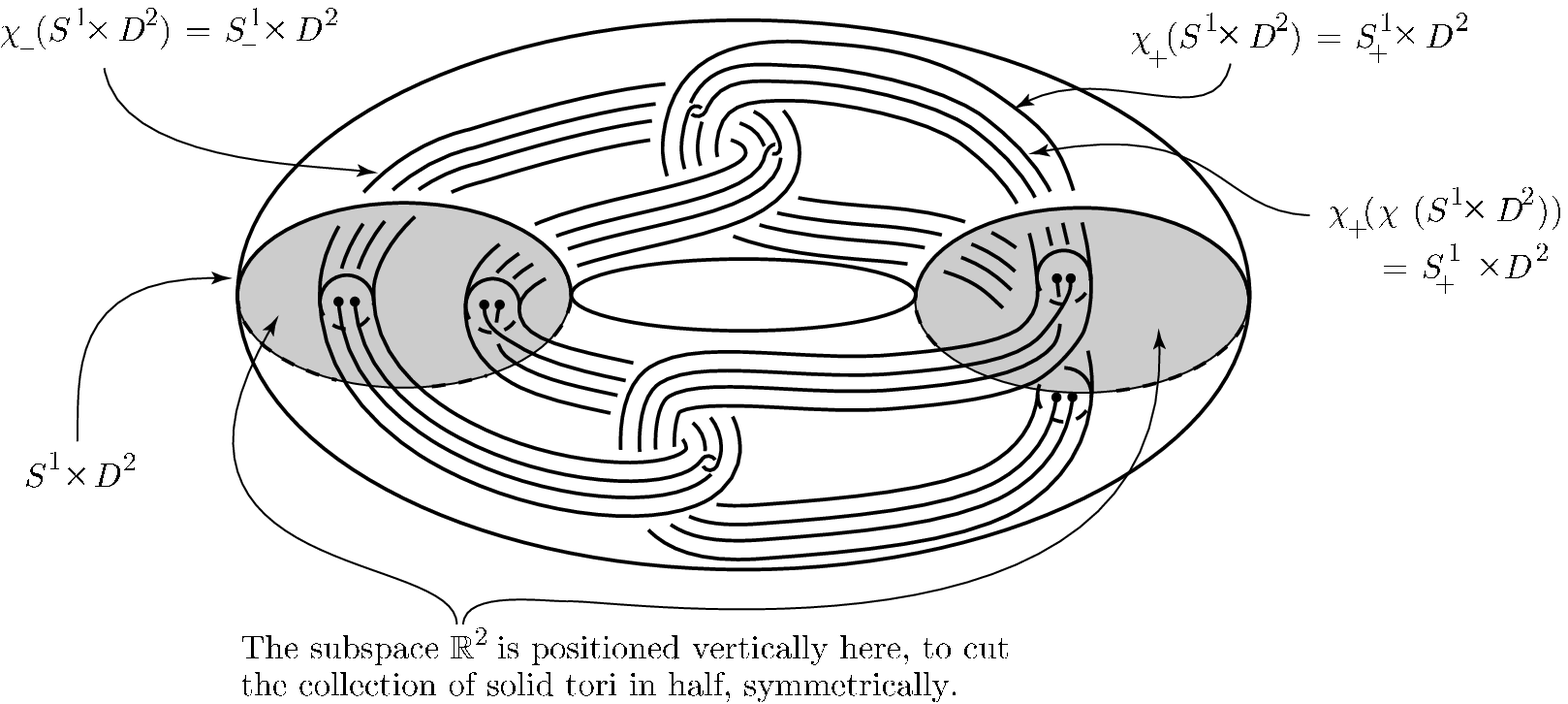}
}\caption{The Bing collection of thickened 1-spheres (i.e., solid tori) in
$\mathbb{R}^{3}$}%
\end{figure}

Starting with $S^{1}\times D^{2}$, let $\chi_{-},\chi_{+}:S^{1}\times
D^{2}\rightarrow S^{1}\times\operatorname*{int}D^{2}$ be two disjoint
embeddings, with images denoted $S_{-}^{1}\times D^{2}$ and $S_{+}^{1}\times
D^{2}$, such that each image is by itself trivially embedded in $S^{1}\times
D^{2}$, and yet the two images are linked in $S^{1}\times D^{2}$ as shown in
Figure I-5. These embeddings can be iterated, to produce for any $p>0$ a
collection of $2^{p}$ solid tori $\{S_{\mu}^{1}\times D^{2}\mid\mu
\in\{-,+\}^{p}\}$ in $S^{1}\times D^{2}$, where%
\[
S_{\mu}^{1}\times D^{2}=\chi_{\mu(1)}(\chi_{\mu(2)}(\ldots(\chi_{\mu(p)}%
(S^{1}\times D^{2}))\ldots)).
\]
We assume that the embeddings $\chi_{-},\chi_{+}$ are chosen so that this
collection is invariant under reflection in $\mathbb{R}^{2}\subset
\mathbb{R}^{3}$, as shown in Figure I-5, where $\mathbb{R}^{2}$ is drawn
vertically in $\mathbb{R}^{3}$.

The collection of $2$-spheres that arises in this step is gotten by
\textquotedblleft spinning\textquotedblright\ this original Bing collection of
1-spheres. Imagine $\mathbb{R}^{4}$ as being gotten from $\mathbb{R}%
^{3}=\mathbb{R}^{3}\times0\subset\mathbb{R}^{4}$, by spinning, or rotating,
$\mathbb{R}^{3}$ in $\mathbb{R}^{4}$ through $360^{0}$ (or $180^{0}$, if you
wish to be economical) about the plane $\mathbb{R}^{2}$, keeping
$\mathbb{R}^{2}$ fixed. Because of the symmetric positioning of the embeddings
$\{\chi_{\mu}\}$, each solid torus $S_{\mu}^{1}\times D^{2}$, when spun,
produces a thickened $2$-sphere $S_{\mu}^{2}\times D^{2}$. This collection
$\{S_{\mu}^{2}\times D^{2}\}$ of thickened $2$-spheres in $S^{2}\times D^{2}$
will be called the \textit{spun Bing collection} of thickened $2$-spheres.

It turns out that in Step 2, as described earlier, the $2^{p}$ $2$-sphere
components of $h_{2}(K\times\mathbb{R}^{1})\cap M\times0$ can be described as
the collection $\{\Sigma_{\mu}^{2}\mid\mu\in\{-,+\}^{p}\}\times0$, where
$\Sigma_{\mu}^{2}$ is the core $\Sigma_{\mu}^{2}\times0$ of the thickened
$2$-sphere $\Sigma_{\mu}^{2}\times D^{2}$ lying in $\Sigma^{2}\times D^{2}$,
all of this data being gotten by corresponding $\Sigma^{2}\times D^{2}$ to
$S^{2}\times D^{2}$.

This can be seen by careful analysis of the earlier description. What follows
is a description of a precise model which may make this clearer; this is the
way I originally perceived the construction. Since the model also will be used
to describe the generalization in Part II, it is presented here in its
generalized context.

The model starts in euclidean 3-space, which will be denoted $\mathbb{E}^{3}$
here so that there will be no erroneous correspondence made with
$\mathbb{R}^{3}$ in the description of spinning given above, or in the
Appendix. In $\mathbb{E}^{3}$ we define the 1-dimensional subsets $L_{-}^{1}$
and $L_{+}^{1}$ shown in Figure I-6b, each consisting of a circle with an
infinitely long tail attached, such that the circles are linked. Using
coordinates $(a,b,c)$ for $\mathbb{E}^{3}$ as shown, their precise
descriptions are as follows (where $\eta>0$ small): $L_{-}^{1}$ is the union
of the set of points in the $ac$-plane which are exactly distance $\eta$ from
the interval $0\times0\times\lbrack-2,0]$, together with the interval
$0\times0\times(-\infty,-2,-\eta]$ (pardon the nonstandard coordinates used in
Figure I-6, but they seem to yield the best pictures).

\begin{figure}[ptb]
\centerline{
\begin{tabular}{c}
\includegraphics{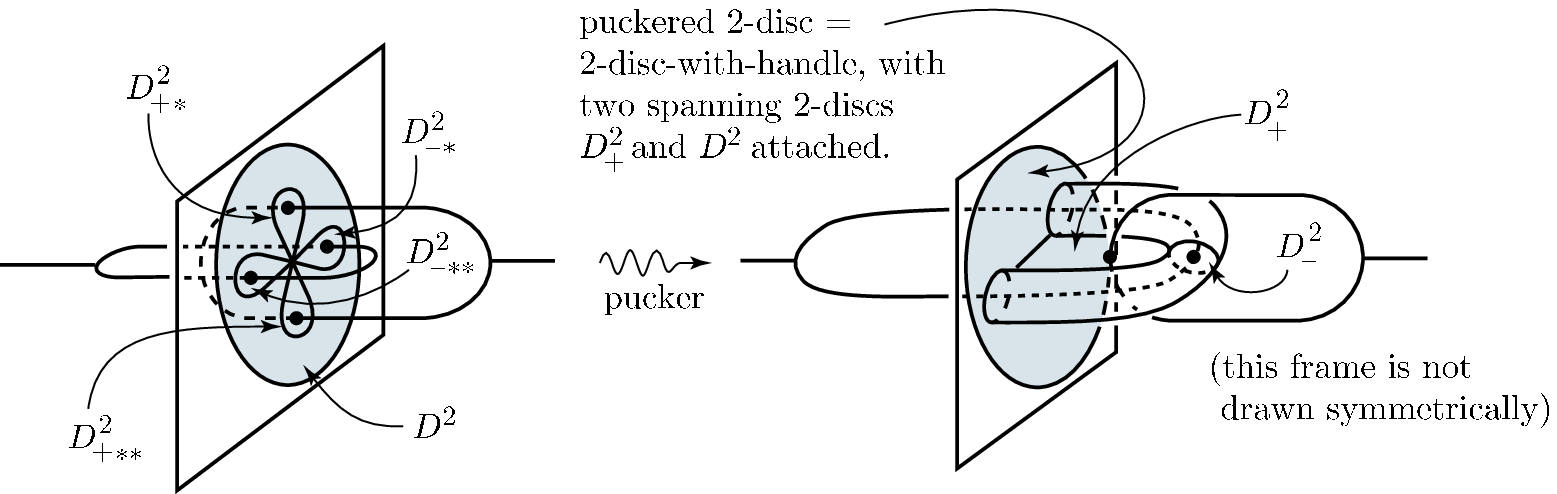}\\
Part (a).  The puckering operation: identify $D^2_{-*}$ to $D^2_{-**}$ to get $D^2_{-}$\\
and identify $D^2_{+*}$ to $D^2_{+**}$ to get $D^2_+$.  Note that the four\\
points of intersection become two points of intersection. \\
\bigskip{}\\
\includegraphics{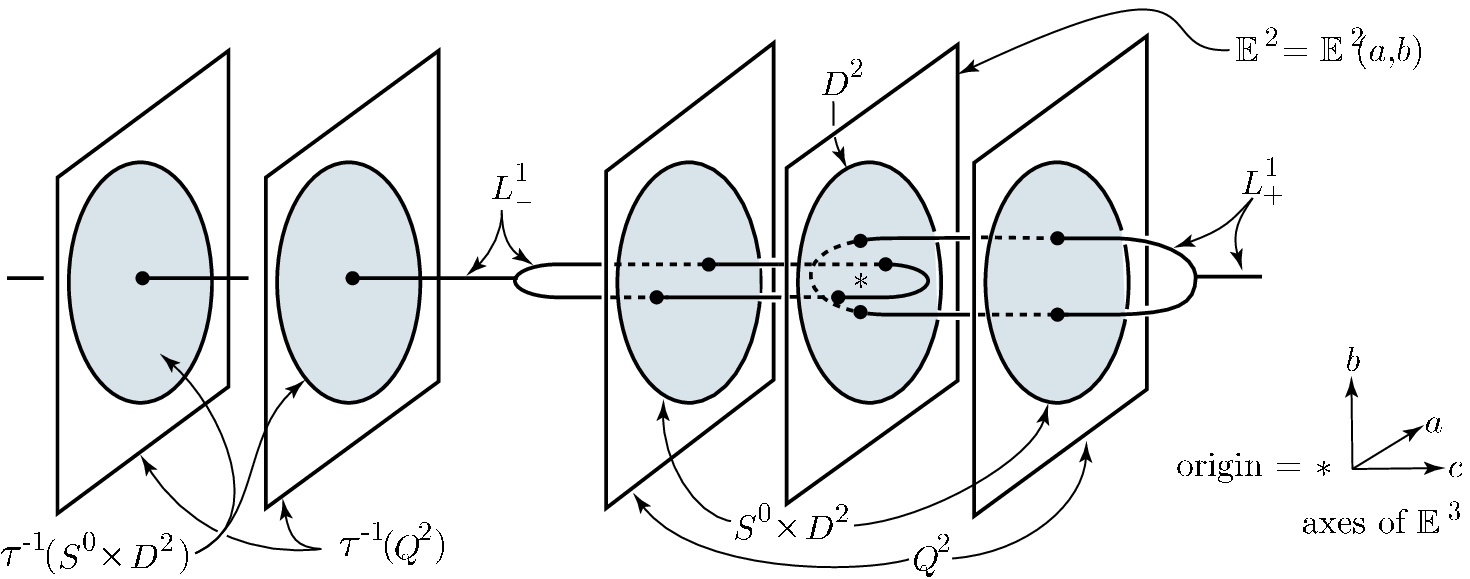}\\
Part (b).  The model 3-space $\mathbb{E}^3$\\
\bigskip{}\\
\includegraphics{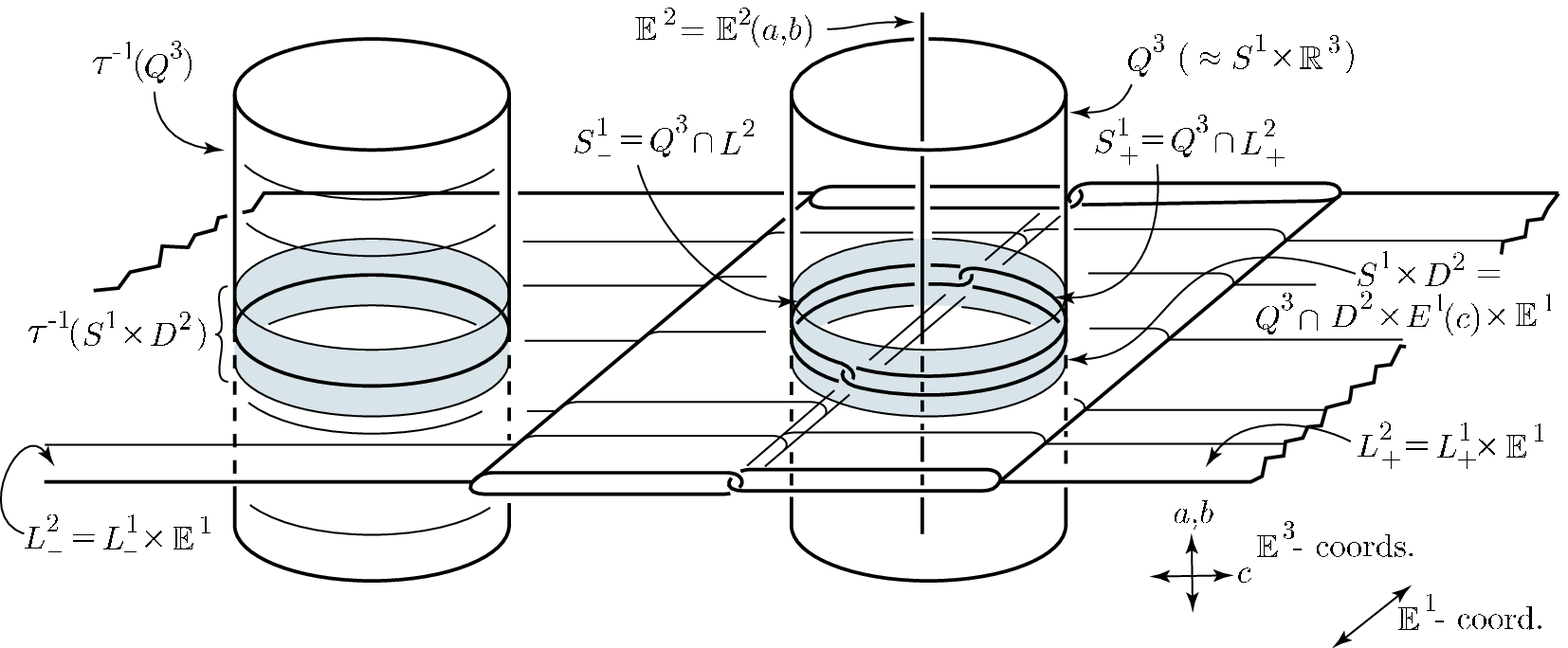}\\
Part (c).  The model 4-space $\mathbb{E}^4=\mathbb{E}^3\times \mathbb{E}^1$
\end{tabular}
}\caption{Some aspects of Step $i=2$}%
\end{figure}

Similarly, $L_{+}^{1}$ is the union of the set of points in the $bc$-plane
which are exactly distance $\eta$ from the interval $0\times0\times
\lbrack0,2]$, together with the interval $0\times0\times\lbrack2+\eta,\infty)$.

In what follows the two first stage $(1+k)$-spheres $S_{-}^{1+k}\cup
S_{+}^{1+k}$ of the $k$-times spun Bing collection of $(1+k)$-spheres in
$S^{1+k}\times D^{2}$ will be described, by means of intersections taking
place in $(4+k)$-space. It is probably best to understand first the $k=0$ case
(i.e. the original Bing case) and the $k=1$ case (which is the case of
interest in Part I). After describing the model, and the operation of
spherical mitosis in the model, it will be shown how the $k=1$ case of the
model corresponds to the current situation in $M^{4}\times\mathbb{R}^{1}$. The
relevant pictures are Figures I-6 b,c.

In euclidean $(4+k)$-space $\mathbb{E}^{4+k}=\mathbb{E}^{3}\times
\mathbb{E}^{1+k}$ (see above comment on the use of $\mathbb{E}$), consider the
sets $L_{-}^{2+k}\equiv L_{-}^{1}\times\mathbb{E}^{1+k};$ $L_{+}^{2+k}\equiv
L_{+}^{1}\times\mathbb{E}^{1+k};$ $\mathbb{E}^{2}=\mathbb{E}^{2}%
\times0=\mathbb{E}^{2}(a,b)\subset\mathbb{E}^{3}=\mathbb{E}^{3}\times
0\subset\mathbb{E}^{4+k}$; and $Q^{3+k}\equiv$ the boundary of the unit
tubular neighborhood of $\mathbb{E}^{2}$ in $\mathbb{E}^{4+k}$. Thus
$\mathbb{E}^{2}$ is a $2$-plane in $\mathbb{E}^{3}$ perpendicular to the
$c$-axis, and $Q^{3+k}$ is naturally homeomorphic to $\mathbb{E}^{2}\times
S^{1+k}$, where we are thinking of $S^{1+k}$ as the unit sphere in
$0\times\mathbb{E}^{1}(c)\times\mathbb{E}^{1+k}$, where $\mathbb{E}^{1}(c)$,
is the $c$-coordinate axis in $\mathbb{E}^{3}=\mathbb{E}^{2}(a,b)\times
\mathbb{E}^{1}(c)$. (Note that $\mathbb{E}^{3}$ itself can be regarded as the
$k=-1$ case of this construction.)

Let $D^{2}$ denote the unit disc in $\mathbb{E}^{2}$, so that $D^{2}\times
S^{1+k}$ (henceforth denoted $S^{1+k}\times D^{2})$ is the intersection of
$Q^{3+k}$ with $D^{2}\times\mathbb{E}^{1}(c)\times\mathbb{E}^{1+k}$ in
$\mathbb{E}^{3}\times\mathbb{E}^{1+k}$. The significance of this model is the following

\textbf{Observation.} \textit{The pair $(S^{1+k}\times D^{2},S^{1+k}\times
D^{2}\cap(L_{-}^{2+k}\cup L_{+}^{2+k}))$ corresponds in a natural manner to
the pair $(S^{1+k}\times D^{2},S_{-}^{1+k}\cup S_{+}^{1+k})$, where $S^{1+k}$
and $S_{+}^{1+k}$ are the two first-stage $(1+k)$-spheres in the $k$-times
spun Bing collection of $(1+k)$-spheres in $S^{1+k}\times D^{2}$.}

Understanding this is a matter of analyzing the model, seeing first the $k=0$
case and then proceeding to higher dimensions.

We wish to describe a process by which these two $(1+k)$-spheres $S_{-}%
^{1+k}\cup S_{+}^{1+k}$ can be gotten from the original core $(1+k)$-sphere
$S^{1+k}\times0\subset S^{1+k}\times D^{2}$ in a continuous manner. Let
$\tau:\mathbb{E}^{3}\times\mathbb{E}^{1+k}\rightarrow\mathbb{E}^{3}%
\times\mathbb{E}^{1+k}$ be translation by 4 units in the positive direction of
the $c$-axis of $\mathbb{E}^{3}$. Let $\tau_{t}$, $0\leq t\leq1$, be the
natural linear isotopy joining $\tau_{0}=identity$ to $\tau_{1}=\tau$. We wish
to focus on the translated copy $\tau^{-1}(S^{1+k}\times D^{2})$ of
$S^{1+k}\times D^{2}$, and to examine its intersection with the always-fixed
subset $L_{-}^{2+k}\cup L_{+}^{2+k}$ of $\mathbb{E}^{4+k}$, as $\tau
^{-1}(S^{1+k}\times D^{2})$ is translated back to its standard position
$S^{1+k}\times D^{2}$ by the isotopy $\{\tau_{t}\}$. (The reason for this
somewhat backward point of view, i.e. moving $S^{1+k}\times D^{2}$ instead of
moving $L_{-}^{2+k}\cup L_{+}^{2+k}$, is that in the real situation in
$M^{4}\times\mathbb{R}^{1}$ this is what is happening, as we will see.)
Examination reveals that the intersection $\tau_{t}(\tau^{-1}(S^{1+k}\times
D^{2}))\cap(L_{-}^{2+k}\cup L_{+}^{2+k})$ starts out at $t=0$ looking like the
core sphere of $\tau^{-1}(S^{1+k}\times D^{2})$, and then as $t$ increases the
intersection undergoes a transformation, dividing into two components, so that
ultimately at time $t=1$ it has become the pair of linked $(1+k)$-spheres
$S_{-}^{1+k}\cup S_{+}^{1+k}$ in $\tau_{1}(\tau^{-1}(S^{1+k}\times
D^{2}))=S^{1+k}\times D^{2}$. This is the process of \textit{spherical
meitosis}. It is interesting to follow the intermediate stages, even in the
original Bing $k=0$ case.

The $k=1$ case of the above-described model is corresponded to our situation
in $M^{4}\times\mathbb{R}^{1}$ in the following manner. The space
$\mathbb{E}^{3}$ is to be thought of as an open subset of the boundary
$\partial D^{4}$ of the 0-handle $D^{4}$ of $M^{4}$, so that $\mathbb{E}%
^{3}\cap L^{1}=L_{-}^{1}\cup L_{+}^{1}$. The fourth coordinate $\mathbb{E}%
^{1}$ of $\mathbb{E}^{3}\times\mathbb{E}^{1+1}$ is to be thought of as the
$(0,1)$-coordinate in $D^{4}-\ast=\partial D^{4}\times\lbrack0,1)$, with the
origin of $\mathbb{E}^{1}$ corresponding say to the point $\frac{1}{2}%
\in\lbrack0,1)$. Finally, the fifth coordinate $\mathbb{E}^{1}$ of
$\mathbb{E}^{3}\times\mathbb{E}^{1+1}$ is to be thought of as the
$\mathbb{R}^{1}$ coordinate of $M^{4}\times\mathbb{R}^{1}$, translated so that
the origin of $\mathbb{E}^{1}$ is at the $2$-level (or later on at the
4-level, or 6-level, or whatever higher level of $\mathbb{R}^{1}$ one is
working in at the time). In summary, the first four coordinates of
$\mathbb{E}^{3}\times\mathbb{E}^{1+1}$ are to be thought of as defining an
open subset of the 0-handle $D^{4}$ of $M^{4}$, and the fifth coordinate is to
be thought of as the vertical coordinate of $M^{4}\times\mathbb{R}^{1}$.
(Aside: Actually, there is no mathematical justification for distinguishing
the last two coordinates of $\mathbb{E}^{3}\times\mathbb{E}^{1+1}$ from each
other, since the construction in the model is symmetric about $\mathbb{E}^{3}%
$. But it is conceptually helpful to make this distinction.)

We now can use this model to understand the construction of $h_{1,1}$ from
$h_{1,0}=h_{1}$. Let $\mathbb{E}_{\#}^{2}=\operatorname*{int}D_{\#}^{2}%
\times\frac{1}{2}\times1$ denote the interior of the $2$-disc part of $A_{1}$.
Then $\mathbb{E}_{\#}^{2}\subset M\times1$, and $\mathbb{E}_{\#}^{2}\cap
K\times\mathbb{R}^{1}$ is the center point $d\times\frac{1}{2}\times1$ of
$\mathbb{E}_{\#}^{2}$. Suppose distances are scaled so that the portion of
$h_{1}^{-1}(M\times0)$ lying near the $2$-sphere component $h_{1}^{-1}%
(\Sigma^{2}\times0)$ of $K\times\mathbb{R}^{1}\cap h_{1}^{-1}(M^{4}\times0)$
looks like the boundary of the unit tubular neighborhood of $\mathbb{E}%
_{\#}^{2}$ in $M^{4}\times\mathbb{R}^{1}$. Forget the 3-stage transformation
of $A_{1,0}=A_{1}$ to $A_{1,1}$; instead let us just isotope $A_{1,0}$ to a
new position by isotoping $\mathbb{E}_{\#}^{2}$ in $M^{4}\times\mathbb{R}^{1}%
$, keeping it fixed near its boundary, so that a neighborhood in
$\mathbb{E}_{\#}^{2}$ of the center point $d\times\frac{1}{2}\times1$ of
$\mathbb{E}^{2}$ is moved up to the $2$-level $M\times2$, there to coincide
with a large compact piece of the plane $\mathbb{E}^{2}=\mathbb{E}^{2}(a,b)$
of the model. During the first part of this isotopy, the intersection of
$\mathbb{E}_{\#}^{2}$ with $K^{2}\times\mathbb{R}^{1}$ is to be kept always
the center point of $\mathbb{E}_{\#}^{2}$. As $\mathbb{E}_{\#}^{2}$ nears the
end of its journey, in a neighborhood of the $M^{4}\times2$ level, we see it
near the 0-level of the model. There we suppose that the end of its isotopy
coincides with the translation by $\tau_{t}$ of the plane $\tau^{-1}%
(\mathbb{E}^{2})$ in the model to its standard position $\mathbb{E}^{2}$.
Hence, near the end of the isotopy the intersection of the moving
$\mathbb{E}_{\#}^{2}$ with the fixed $K\times\mathbb{R}^{1}$ changes from a
single point to four points. But that is not important. What is important is
that the movement of $h_{1,0}^{-1}(M\times0)$ to $h_{1,1}^{-1}(M\times0)$ can
be viewed as the movement of the frontier of a regular neighborhood of this
moving plane $\mathbb{E}_{\#}^{2}\subset A_{1,0}$, where at the end of the
isotopy the portion of this regular neighborhood we see in the model is the
1-neighborhood of the moving plane $\tau_{t}(\tau^{-1}(\mathbb{E}^{2}))$. So
the point is this: qualitatively, the change in the intersection of
$h_{1,0}^{-1}(M\times0)$ with $K\times\mathbb{R}^{1}$, as it is isotoped to to
$h_{1,1}^{-1}(M\times0)$, is the same as the change of the intersection of
$\tau^{-1}(Q^{4})$ with $L_{-}^{3}\cup L_{+}^{3}$ in the model $\mathbb{E}%
^{3}\times\mathbb{E}^{2}$, as $\tau^{-1}(Q^{4})$ is isotoped by the
translations $\{\tau_{t}\}$ back to its standard position $Q^{4}$. Hence the
change in the intersection amounts to spherical meitosis.

From this, the relation of the two $2$-sphere components of $h_{1,1}%
(K\times\mathbb{R}^{1})\cap M\times0$ to the single $2$-sphere component of
$h_{1,0}(K\times\mathbb{R}^{1})\cap M\times0$ can be seen to be modelled on
the spun Bing construction. (As an incidental remark, to change the above
repositioned $A_{1,0}$ into the earlier-described $A_{1,1}$, one only has to
pucker (cf. earlier) the newly positioned $\mathbb{E}_{\#}^{2}$, to reduce its
four pints of intersection with $K^{2}\times\mathbb{R}^{1}$ to two points of
intersection. See Figure I-6a.)

This same analysis works for the later stages of this step, i.e. going from
$h_{1,1}$ to $h_{1,2}$, etc.

\textbf{Step i=3 }(second round). In this step, there are no qualitative
changes made in $h_{2}$ to get $h_{3}$, but instead the diameters of the
various $2$-sphere components of $h_{2}(K\times\mathbb{R}^{1})\cap M\times0$
are made $\delta$-small, by following $h_{2}$ by a level-preserving
homeomorphism of $M\times\mathbb{R}^{1}$. This level-preserving homeomorphism
is obtained in the obvious manner from the ambient isotopy of $\Sigma
^{2}\times D^{2}\;rel\;\partial$ used to shrink small the $2^{p}$
$p^{\text{th}}$ stage $2$-spheres in the spun Bing collection of $2$-spheres.
This shrinking of the spun Bing collection is explained in the Appendix to
Part II. It is for this shrinking argument that $p=p(\delta)$ must be chosen
large. Note that $p$ can in fact be chosen at the start of Step 2 (as it must
be) because, in addition to $\delta$, it depends only on how $\Sigma^{2}\times
D^{2}$ is embedded in $M$, and that embedding is chosen after the construction
of $h_{1}$.

\textbf{Step i=4 }(second round). The effect of this step is to take each
$\delta$-small $2$-sphere component $\Sigma_{\mu}^{2}\times0$ of
$h_{3}(K\times\mathbb{R}^{1})\cap M\times0$, and to replace it by two copies
of duncehats, $h_{4}(K^{2}\times(t_{\mu}\pm\lambda))$, as suggested earlier.
See Figure I-7.

\begin{figure}[th]
\centerline{
\includegraphics{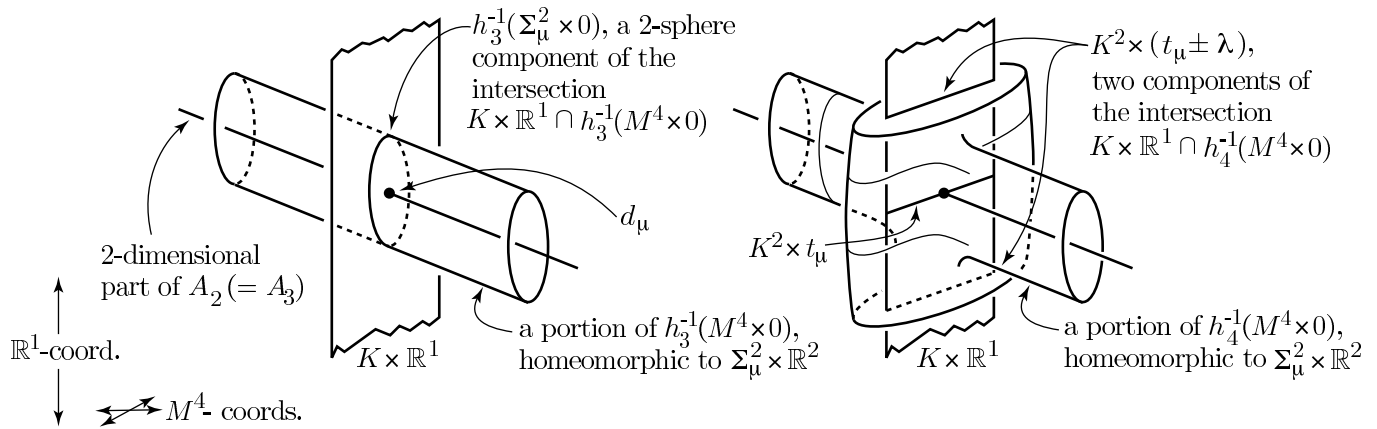}
}\caption{Step $i=4$: the passage from $h_{3}$ to $h_{4}$}%
\end{figure}

The important thing is that these two new duncehat components of intersection
can be chosen to lie in any arbitrarily small neighborhood $\Sigma_{\mu}%
^{2}\times D^{2}\times0$ of $\Sigma_{\mu}^{2}\times0$ in $M^{4}\times0$, hence
their size is automatically controlled. If one examines the isotopy of $h_{3}$
to $h_{4}$, engendered by the expansion-collapse of $A_{3}$ to $A_{4}$, one
can see the component-of-intersection $\Sigma_{\mu}^{2}\times0$ undergo
meiosis in $\Sigma_{\mu}^{2}\times D^{2}\times0$, becoming two duncehats which
are linked there.

This completes the proof of the Basic Lemma.
\end{proof}

\indent Given the Basic Lemma, the remainder of the proof of the Shrinking
Proposition is patterned on \cite[Section 3]{EM}. Figures 1, 3, and 4 there
are meaningful here, too. More precise details can be gotten from that paper.
Recall that $B^{2}\times D^{2}$ denotes the $2$-handle of $M^{4}$.

\textsc{Window Building Lemma. }\emph{Given }$\delta>0$\emph{, there is a
homeomorphism }$h_{\#}:M\times\mathbb{R}^{1}\rightarrow M\times\mathbb{R}^{1}%
$\emph{, fixed on }$\partial M\times\mathbb{R}^{1}$\emph{, such that}

\begin{enumerate}
\item \emph{for each }$j\in\mathbb{Z}$\emph{ and each }$t\in\lbrack j-1,j+1]$,
$h_{\#}(M\times t)\subset M\times\lbrack j-1,j+1],$\emph{ and }

\item \emph{for each }$t\in\mathbb{R}^{1}$\emph{, if }$h_{\#}(K\times t)\cap
B^{2}\times D^{2}\times\lbrack j-1+\delta,j+1-\delta]\neq\emptyset$\emph{ for
any }$j\in2\mathbb{Z}$\emph{, then }$\operatorname*{diam}h_{\#}(K\times
t)<\delta$\emph{.}\textit{ }
\end{enumerate}

\begin{proof}
This follows quickly from the Basic Lemma. We can suppose that the
homeomorphism $h_{4}$ constructed there has the additional property that the
punctured duncehat component of $h_{4}(K^{2}\times\mathbb{R}^{1})\cap
M^{4}\times0$ does not intersect $B^{2}\times D^{2}\times0$. This either can
be verified by construction (e.g., it would be natural to make $h_{1}$ have
this property), or one can argue that the pierced duncehat component, being
$\delta$-small, can be assumed to lie in some ball in $M^{4}\times0$, and
therefore it can be isotoped off of $B^{2}\times D^{2}\times0$ in $M^{4}%
\times0$. Supposing this, then by applying after $h_{4}$ a vertical
homeomorphism of $M^{4}\times\mathbb{R}^{1}$ which basically expands
$B^{2}\times D^{2}\times\lbrack-\eta,\eta]$ to $B^{2}\times D^{2}\times
\lbrack-1+\delta,1-\delta]$ (i.e., \textquotedblleft opens the window
wider\textquotedblright), for some small $\eta>0$, one can obtain a
homeomorphism $h_{\ast}:M\times\mathbb{R}^{1}\rightarrow M\times\mathbb{R}%
^{1}$, fixed on $\partial M\times\mathbb{R}^{1}$ and having compact support,
such that $h_{\ast}$ satisfies condition (2) above for the value $j=0$. To
complete the proof, one conjugates $h_{\ast}$ to make it have the additional
property that its support lies in $M\times\lbrack-1,1]$, and then one obtains
$h_{\#}$ by taking the sum of an infinite number of translates of $h_{\ast}$
stacked on top of each other.
\end{proof}

\begin{proof}
[Proof of Shrinking Proposition from the Window Building Lemma]Given
$\epsilon>0$, choose a small ball $B^{4}$ in $\operatorname*{int}M^{4}$, with
$\operatorname*{diam}B^{4}<\epsilon/2$. We indicate how to construct a
homeomorphism $h:M\times\mathbb{R}^{1}\rightarrow M\times\mathbb{R}^{1}$,
fixed on $\partial M\times\mathbb{R}^{1}$, which satisfies the following
weakened versions of the conditions from the Shrinking Proposition: for each
$t\in\mathbb{R}^{1}$,

\begin{itemize}
\item[1$^{\prime}$.] $h(M\times t)\subset M\times\lbrack t-3,t+3]$, and

\item[2$^{\prime}$.] either

\begin{itemize}
\item[a.] $h(K^{2}\times t)\subset B^{4}\times\lbrack t-3,t+3]$, or

\item[b.] $\operatorname*{diam}h(K^{2}\times t)<\epsilon$
\end{itemize}
\end{itemize}

From this weaker version of the Shrinking Proposition it is clear that the
original version follows, simply by rescaling the vertical coordinate.

\begin{itemize}
\item[1.] To construct this $h$, first one constructs a uniformly continuous
homeomorphism $g:M\times\mathbb{R}^{1}\rightarrow M\times\mathbb{R}^{1}$,
fixed on $\partial M\times\mathbb{R}^{1}$, such that for each $t\in
\mathbb{R}^{1},g(M\times t)\subset M\times\lbrack t-1,t+1]$, and

\item[2.] the image under $g$ of $(M_{1}^{4}\cap(D^{4}\cup B^{1}\times
D^{3}))\times\mathbb{R}^{1}\cup\bigcup\{M_{1}\times(j+1)\mid j\in
2\mathbb{Z}\}$ lies in $B^{4}\times\mathbb{R}^{1}$, where $M_{1}^{4}%
=M^{4}-\partial M\times\lbrack0,1)$ for some collar $\partial M\times
\lbrack0,2)$ of $\partial M$ in $M^{4}$, and where $D^{4}\cup B^{1}\times
D^{3}$ is the union of the 0-handle and the 1-handle of $M^{1}$.
\end{itemize}

The details for $g$ (which are simple) are omitted, since a more general such
$g$ is constructed in Part II. Given $g$, then one can let $h=gh_{\#}$, where
$h_{\#}$ is provided by the Window Building Lemma for some sufficiently small
value of $\delta=\delta(\epsilon,g)$, and where we are assuming without loss
that $h_{\#}(K\times\mathbb{R}^{1})\subset M_{1}\times\mathbb{R}^{1}$. This
completes the proof of the Shrinking Proposition, and hence Part I.
\end{proof}

\part{The double suspension of any homology $n$--sphere which bounds a
contractible $\left(  n-1\right)  $--manifold is a sphere}

\indent The purpose of this part is to generalize Part I to prove

\begin{theorem}
\label{double-suspension-of-boundary}\textit{The double suspension $\Sigma
^{2}H^{n}$ of any homology sphere $H^{n}$ which bounds a contractible
topological manifold is homeomorphic to a sphere. In particular, if $n\geq4$,
this conclusion holds.}
\end{theorem}

\indent The last sentence is justified by the following Proposition (which
basically is known; see proof). It is stated in sufficient generality for use
in the Prologue, Section II.

\begin{proposition}
\label{homology-spheres-bound}Suppose $H^{n}$ is a compact space such that
$H^{n}\times S^{1}$ is a topological manifold-without-boundary, and such that
$H^{n}$ has the integral homology groups of the $n$-sphere. If $n\geq4$, then
there is an embedding of $H^{n}\times\mathbb{R}^{1}$ into some open
contractible topological manifold $M^{n+1}$ so that $M^{n+1}-(H^{n}%
\times\mathbb{R}^{1})$ is compact (i.e., $H^{n}\times\mathbb{R}^{1}$ is a
neighborhood of the end of $M^{n+1}$).
\end{proposition}

\indent The Proposition is unknown for $n=3$. If $n\geq4$, and $H^{n}$ is a
$PL$ manifold homology $n$-sphere, then Kervaire proved in \cite[Cor, p.
71]{Ke} that $H^{n}$ bounds a contractible $PL$ manifold. Using the post-1968
knowledge that any (topological manifold) homology $n$-sphere, $n\geq5$, is a
$PL$ manifold, then Kervaire's proof covers this case, too. (Aside: it is
unknown whether every homology $4$-sphere is a $PL$ manifold, but topological
immersion theory readily establishes, by immersing $H^{n}%
-\operatorname*{point}$ into $\mathbb{R}^{4}$, that $H^{n}%
-\operatorname*{point}$ is a $PL$ manifold).

\indent The virtue of the following proof is its brevity.

\begin{proof}
[Proof of Proposition \ref{homology-spheres-bound}]It is rudimentary homotopy
theory that one can attach to $H^{n}$ a finite number of $2$-cells, and then
an equal number of $3$-cells, to make $H^{n}$ homotopically a sphere. Doing
these attachings as surgeries in a band $H^{n}\times(0,1)$ in $H^{n}\times
S^{1}$, one produces from $H^{n}\times S^{1}$ a new manifold, say $G$, which
is homotopically equivalent to $S^{n}\times S^{1}$. (The one nontrivial aspect
of this argument is embedding the attaching $1$- and $2$-spheres to have
product neighborhoods (i.e., trivial normal bundles). Perhaps this is most
easily handled by using the above-mentioned immersion proof to put a smooth
structure on $(H^{n}-\operatorname*{point})\times(0,1)$, and then doing all
surgeries there, using well-known smooth arguments). Now the elementary
argument in \cite[Appendix I]{Si6} establishes that $\widetilde{G}$, the
universal cover of $G$, is homeomorphic to $S^{n}\times\mathbb{R}^{1}$. Let
$F=\widetilde{G}\cup\infty\approx\mathbb{R}^{n+1}$ be gotten by compactifying
one end of $\widetilde{G}$. Then we can let $M^{n+1}\subset F$ be the bounded
open complementary domain of $F-H^{n}\times0$, where $H^{n}\times0$ denotes
any of the infinitely many natural copies of $H^{n}$ in $F$.
\end{proof}

\indent As indicated in the Preliminaries, Theorem
\ref{double-suspension-of-boundary} follows from the more general

\begin{theorem}
[single cell-like set version]\label{single-cell-like-set}\textit{Suppose $X$
is a cell-like set in a manifold-without-boundary $M^{m}$. If $m\geq4$, then
the stabilized quotient map $\pi\times id_{\mathbb{R}^{1}}:M\times
\mathbb{R}^{1}\rightarrow M/X\times\mathbb{R}^{1}$ is approximable by
homeomorphisms.}
\end{theorem}

\indent The $m=3$ version of this theorem was established in \cite{EP}
\cite{EM} under the additional assumption that $X$ has a irreducible
3-manifold neighborhood in $M$, thus avoiding the Poincar\'{e} conjecture. The
$m=2$ version, which does not require stabilization, is considered classical.
From the above-stated version of Theorem \ref{single-cell-like-set}, one
readily deduces the somewhat more general \textquotedblleft closed-$0$%
-dimensional\textquotedblright\ version, as in \cite{EP} \cite{EM}, explained
again in \cite{Ed2}, but there is no reason to elaborate this here.

\indent The reader will recognize the construction below as a straightforward
generalization of that in Part I, a fact which I realized in the month
(January 1975) following the completion of Part I. But the $m=4$ case was
elusive, and it wasn't until several months later (August 1975) that I
realized that the solution there was to use the freedom of the extra
$\mathbb{R}^{1}$-coordinate.

\indent In case $X$ happens to be a codimension 2 polyhedron in $M$ (e.g., if
$M$ is a compact contractible $PL$ manifold, $\dim M \geq5$, then $M$ has such
a spine $X$), then several steps in the following proof become trivial, so
this is a good case to keep in mind.

\indent Theorem \ref{single-cell-like-set} is proved by showing that the Bing
Shrinking Criterion is satisfied, i.e., by proving\medskip

\textsc{Shrinking Proposition. }\emph{Suppose }$X$\emph{ is a cell-like set in
a manifold-without-boundary }$M^{m}$\emph{. If }$m\geq4$\emph{, then given any
neighborhood }$U$\emph{ of }$X$\emph{ and any }$\varepsilon>0$\emph{, there is
a homeomorphism }$h:M\times\mathbb{R}^{1}\rightarrow M\times\mathbb{R}^{1}%
$\emph{, fixed on }$(M-U)\times\mathbb{R}^{1}$\emph{, such that for each
}$t\in\mathbb{R}^{1}$\emph{,}

\begin{itemize}
\item[1.] $h(U\times t)\subset U\times\lbrack t-\epsilon,t+\epsilon]$\emph{,
and}

\item[2.] $\operatorname*{diam}h(X\times t)<\epsilon$\emph{.\medskip}
\end{itemize}

Succinctly stated, the idea of the proof is to produce a neighborhood basis
for $X$ in $M$ which is sufficiently standard in some sense, so that the
constructions from Part I can be applied. Details follow.

\indent Since $X$ is cell-like, then $X$ has a $PL$ triangulable neighborhood
in $M^{m}$ (even when $m=4)$. The easiest proof of this is the one which uses
topological immersion theory to immerse some neighborhood of $X$ into
$\mathbb{R}^{m}$ (cf. proof of Proposition above). So without loss $M$ is a
$PL$ manifold.

\indent It is a standard fact that $X$ has an arbitrarily small compact $PL$
manifold neighborhood $N$ in $M$ such that $N$ has an $(m-2)$-dimensional
spine. This is proved by taking an arbitrary compact $PL$ manifold
neighborhood $N_{\ast}$ of $X$, and by isotoping its dual $1$-skeleton off of
$X$, by an ambient isotopy supported in $\operatorname*{int}N_{\ast}$, using
the fact that $X$ is cell-like and $N_{\ast}-X$ has one end. Then a small
neighborhood of the repositioned dual 1-skeleton of $N_{\ast}$ can be deleted
from $N_{\ast}$ to produce $N$. This argument used the dimension restriction
$m\geq4$ in constructing the isotopy. (Recall that when $m=3$ the existence of
such a neighborhood for $X$ a contractible $2$-dimensional polyhedron would
imply the Poincar\'{e} conjecture. However, in this dimension, if one
hypothesizes in addition that $X$ has an irreducible $3$-manifold
neighborhood, then such $1$-spine neighborhoods exist; see \cite[Lemma
1]{McM1}.)

\indent Letting $N$ be such a neighborhood of $X$, we can write $N=L\cup
(\bigcup_{\alpha=1}^{r}H_{\alpha})$ (see Figure II-1), where $L$ is a compact
manifold with an $(m-3)$-dimensional spine, and where the $H_{\alpha}$'s are
disjoint handles of index $m-2$, attached to $L$ so that for each $\alpha$,
$L\cap H_{\alpha}=\delta H_{\alpha}\subset\partial L$, where $\delta
H_{\alpha}$ is the \emph{attaching-boundary} of $H_{\alpha}$, defined by
$\delta H_{\alpha}=\partial D_{\alpha}^{m-2}\times D^{2}\subset D_{\alpha
}^{m-2}\times D^{2}=H_{\alpha}$. The \emph{core} of the handle is the
$(m-2)$-cell $D_{\alpha}^{m-2}=D_{\alpha}^{m-2}\times0\subset D_{\alpha}%
^{m-2}\times D^{2}=H_{\alpha}$. See Figure II-1.

\begin{figure}[th]
\centerline{
\includegraphics{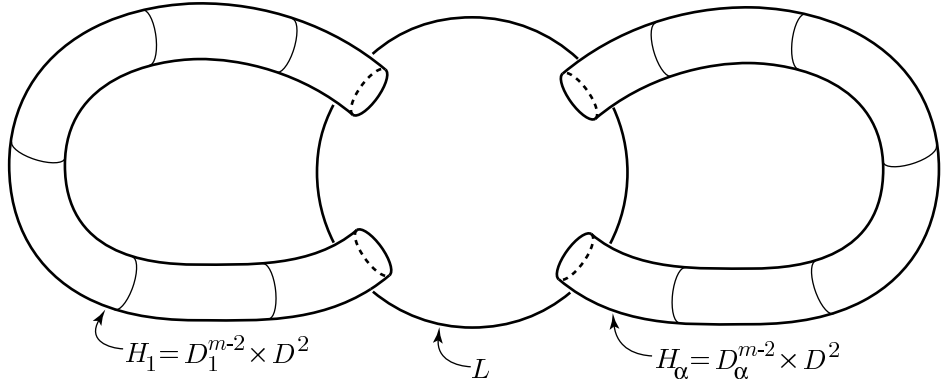}
}\caption{The neighborhood $\displaystyle N=L\cup\bigcup_{\alpha=1}%
^{r}H_{\alpha}$}%
\end{figure}

\indent This neighborhood $N$ can be thought of as a higher dimensional
analogue of the familiar 3-dimensional cube-with-handles neighborhood.

\indent The primary goal is to establish

\begin{lemma}
[Building a Single Window]\label{building-a-window}Suppose the data $X\subset
N^{m}$ as above. Let $H_{a}$ be any one of the $\left(  m-2\right)  $-handles
of $N$. Then given any $\delta>0$, there is a homeomorphism $h_{\alpha
}:N\times\mathbb{R}^{1}\rightarrow N\times\mathbb{R}^{1}$ such that

\begin{enumerate}
\item $h_{\alpha}$ \textit{has compact support in $\operatorname*{int}%
N\times\mathbb{R}^{1}$, and}

\item \textit{for each $t\in\mathbb{R}^{1}$, if $h_{\alpha}(X\times t)\cap
H_{\alpha}\times\lbrack-1,1]\neq\emptyset$, then} $\operatorname*{diam}%
h_{\alpha}(X\times t)<\delta$.
\end{enumerate}
\end{lemma}

\textsc{Note. } $h_{\alpha}$ can be regarded as the analogue of the
homeomorphism $h_{\ast}$ constructed during the proof of the Window Building
Lemma, Part I (at least in the case $N$ has only one $(m-2)$-handle
$H_{\alpha}$).

\indent The first task is to describe a certain \emph{model handle} sequence
which will be used in the proof of Lemma 1. The 3-dimensional version of the
model is a solid cylinder $H^{3}(0)=D^{1}\times D^{2}$ (to be thought of as a
1-handle), containing the familiar infinite sequence of linked sub-1-handles
as shown in Figure II-2. That is, $H^{3}(0)$, together with its subhandles,
amount to the Bing collection of solid tori, cut in half by the vertical plane
shown in Figure I-5. The \emph{attaching-boundary} of $H^{3}(0)$ is the union
of the two end discs, $\delta H^{3}(0)=\partial D^{1}\times D^{2}$. The union
of the $2^{p}$ sub-1-handles of $H^{3}(0)$ at stage $p$ is denoted $H^{3}(p)$,
with attaching-boundary $\delta H^{3}(p)=H^{3}(p)\cap\delta H^{3}(0)$.

\begin{figure}[th]
\centerline{
\includegraphics{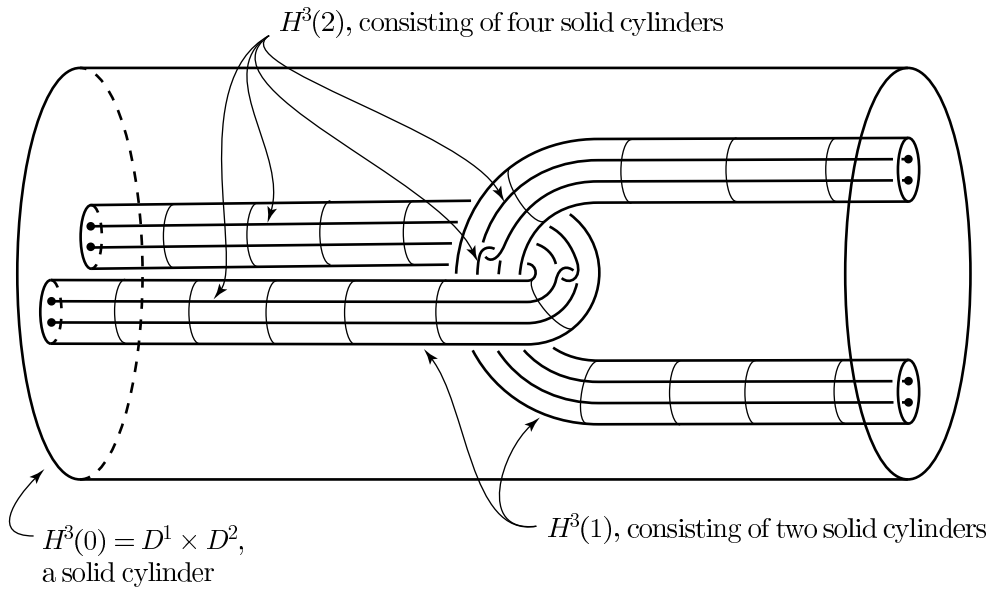}
}\caption{The model handle sequence $H^{3}(0) \supset H^{3}(1) \supset
H^{3}(2) \supset\cdots$}%
\end{figure}

\indent In higher dimensions, the initial \emph{model handle} is
$H^{m}(0)\equiv H^{3}(0)\times D^{m-3}$ (which we really want to think of as
$(D^{1}\times D^{m-3})\times D^{2})$, and the \emph{attaching-boundary} of
$H^{m}(0)$ is $\delta H^{m}(0)=\delta H^{3}(0)\times D^{m-3}\cup
H^{3}(0)\times\partial D^{m-3}$ (which can then be thought of as the thickened
$(m-3)$-sphere $\partial(D^{1}\times D^{m-3})\times D^{2})$. Let
$H^{m}(p)=H^{3}(p)\times D^{m-3}$ and let $\delta H^{m}(p)=H^{m}(p)\cap\delta
H^{m}(0)$. It is important to realize that $\delta H^{m}(p)$ is \textbf{not}
just $\delta H^{3}(p)\times D^{m-3}$, but also includes $H^{3}(p)\times
\partial D^{m-3}$. \emph{Hence for }$m\geq1$\emph{, the subset }$\delta
H^{m}(p)$\emph{ of }$\delta H^{m}(0)$\emph{ looks like the }$p^{th}$\emph{
stage of the spun Bing collection of thickened }$(m-3)$\emph{-spheres in
}$S^{m-3}\times B^{2}$\emph{.}

\indent We will need ramified versions of these models, produced in the spirit
of R. Daverman's ramified cantor set constructions \cite[Lemma 4.1]{Da1}. Any
stage in the construction of the sequence $\{H^{3}(p)\}$, starting with the
$0^{\text{th}}$ stage, may be \emph{ramified}, which means that each handle of
that stage is replaced by several adjacent, parallel copies of itself.
Examples are shown in Figure II-3.

\begin{figure}[th]
\centerline{
\includegraphics{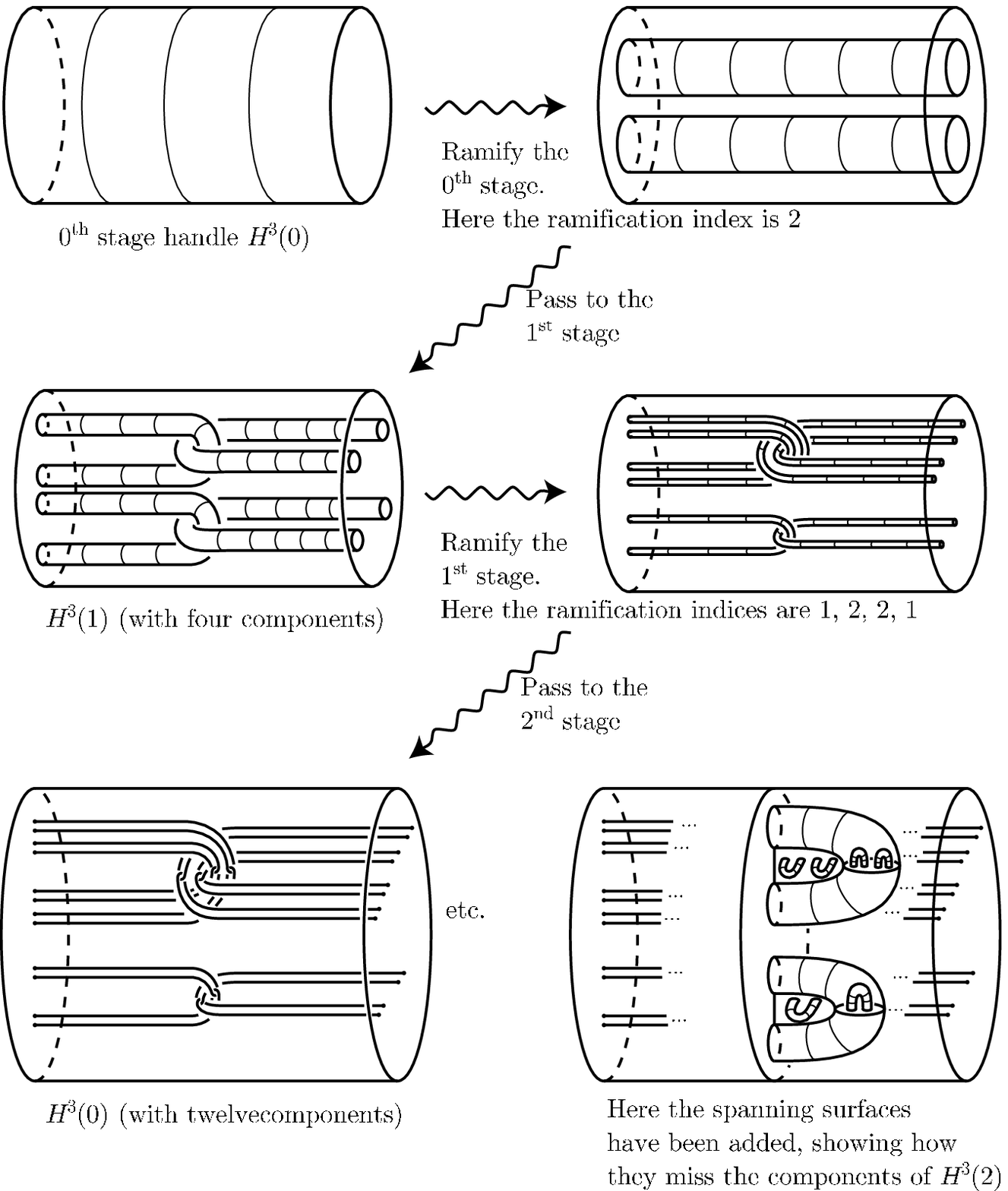}
}\caption{Producing a model ramified handle sequence}%
\end{figure}

There may be a different number of ramifications performed in each
handle-component at each stage; the only restriction is that the number of
ramifications at each stage be finite. These ramification indices will not be
recorded in any way, and our notation for a ramified handle will be the same
as that for an unramified handle. It will become clear why these ramifications arise.

\indent A model ramified $(m-2)$-handle sequence $\{H^{m}(p)\}$ is defined to
be $\{H^{3}(p)\times D^{m-3}\}$, as before, where now $\{H^{3}(p)\}$ is a
model ramified 1-handle sequence. In particular there is no ramification in
the $D^{m-3}$-coordinate (whatever that would mean).

\indent The reason for these model handles is for use in the following Lemma.
From now on, we restrict attention to the $\dim M \geq5$ case, returning to
the $\dim M = 4$ case at the end of the proof.

\begin{lemma}
[Repositioning]\label{repositioning}\textit{Suppose $X,N^{m}$ and $H_{\alpha}$
are as in Lemma 1, and suppose $m\geq5$. Then for any $p\geq0$, there is a
repositioning of $X$ in $\operatorname*{int}N$ (by an unnamed ambient isotopy
of $N$ $rel$ $\partial$) and there is a sequence of compact neighborhoods of
the repositioned $X,$ $N=N_{0}\supset N_{1}\supset\cdots\supset N_{p}$, such
that the $(p+1)$-tuple of pairs $(H_{\alpha},\delta H_{\alpha})\cap
(N_{0},N_{1},\ldots,N_{p})$ is homeomorphic, respecting this filtered
structure, to the first $p$ stages of some model ramified pair $(H^{m}%
(0),\delta H^{m}(0))$, i.e., to}
\[
((H^{m}(0),H^{m}(1),\cdots,H^{m}(p)),(\delta H^{m}(0),\delta H^{m}%
(1),\cdots,\delta H^{m}(p))).
\]

\end{lemma}

\indent We are really only interested in the $p^{th}$ stage itself, but the
intermediate stages help clarify the picture.

\indent The Repositioning Lemma will be deduced by repeated application of the
following Lemma, which was inspired by a similar construction in \v{S}tanko's
work \cite[\S \ 3.1]{St1}, \cite[\S \ 6.1]{St2}, cf. \cite[Fundamental
Lemma]{Ed1}.

\begin{lemma}
\label{Stanko-type-lemma}Suppose $N$ is a compact manifold neighborhood of a
cell-like set $X$, $\dim N\geq5,$ and suppose $D_{1}^{2},\cdots,D_{q}^{2}$ are
disjoint 2-discs in $N$ such that for each $j,$ $D_{j}^{2}\cap\partial
N=\partial D_{j}^{2}$ (with all embeddings nice, e.g. $PL$). Then for each
$j$, there is a disc-with-handles $F_{j}^{2}$, gotten from $D_{j}^{2}$ by
adding some (unspecified) number of standard handles, which are arbitrarily
small and arbitrarily close to $D_{j}^{2}\cap X$, and there is a repositioning
of $X$, with support arbitrarily close to $X$, such that when done $X\cap(%
{\textstyle\bigcup\nolimits_{j=1}^{q}}
F_{j}^{2})=\emptyset$.
\end{lemma}

\begin{proof}
[Proof of Lemma \ref{Stanko-type-lemma}]For simplicity, we treat only the
$q=1$ case, writing $D^{2}$ for $D_{1}^{2};$ the general case is the same.
Also, we omit the precise treatment of epsilons. Let $A^{2}\subset
\operatorname*{int}D^{2}$ be a small compact $2$-manifold neighborhood of
$X\cap D^{2}$ in $D^{2}$. The key fact is: for each component $A_{\ell}^{2}$
of $A^{2}$, there is a connected oriented surface $\widetilde{A}_{\ell}^{2}$
in $N-X$, lying close to $X$, with $\partial\widetilde{A}_{\ell}^{2}=\partial
A_{\ell}^{2}$. This is because, by Alexander duality, loops in $N-X$ near $X$
are null-homologous in $N-X$ near $X$. Let $\widetilde{F}^{2}=(D^{2}%
-A^{2})\cup\widetilde{A}^{2}$, where $\widetilde{A}^{2}$ is the union of the
$\widetilde{A}_{\ell}^{2}$'s, general-positioned to be disjoint as necessary.
It remains to isotope $\widetilde{F}^{2}$ to coincide with the
\textquotedblleft standard\textquotedblright\ surface $F^{2}$, which is gotten
from $D^{2}$ by adding to each component of $A_{\ell}^{2}$ a number of small
handles equal to the genus of $\widetilde{A}_{\ell}^{2}$. This is where $X$ is
moved. Since $X$ is cell-like, there is a homotopy of $\widetilde{F}^{2}$
inside of a small neighborhood of $X$, carrying id: $\widetilde{F}%
^{2}\rightarrow\widetilde{F}^{2}$ to a homeomorphism $\widetilde{F}%
^{2}\rightarrow F^{2}$. Since $\dim M\geq5$, this homotopy can be converted
into an isotopy (if $\dim N=5$, use general position and the Whitney trick to
embed this 3-dimensional homotopy in $N\times I$ near $X\times I$, and then
invoke concordance $\Rightarrow$ isotopy). This completes the proof of Lemma 3.
\end{proof}

\begin{proof}
[Proof of Lemma \ref{repositioning} from Lemma \ref{Stanko-type-lemma}]First
we produce $N_{1}$. Let $D_{\alpha}^{2}$ be the cocore $2$-disc of the handle
$H_{\alpha}$, i.e., $D_{\alpha}^{2}=0\times D^{2}\subset D_{\alpha}%
^{m-2}\times D^{2}=H_{\alpha}^{m}$. Applying Lemma 3, we can assume that $X$
has been repositioned so that $X\cap F_{\alpha}^{2}=\emptyset$, where
$F_{\alpha}^{2}$ is gotten from $D_{\alpha}^{2}$ by adding some number of
small handles, say $k$ of them. Let $H^{m}(1)$ denote the first stage of a
model $(m-2)$-handle $H^{m}(0)$ in which the 0$^{\text{th}}$ stage
ramification index is $k$, that is, $H^{m}(1)$ consists of $k$ pairs of linked
handles, instead of just one pair of linked handles as in the unramified model
(caution: by our convention, handles from distinct pairs in $H^{m}(1)$ do not
link; see Figure II-3, third frame.) Let $N_{1}$ be gotten from $N$ by
replacing $H_{\alpha}^{m}$ with a copy of $H^{m}(1)$, positioned so that
$N_{1}\cap F_{\alpha}^{2}=\emptyset$, and so that $N_{1}$ is a spine of
$N-F_{\alpha}^{2}$. This is best seen in the $m=3$ case; see Figure II-3,
third and sixth $f$. To complete this step, isotope $X$ into
$\operatorname*{int}N_{1}$.

Next, $N_{2}$ is produced inside of $N_{1}$, just as $N_{1}$ was produced in
$N$, this time using the $2k$ $2$-discs in $N_{1}$ which are the natural
cocore $2$-discs of the $2k$ $(m-2)$-handles of $N_{1}$ lying in $H_{\alpha}$.
See Figure II-3, last four frames. In applying Lemma 3 to the set of $2$-discs
in $N_{1}$, each $2$-disc may have a different number of handles added, and so
the ramification indices of these various first stage handles may be different
(of course, one could add dummy handles, to make the indices all the same).
Continuing this way one completes the proof of the Repositioning Lemma.

From this point on, the outline of the proof is the same as that of Part I,
and we will concentrate only on the nontrivial differences.
\end{proof}

\begin{proof}
[Proof of Lemma \ref{building-a-window}]This is modeled on the Basic Lemma of
Part I. We will describe the analogous successive homeomorphisms $h_{1}%
,h_{2},h_{3},h_{4}$, and finally the homeomorphism $h_{\alpha}$, without
precisely stating the full list of properties each has, to save unnecessary repetition.

During the course of this construction, there will be a great deal of
\textquotedblleft repositioning\textquotedblright\ of $X$ in $N$, often
unnamed, as done earlier in Lemmas 2 and 3, and this repositioning is to be
built into the $h_{i}$'s. This repositioning is always to be thought of as
being done in the source copy of $N\times\mathbb{R}^{1}$, that is, such
motions will always be put in front of any already-constructed motions. It
will be understood that any such repositionings in $N\times\mathbb{R}^{1}$
will be level-preserving, and will be damped to be the identity out near
$N\times\pm\infty$, away from all of the essential activity.

The goal of the first three steps is to produce a homeomorphism $h_{3}$ which
has nicely controlled behavior at the target level $N\times0$, as before. But
unlike in Part I, there will be no discussion until Step 4 of the
source-isolation of the components of $h_{i}(X\times\mathbb{R}^{1})\cap
N\times0$.

\textbf{Step 1.} Suppose $X$ has been repositioned so that $N_{1}\subset N$
exists, as described in Lemma 2. For the moment, suppose the intersection
$N_{1}\cap H_{\alpha}$ is a single, unramified pair of linked handles. The
purpose of $h_{1}$ is to \textquotedblleft pierce\textquotedblright\ the
intersection $N_{1}\times\mathbb{R}^{1}\cap N\times0=N_{1}\times0$ by moving
$N_{1}\times\mathbb{R}^{1}$ so that $h_{1}(N_{1}\times\mathbb{R}^{1})\cap
N\times0$ consists of two components, one of them a pierced copy of $N_{1}$,
and the other a thickened $(m-2)$-sphere. This is most easily described in the
source, using a polyhedral pair $(A_{1},B_{1})$, as in Part I. This time the
important $2$-dimensional part of $A_{1}$ (i.e., $A_{1}-N\times(-\infty,0])$
is more like a $2$-dimensional finger rather than a 1-dimensional feeler with
a $2$-disc attached to its end. See Figure II-4.

\begin{figure}[th]
\centerline{
\includegraphics{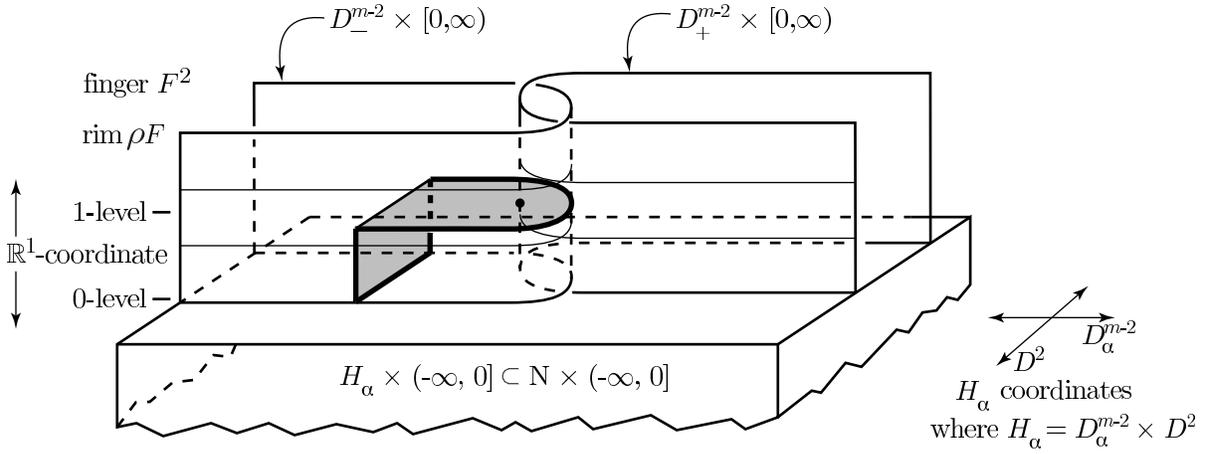}
}\caption{The 2-dimensional finger part $F^{2}$ of $A_{1}$, in the
construction of $h_{1}$}%
\end{figure}

If we were restricting the construction to the essential part of the linked
handle pair $N_{1}\cap H_{\alpha}$, namely the $(m-2)$-dimensional cores
$D_{-}^{m-2}\cup D_{+}^{m-2}\subset N_{1}\cap H_{\alpha}$, then we would
choose $A_{1}\equiv N\times(-\infty,0]\cup F^{2}$, where $F^{2}$ is a
$2$-dimensional finger (topologically a $2$-disc) which intersects
$N\times(-\infty,0]$ in (the bottom) half of its 1-dimensional circle
boundary, and intersects $D_{-}^{m-2}\times\lbrack0,\infty)$ in the other half
of its boundary $(\approx$ interval, called its \emph{rim} and denoted
$\rho_{-}^{F}$), and is such that $F^{2}\cap(D_{+}^{m-2}\times\mathbb{R}%
^{1})=$ point. If the cores $D_{-}^{m-2}\cup D_{+}^{m-2}$ are thickened by
crossing them with $D^{2}$ to form $N_{1}\cap H_{\alpha}$, then $F^{2}$ should
be thickened at its rim $\rho F$ by the same amount Given this new, suitably
thickened $A_{1}=N\times(-\infty,0]\cup F^{2}\cup\rho F\times D^{2}$, let
$h_{1}:N\times\mathbb{R}^{1}\rightarrow N\times\mathbb{R}^{1}$ be defined as
in Part I, so that $h_{1}^{-1}(N\times(-\infty,0])$ is a regular neighborhood
rel $\partial N\times\mathbb{R}^{1}$ of $A_{1}$ in $N\times\mathbb{R}^{1}$,
and so that $N_{1}\times\mathbb{R}^{1}\cap h_{1}^{-1}(N\times0)$ is the
frontier in $N_{1}\times\mathbb{R}^{1}$ of a regular neighborhood rel
$\partial N_{1}\times\mathbb{R}^{1}$ of
\[
B_{1}\equiv A_{1}\cap N_{1}\times\mathbb{R}^{1}=N_{1}\times(-\infty,0]\cup\rho
F\times D^{2}\cup(\operatorname*{point}\times D^{2}).
\]
It follows that this set $N_{1}\times\mathbb{R}^{1}\cap h_{1}^{-1}(N\times0)$
consists of the two components being sought, namely a pierced copy of $N_{1}$
and a thickened $(m-2)$-sphere.

In the general ramified situation where the intersection $N_{1}\cap H_{\alpha
}$ consists of $k$ pairs of linked handles instead of just one pair, one must
use $k$ fingers in $A_{1}$ instead of one, so that $h_{1}(N_{1}\times
\mathbb{R}^{1})\cap N\times0$ will have $k+1$ components, one of them a
$k$-times-pierced copy of $N_{1}$, denoted $N_{1\#}$, and the other $k$ of
them being thickened $(m-2)$-spheres.

To complete Step 1, one follows the above homeomorphism $h_{1}$ by a
level-preserving homeomorphism of (the target) $N\times\mathbb{R}^{1}$, which
moves $N_{1\#}\times0$ in $N\times0$ to be disjoint from the handle
$H_{\alpha}\times0$. It is the newly made holes in $N_{1\#}$ that allow this
move, just as in the pierced duncehat case in Part I.

\textbf{Step 2.} The homeomorphism $h_{2}$ will be gotten from $h_{1}$ by
preceding $h_{1}$ by a level-preserving repositioning of $X$ in the source.
Qualitatively this will be a mildly different point of view form that in Part
I, but it will accomplish the same thing (Step 2 of Part I could have been
done this way, but it didn't seem worth the effort in that simple situation).

Suppose $p=p(\delta)$ has been chosen, as it will be in Step 3 (caution: this
$p$ will correspond to $p-1$, not $p$, of Part I. Relabeling could be done to
avoid this, but that probably would cause more confusion that it would
prevent.) The motion needed to get $h_{2}$ from $h_{1}$ is simply this: $X$ is
repositioned in $N_{1}$ so that a sequence $N\supset N_{1}\supset\ldots\supset
N_{p}$ of neighborhoods of $X$ exists, as described in Lemma 2 (where $N_{1}$
is from the preceding Step 1). Then it turns out that the intersection
$h_{2}(N_{p}\times\mathbb{R}^{1})\cap N\times0$, which lies in the previously
constructed intersection $h_{1}(N_{1}\times\mathbb{R}^{1})\cap N\times0$,
consists of a copy of $N_{p}$ that has been pierced many times, plus a
collection of thickened $(m-2)$-spheres which is equivalent to the collection
of $(p-1)^{\text{st}}$ stage thickened $(m-2)$-spheres in a ramified spun Bing
collection of $(m-2)$-spheres. If there were no ramifications done in
constructing the sequence $N_{1},\ldots,N_{p}$, then there would be exactly
$2^{p-1}$ thickened $(m-2)$-spheres in this collection; in general, there are
many more.

\textbf{Step 3.} The homeomorphism $h_{3}$ is gotten by following $h_{2}$ by a
level-preserving homeomorphism of (the target) $N\times\mathbb{R}^{1}$, which
shrinks to diameter $<\delta$ the thickened $(m-2)$-sphere components of the
intersection $h_{2}(N_{p}\times\mathbb{R}^{1})\cap N\times0$. This uses the
shrinkability of the spun Bing collection of $(m-2)$-spheres, explained in the
Appendix to Part II. Here we assume that $p=p(\delta)$ was chosen so large
that this shrinking to size less than $\delta$ is possible. At the end of this
step, it has been arranged that each component of $h_{3}(N_{p}\times
\mathbb{R}^{1})\cap N\times0$ which intersects $H_{\alpha}\times0$ has
diameter $<\delta$.

\textbf{Steps 4 and }$\alpha$\textbf{.} These steps are combined, because the
generality of the present setting makes it difficult\footnote{In fact
impossible, by the example in [McM2].} to construct a pure analogue of the
homeomorphism $h_{4}$ of the duncehat case. The complicating factors will be
explained after establishing some geometry.

The important thing to realize for this step is that near any thickened
$(m-2)$-sphere component $\Sigma_{\mu}^{m-2}\times D^{2}$ of $N_{p}%
\times\mathbb{R}^{1}\cap h_{3}^{-1}(N\times0)$ (working now in the source),
the preimage $h_{3}^{-1}(N\times0)$ looks like the boundary of a tubular
neighborhood in $N\times\mathbb{R}^{1}$ of some $2$-dimensional plane, call it
$\mathbb{E}_{\mu}^{2}\times1$ (which we will assume lies in the level
$N\times1$), where $\mathbb{E}_{\mu}^{2}$ cuts transversally across one of the
$(m-2)$-handles of $N_{p}$ in $N$, say $D_{\gamma}^{m-2}\times D^{2}$, hence
$\mathbb{E}_{\mu}^{2}\cap D_{\gamma}^{m-2}\times D^{2}=\operatorname*{point}%
\times D^{2}.$ The plane $\mathbb{E}_{\mu}^{2}$ is not a closed subset of $N$,
so when we talk about a tubular neighborhood of $\mathbb{E}_{\mu}^{2}\times1$
in $N\times\mathbb{R}^{1}$, it should be understood that we are restricting
attention to a neighborhood of $N_{p}\times\mathbb{R}^{1}$. The plane
$\mathbb{E}_{\mu}^{2}$ is not immediately apparent from the description given
so far, but it could be obtained from the finger $F^{2}$ by puckering $F^{2}$,
just as the corresponding plane $\operatorname*{int}D_{\mu}^{2}$ in Part I was
obtained by puckering the $2$-disc part of $A_{1}$. The only real difference
between the present situation and that of Part I is the thickness of the
various components of the picture, i.e., the fact that some of the components
from Part I have been producted with $D^{2}$ to obtain the components here.

Consider the pair $(N(\mathbb{E}_{\mu}^{2}\times1),\partial N(\mathbb{E}_{\mu
}^{2}\times1))\cap N_{p}\times\mathbb{R}^{1}$, i.e. the intersection of the
tubular neighborhood $N(\mathbb{E}_{\mu}^{2}\times1)$ of $\mathbb{E}_{\mu}%
^{2}\times1$ in $N\times\mathbb{R}^{1}$, and its boundary $\partial
N(\mathbb{E}_{\mu}^{2}\times1)$, with $N_{p}\times\mathbb{R}^{1}$. It can be
regarded as a pair $(\Delta_{\mu}^{m-1},\Sigma_{\mu}^{m-2}))\times D^{2}$,
where the (ball, boundary sphere) pair $(\Delta_{\mu}^{m-1},\Sigma_{\mu}%
^{m-2})$ is the intersection of the pair $(N(\mathbb{E}_{\mu}^{2}%
\times1),\partial N(\mathbb{E}_{\mu}^{2}\times1))$ with $D_{\gamma}%
^{m-2}\times\mathbb{R}^{1}$, where $D_{\gamma}^{m-2}$ is as above. As a
consequence, the intersection $N_{p}\times\mathbb{R}^{1}\cap h_{3}%
^{-1}(N\times\lbrack-\epsilon,\epsilon])$, for small $\epsilon>0$, can be
regarded as being $\Sigma_{\mu}^{m-2}\times\lbrack-\epsilon,\epsilon]\times
D^{2}$, where $\Sigma_{\mu}^{m-2}\times\lbrack-\epsilon,\epsilon]$ denotes a
small collar neighborhood of $\Sigma_{\mu}^{m-2}$ in $D_{\gamma}^{m-2}%
\times\mathbb{R}^{1}$, and where the correspondence preserves the
$[-\epsilon,\epsilon]$-coordinate. By a simple reparametrization of this
collar coordinate, combined with an expansion in the target taking
$N\times\lbrack-\epsilon,\epsilon]$ onto $N\times\lbrack-1,1]$, we can assume
without loss that $\epsilon=1$ above, which we do from now on.

As in Part I, the components of the intersection $N_{p}\times\mathbb{R}%
^{1}\cap h_{3}^{-1}(N\times\lbrack-1,1])$ can be made (source-)isolated, i.e.,
their projections to the $\mathbb{R}^{1}$ coordinate can be arranged to be
disjoint, because the various disjoint thickened $(m-1)$-cells $\{\Delta_{\mu
}^{m-1}\times D^{2}\}$ corresponding to the spherical components of
intersection can be slid vertically to have nonoverlapping $\mathbb{R}^{1}%
$-coordinate values. Another way of saying this is that the planes
$\{\mathbb{E}_{\mu}^{2}\times1\}$ in $N\times\mathbb{R}^{1}$ can be vertically
repositioned to lie at different levels $\{\mathbb{E}_{\mu}^{2}\times t_{\mu
}\}$. We assume that this has been done (as it was in Part I).

In constructing the homeomorphism $h_{\alpha}$ there are two aspects of this
general situation which serve to make this step more complicated than that of
Part I: the $2$-plane $\mathbb{E}_{\mu}^{2}$ may intersect $X$ in more than a
single point, and also $X\times t_{\mu}$ may not have a ball regular
neighborhood in $N\times\mathbb{R}^{1}$ whose boundary slices $X\times
\mathbb{R}^{1}$ at precise, entire levels (this was guaranteed in Part I by
$K\times I$ being collapsible). Before dealing with these difficulties, it is
worth noting that if in fact these nice conditions prevail, then the
construction in Part I for $h_{4}$ and $h_{\ast}$ also works here to produce
$h_{\alpha}$.

The following adaptation of the Part I construction will take care of both
difficulties at the same time. Suppose for the moment that $C_{\mu}%
\equiv\mathbb{E}_{\mu}^{2}\cap X$ is $0$-dimensional, e.g., a cantor set. Let
$\phi_{\mu}:N\times\mathbb{R}^{1}\longrightarrow N\times\mathbb{R}^{1}$ be a
near-homeomorphism (i.e. a limit of homeomorphisms), supported arbitrarily
near $X\times t_{\mu}$, such that for each $t\in\mathbb{R}^{1}$, if $\phi
_{\mu}(X\times t)\cap\mathbb{E}_{\mu}^{2}\times t_{\mu}\neq\emptyset$, then
$\phi_{\mu}(X\times t)=\operatorname*{point}\in\mathbb{E}_{\mu}^{2}\times
t_{\mu}$, and furthermore these $X\times t$'s are the only nontrivial
point-inverses of $\phi_{\mu}$ (see Figure II-5). To get the map $\phi_{\mu}$,
basically one uses the idea that any cantor set's worth of $X\times t$'s in
$N\times\mathbb{R}^{1}$ is shrinkable. Arguing more precisely, one first can
do a vertical perturbation of $\mathbb{E}_{\mu}^{2}\times t_{\mu}$ to make it
intersect each $X\times t$ in at most a single point, and then one can use
engulfing to shrink these $X\times t$'s to the points of intersection, keeping
fixed the repositioned $\mathbb{E}_{\mu}^{2}\times t_{\mu}$. After shrinking,
one brings the repositioned $\mathbb{E}_{\mu}^{2}\times t_{\mu}$ back to its
original position, using the inverse of the original vertical perturbation.

\begin{figure}[th]
\centerline{
\includegraphics{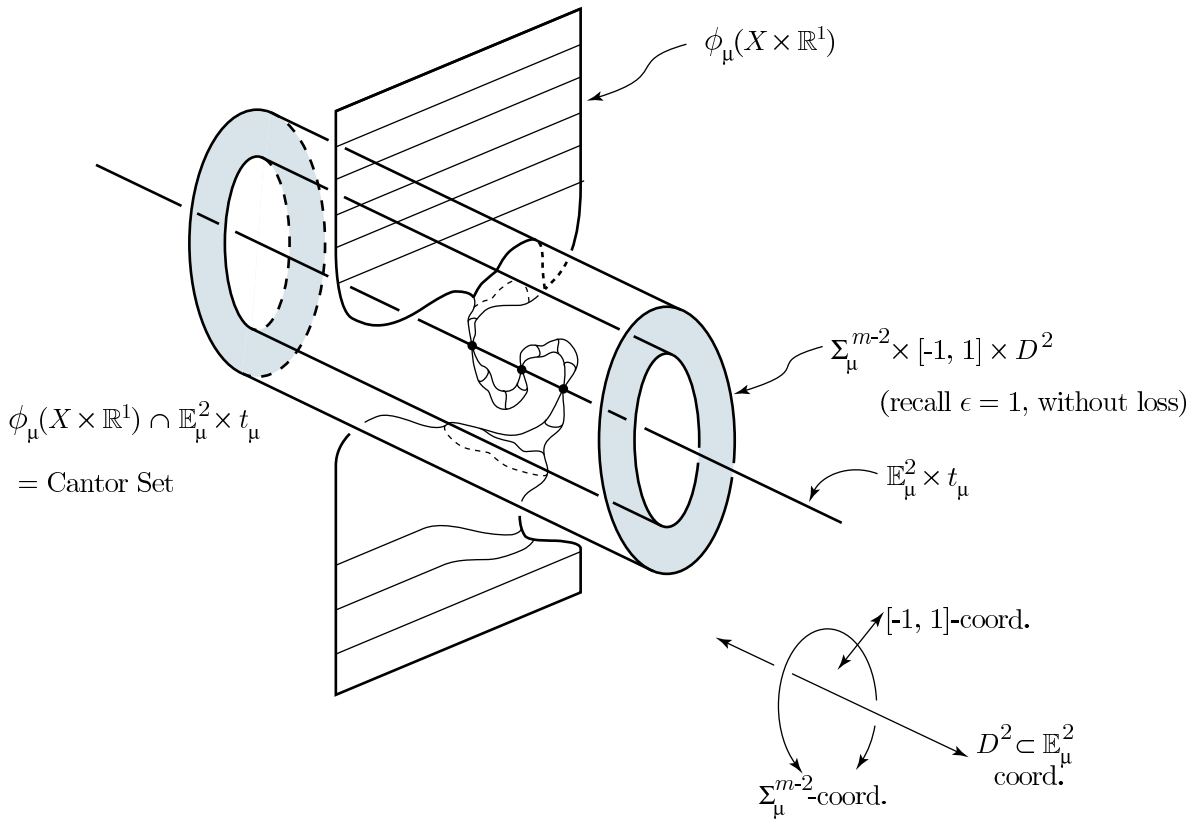}
}\caption{The construction of $h_{4}$ and $h_{\alpha}$}%
\end{figure}

Given $\phi_{\mu}$, the idea now is to do a sort of expansion-meshing
operation, by taking the earlier chosen collar $\Sigma_{\mu}^{m-2}%
\times\lbrack-1,1]$ in $D_{\gamma}^{m-2}\times\mathbb{R}^{1}$ (recall
$\epsilon=1$) and isotoping it in $D_{\gamma}^{m-2}\times\mathbb{R}^{1}$ in a
certain manner. Because of the correspondence made above of $\Sigma_{\mu
}^{m-2}\times\lbrack-1,1]\times D^{2}$ with $N_{p}\times\mathbb{R}^{1}\cap
h_{3}^{-1}(N\times\lbrack-1,1])$, this isotoping may be regarded as producing
a modification of the homeomorphism $h_{3}$. The goal of isotoping the collar
is to achieve, for an arbitrary preassigned $\eta>0$, that for each
$t\in\mathbb{R}^{1}$, $\phi_{\mu}(X\times t)\cap\Sigma_{\mu}^{m-2}%
\times\lbrack-1,1]\times D^{2}\subset\Sigma_{\mu}^{m-2}\times\lbrack
s-\eta,s+n]\times D^{2}$ for some $s\in\lbrack-1,1]$. This is done as follows.
First, one takes the originally chosen collar $\Sigma_{\mu}^{m-2}\times
\lbrack-1,1]$ and, keeping the outer boundary $\Sigma_{\mu}^{m-2}\times1$
fixed, one isotopes the band $\Sigma_{\mu}^{m-2}\times\lbrack-1,1-\eta]$ so
close to the plane $\mathbb{E}_{\mu}^{2}\times t_{\mu}$ that no image
$\phi_{\mu}(X\times t)$, $t$ arbitrary, intersects both $\Sigma_{\mu}%
^{m-2}\times1$ and $\Sigma_{\mu}^{m-2}\times\lbrack-1,1-\eta]$. Next, keeping
fixed the repositioned band $\Sigma_{\mu}^{m-2}\times\lbrack1-\eta,1]$ one
isotopes the band $\Sigma_{\mu}^{m-2}\times\lbrack-1,1-2\eta]$ much closer to
the plane $\mathbb{E}_{\eta}^{2}\times t_{\mu}$, so that no image $\phi_{\mu
}(X\times t)$ intersects both $\Sigma_{\mu}^{m-2}\times\lbrack1-\eta,1]$ and
$\Sigma_{\mu}^{m-2}\times\lbrack-1,1-2\eta]$. Continuing this way, the desired
degree of control is achieved.

Having done this expansion-meshing operation, independently and disjointly for
each $\mu$, then one can define $h_{\alpha}=\overline{h}_{3}\overline{\phi}$,
where $\overline{\phi}$ is a homeomorphism closely approximating the union
$\phi$ of all of the $\phi_{\mu}$'s that were chosen above, and where
$\overline{h}_{3}$ denotes the homeomorphism obtained by modifying $h_{3}$ by
the expansion-meshing operations just described.

It remains to discuss how the 0-dimensionality of $C_{\mu}$ can be achieved.
If $X$ had dimension $\leq m-2$, i.e. if $X$ had codimension $\geq2$ in $N$,
then achieving this would be an easy matter of putting $X$ in topological
general position with respect to $\mathbb{E}_{\mu}^{2}$. This suggests a
solution: make $X$ have dimension $\leq m-2$. That is, replace $X$ with a new
cell-like space $X_{\ast}$ lying near $X$ and obtained from $X$ by a certain
limiting process, such that $\dim X_{\ast}\leq m-2$ and $M-X_{\ast}\approx
M-X$. The most convenient time to do this is at the very start of the proof,
using the fact that $X$ has arbitrarily small neighborhoods with spines of
dimension $\leq m-2$. For if one takes a nested basis of such neighborhoods,
and squeezes the first one close to its spine, and then the (repositioned)
second one very close to its (repositioned) spine, etc., one produces in the
limit the desired $X_{\ast}$. This process is enlarged upon in the Postscript
below, where a replacement $X_{\ast}$ having much nicer properties is
produced. This completes the proof of Lemma 1, in the $m\geq5$ case.

In the $m=4$ case, the above proof breaks down for one essential reason: Lemma
\ref{repositioning} is unknown. The trouble comes in Lemma
\ref{Stanko-type-lemma}, in constructing the surfaces $\widetilde{A}_{\ell
}^{2}$ (which can in fact be done), and in trying to move these surfaces back
to standard position by isotopy (this is a fundamental problem). The way to
circumvent this difficulty is to use the extra freedom provided by the
$\mathbb{R}^{1}$ coordinate, which in effect turns the problem from a
4-dimensional problem into a 5-dimensional problem. An outline of this rescue
operation follows.

The idea is that one can at least do the already-described motions of Lemma
\ref{Stanko-type-lemma} in $N\times\mathbb{R}^{1}$, if not in $N$. This lets
one prove a weaker version of Lemma \ref{repositioning}, which says that
$X\times\mathbb{R}^{1}$ can be repositioned in $N\times\mathbb{R}^{1}$ so that
$X\times\mathbb{R}^{1}\cap N\times\lbrack-a,a]\subset N_{p}\times\lbrack-a,a]$
for any preassigned large $a$ (of course, the vertical movement of
$X\times\mathbb{R}^{1}$ may have to be as large as $a$, but that is not
important). Now, one can do Steps 1, 2, and 3 above without change. Step
4-$\alpha$ works also (let us assume, as justified above, that $\dim X\leq2)$,
even though $X\times\mathbb{R}^{1}$ has been grossly perturbed. One can still
make the intersection of $\mathbb{E}_{\mu}^{2}\times t_{\mu}$ with the
perturbed $X\times\mathbb{R}^{1}$ be 0-dimensional, with the intersection
points all having distinct $\mathbb{R}^{1}$-coordinates in $X\times
\mathbb{R}^{1}$, and then one can shrink these particular $X\times t$ levels
to these intersection points. Then the expansion-meshing works as described,
to produce $h_{\alpha}$.

The proof of Lemma 1 is now complete, for all $m\geq4$, as claimed.
\end{proof}

\indent From this point on, the proof is modeled on \cite[pp. 201,202]{EM},
just as the proof in Part I. The first step is to establish the\bigskip

\textsc{Window Building Lemma. }\emph{Suppose the data }$X\subset N^{m}$\emph{
as earlier (explained before Lemma 1), and let }$H_{\#}=\cup_{\alpha=1}%
^{r}H_{\alpha}$\emph{ be the union of the }$(m-2)$\emph{-handles of }%
$N$\emph{. Then given any }$\delta>0$\emph{, there is a homeomorphism }%
$h_{\#}:N\times\mathbb{R}^{1}\rightarrow N\times\mathbb{R}^{1}$\emph{, fixed
on }$\partial N\times\mathbb{R}^{1}$\emph{, such that}

\begin{itemize}
\item[1.] \emph{for each }$j\in2\mathbb{Z}$\emph{ and each }$t\in\lbrack
j-1,j+1],$ $h_{\#}(N\times t)\subset N\times\lbrack j-1,j+1]$\emph{, and}

\item[2.] \emph{for each }$t\in\mathbb{R}^{1},$\emph{ if }$h_{\#}(X\times
t)\cap H_{\#}\times\lbrack j-1+\delta,j+1-\delta]\neq\emptyset$\emph{ for any
}$j\in2\mathbb{Z}$\emph{, then }$\operatorname*{diam}h_{\#}(X\times t)<\delta
$\emph{.\bigskip}
\end{itemize}

\indent This is easily proved from the model single window version, Lemma 1.
The idea is first to construct a homeomorphism $h_{\ast}:N\times\mathbb{R}%
^{1}\rightarrow N\times\mathbb{R}^{1}$, compactly supported in
$\operatorname*{int}N\times\mathbb{R}^{1}$, such that $h_{\ast}$ satisfies
condition (2) above for the value $j=0$. This $h_{\ast}$ is gotten by first
selecting some nonoverlapping vertical translates of the \textquotedblleft
window blocks\textquotedblright\ $H_{\alpha}\times\lbrack-1+\delta
,1-\delta],1\leq\alpha\leq r$ then applying Lemma 1 separately to
\textquotedblleft make a window\textquotedblright\ in each of these completely
disjoint blocks, using for this the different $h_{\alpha}$'s, and finally
translating these windows back to their original positions. Next, $h_{\ast}$
can be conjugated by a vertical homeomorphism of $N\times\mathbb{R}^{1}$ to
arrange that the support of $h_{\ast}$ lies in $\operatorname*{int}%
N\times(-1,1)$. From this $h_{\ast}$, one builds the desired $h_{\#}$ by
stacking vertical translates of $h_{\ast}$ on top of each other.

\begin{proof}
[Proof of Shrinking Proposition from the Window Building Lemma]Given a
neighborhood $U$ of $X$ in $M^{m}$ and an $\epsilon>0$, choose a small ball
$B^{m}$ in $U$, with $\operatorname*{diam}B^{m}<\epsilon/2$ and
$\operatorname*{int}B^{m}\cap X\neq\emptyset$. We show how to construct a
homeomorphism $h:M\times\mathbb{R}^{1}\rightarrow M\times\mathbb{R}^{1}$,
fixed on $(M-U)\times\mathbb{R}^{1}$, which satisfies the following weakened
versions of the conditions from the Shrinking Proposition: for each
$t\in\mathbb{R}^{1}$,

\begin{itemize}
\item[1$^{\prime}$.] $h(U\times t)\subset U\times\lbrack t-3,t+3]$, and

\item[2$^{\prime}$.] either

\begin{itemize}
\item[a.] $h(X\times t)\subset B^{m}\times\lbrack t-3,t+3]$, or

\item[b.] $\operatorname*{diam}h(X\times t)<\epsilon$.
\end{itemize}
\end{itemize}

From this it is clear that the original Shrinking Proposition follows, simply
by rescaling the vertical coordinate.

To construct this $h$, first one constructs a certain auxiliary homeomorphism
$g=g_{2}g_{1}:M\times\mathbb{R}^{1}\rightarrow M\times\mathbb{R}^{1}$. Let
$\ast\in\operatorname*{int}B^{m}\cap X$ be a basepoint. First one chooses a
uniformly continuous homeomorphism $g_{2}$, which is supported in
$\bigcup\{U\times\lbrack j+\frac{1}{2},j+\frac{3}{2}]\mid j\in2\mathbb{Z}\}$
and fixed on $\ast\times\mathbb{R}^{1}$, such that for each $j\in
2\mathbb{Z},g_{2}(X\times(j+1))\subset\operatorname*{int}B^{m}\times
\mathbb{R}^{1}$. The existence of $g_{2}$ follows from the cellularity of
$X\times t$ in $M\times\mathbb{R}^{1}$. Let $V$ be a compact neighborhood of
$X$ in $U$ such that $g_{2}(V\times(j+1))\subset\operatorname*{int}B^{m}%
\times\mathbb{R}^{1}$ for each $j\in2\mathbb{Z}$, and let $C^{m}$ be a compact
neighborhood of $\ast$ in $B^{m}$ such that $g_{2}(C^{m}\times\mathbb{R}%
^{1})\subset\operatorname*{int}B^{m}\times\mathbb{R}^{1}$. By codimension 3
engulfing, there is a (small) neighborhood $N$ of $X$ in $V$, with structure
$N=L\cup\bigcup_{\alpha=1}^{r}H_{\alpha}$ as discussed earlier (so in
particular $L$ has an $(m-3)$-spine), and there is a homeomorphism
$g_{0}:M\rightarrow M$, supported in $V$, such that $g_{0}(L)\subset C$. Let
$g_{1}=g_{0}\times id_{\mathbb{R}^{1}}:M\times\mathbb{R}^{1}\rightarrow
M\times\mathbb{R}^{1}$. Then $g=g_{2}g_{1}:M\times\mathbb{R}^{1}\rightarrow
M\times\mathbb{R}^{1}$ is such that

\begin{itemize}
\item[1.] $g=$ identity on $(M-U)\times\mathbb{R}^{1}$, and for each
$t\in\mathbb{R}^{1}$, $g(U\times t)\subset U\times\lbrack t-1,t+1],$ and

\item[2.] the image under $g$ of $L\times\mathbb{R}^{1}\cup\bigcup
\{N\times(j+1)\mid j\in2\mathbb{Z}\}$ lies in $B^{m}\times\mathbb{R}^{1}$.
\end{itemize}

Given $g$, then one can let $h=gh_{\#}$, where $h_{\#}$, is provided by the
Window Building Lemma for some sufficiently small value of $\delta
=\delta(\epsilon,g)$, and where $h_{\#}$ is assumed to be extended via the
identity over $(M-N)\times\mathbb{R}^{1}$. This completes the proof of the
Shrinking Proposition, and hence Theorem 2.

\newpage

\centerline{\bf Postscript to Part II}\bigskip
\end{proof}

\indent Embedded in the preceding proof is what I call the\bigskip

\textsc{Replacement Principle for Cell-Like Sets. }\emph{Given any cell-like
compact set }$X$\emph{ in a manifold-without-boundary }$M,$ $\dim M\geq
5$\emph{, then }$X$\emph{ can be replaced by a nearby 1-dimensional cell-like
set }$X_{\ast}$\emph{, such that }$M-X_{\ast}$\emph{ is homeomorphic to
}$M-X,$ \emph{while }$X_{\ast}$\emph{ is homeomorphic to a reduced cone on a
cantor set (defined more precisely below).}

\emph{In fact, the replacement operation can be done in the following
continuous manner. There is an isotopy of embeddings }$g_{t}:M-X\rightarrow
M,$ $0\leq t\leq1$\emph{, starting at }$g_{0}=identity$\emph{ and fixed off of
an arbitrarily small neighborhood of }$X$\emph{, such that }$g_{1}%
(M-X)=M-X_{\ast}$\emph{. Thus, in particular, }$X$\emph{ is transformed to
}$X_{\ast}$\emph{ in a semicontinuous fashion, through the intermediate
cell-like sets }$\{M-g_{t}(M-X)\}$\emph{.}\textit{\bigskip}

\indent I first enunciated this principle at the time of the Park City
Conference (February 1974), when I proved it for (wildly embedded) cell-like
sets $X$ of codimension $\geq3$, using my variation of the proof of
\v{S}tanko's Approximation Theorem (unpublished notes; cf. \cite{BL2}). In
certain special cases one can find better replacements for $X$, e.g. Giffen's
replacement of the duncehat spine of Mazur manifold by a $2$-disc \cite{Gi}
(which can then be made into an arc (see \cite{DE} and \cite[\S 3]{She})), but
as general rule, the above principle seems the best possible, even for $X$ a polyhedron.

\begin{proof}
[Proof of the Replacement Principle for Cell-Like Sets ](Using the machinery
of Part II). The basic tool is the following variation of Lemma
\ref{repositioning} (see Figure IIP-1). Recall $m\geq5$.\bigskip
\end{proof}

\setcounter{figure}{0} \renewcommand{\thefigure}{\Roman{part}P-\arabic{figure}}

\begin{figure}[th]
\centerline{
\includegraphics{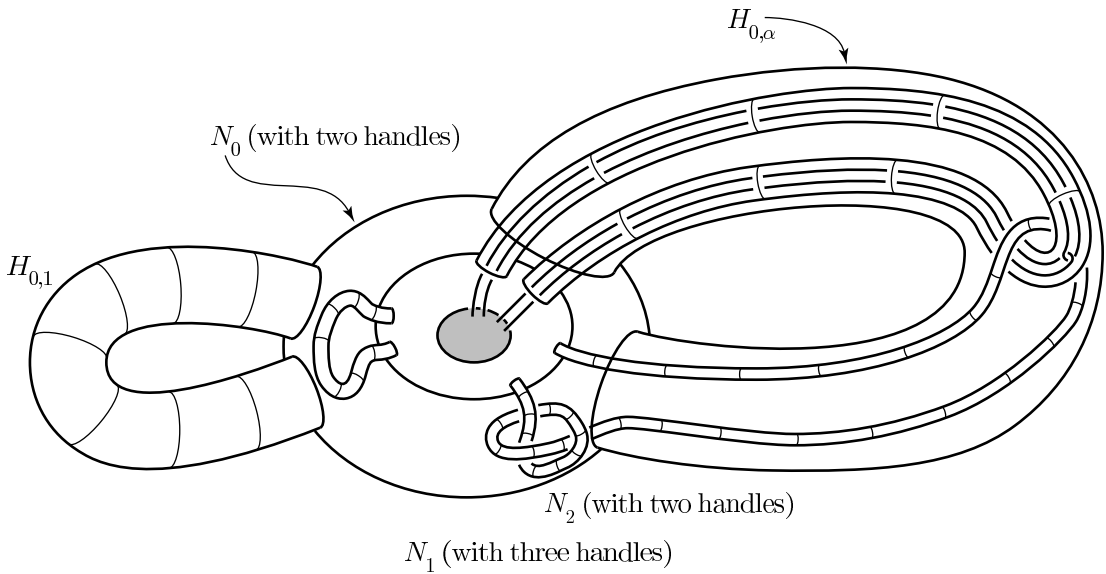}
}\caption{The nested neighborhood sequence $\{N_{i}\}$, from the Postscript.}%
\end{figure}

\textsc{Lemma 2}$^{\prime}$\textsc{. }\textit{Given any neighborhood
$N_{0}^{m}=L_{0}\cup\bigcup\{H_{0,\alpha}\mid1\leq\alpha\leq r_{0}\}$ of $X$
with structure as described earlier \textrm{(i.e.,} $L_{0}$ has an
$(m-3)$-spine, and the $H_{0,\alpha}$'s are $(m-2)$-handles attached to
$L_{0}$, then there is an arbitrarily small neighborhood $N_{1}=L_{1}%
\cup\{H_{1,\beta}\mid1\leq\beta\leq r_{1}\}$ of $X$ in $\operatorname*{int}%
N_{0}$, and a repositioning of $N_{1}$ in $N_{0}$ (by an unnamed isotopy of
$N_{0}\;rel\;\partial N_{0}$), such that after repositioning,}

\begin{itemize}
\item[1.] $L_{1}\subset\;\operatorname*{int}L_{0}$,

\item[2.] \textit{for each handle $H_{0,\alpha}$ of $N_{0}$, the triple
$(H_{0,\alpha},\delta H_{0,\alpha},H_{0,\alpha}\cap N_{1})$ is homeomorphic to
a standard triple $(H(0),\delta H(0),H(1))$, where $H(1)$ is a ramified
collection of 1st stage $(m-2)$-handles in $H(0)$ (as defined earlier; see
Figure II-3, third frame), and }

\item[3.] \textit{each handle $H_{1,\beta}$ of $N_{1}$ intersects at most one
handle $H_{0,\alpha}$ of $N_{0}$, and there is at most one component of
intersection (which by (2) is a single subhandle of $H_{0,\alpha})$.\bigskip}
\end{itemize}

\begin{proof}
[Proof of Lemma 2$^{\prime}$]First choose an arbitrarily small neighborhood
$N^{m}$ of $X$ having $(m-2)$-spine, and position $N$ in $N_{0}$, using simple
general positioning, so that for each handle $H_{0,\alpha}$ of $N_{0}$, the
intersection $N\cap H_{0,\alpha}$ goes straight through the handle, looking
like several parallel, smaller copies of $H_{0,\alpha}$ (i.e., like a ramified
version of $H_{0,\alpha}$, as in Figure II-3, second frame). Now, using the
earlier Repositioning Lemma 2, applied to $N$ and the $2$-discs which are the
central cocores of the handle-components of $N\cap H_{0,\alpha}$, for all
$\alpha$, find a repositioning of $X$ in $N$ and a neighborhood $N_{1}$ having
$(m-2)$-spine such that for each handle $H_{0,\alpha}$, the intersection
$H_{0,\alpha}\cap N_{1}$ is as described in (2) above. It remains to achieve
(1) and (3). Let the $(m-2)$-dimensional spine $K_{1}$ of $N_{1}$ be
positioned, and subdivided if necessary, so that $K_{1}^{(m-3)}\subset
\operatorname*{int}L_{0}$, and so that each $(m-2)$-simplex of $K_{1}$
intersects at most one $(m-2)$-handle $H_{0,\alpha}$ of $N_{0}$, in at most a
single $(m-2)$-disc. Now squeeze $N_{1}$ close to its newly positioned spine
$K_{1}$, and let $L_{1}$ be the \textquotedblleft
restriction\textquotedblright\ of $N_{1}$ to $K_{1}^{(m-3)}$, This completes
the proof of Lemma 2$^{\prime}$.
\end{proof}

\begin{proof}
[Proof of Replacement Principle (continued)]Given Lemma 2$^{\prime}$, one
applies it repeatedly to produce an infinite sequence $\{N_{i}=L_{i}%
\cup\bigcup\{H_{i,\gamma}\mid1\leq\gamma\leq r_{i}\}\}$ of neighborhoods of
$X$, repositioning $X$ at each step, such that each $N_{i}$ is null homotopic
in its predecessor $N_{i-1}$, and is positioned there as described by Lemma
2'. This produces in the limit a new cell-like compactum $X_{\infty}%
=\cap\{N_{i}\}$. We assume in addition that in the sequence $\{N_{i}\}$ the
handle alignment is maintained from step to step as shown in Figure IIP-1 so
that for any handle $H_{i,\gamma}$ of any neighborhood $N_{i}$, the triple
$(H_{i,\gamma},\delta H_{i,\gamma},H_{i,\gamma}\cap X_{\infty})$ is
homeomorphic to a model ramified triple $(H(0),\delta H(0),H(\infty))$, where
$H(\infty)\equiv\cap_{p=0}^{\infty}H(p)$, which is a ramified collection of
$(m-2)$-dimensional Bing cells. Elaborating this, letting $\delta
H(\infty)\equiv\cap_{p=0}^{\infty}\delta H(p)=H(\infty)\cap\delta H(0)$, then
$\delta H(\infty)$ is a cantor set's worth of $(m-3)$-cells (namely the
$(m-4)$-times-spun Bing collection of arcs) and $(H(\infty),\delta
H(\infty))\approx\delta H(\infty)\times([0,1],0)$. Hence in particular
$X_{\infty}-\operatorname*{int}H_{i,\gamma}$ is a strong deformation retract
of $X_{\infty}$.

For each handle $H_{i,\gamma}$, let $G_{i,\gamma}$ denote the union of all
$(m-2)$-handles $H_{j,\sigma},j>1$, which intersect $H_{i,\gamma}$, and let
$F_{i,\gamma}=G_{i,\gamma}\cap X_{\infty}$. By our assumptions, if
$H_{i,\gamma}$ and $H_{k,\tau}$ are two handles and if $G_{i,\gamma}\cap
G_{k,\tau}\neq\emptyset$ then $k>i$ and $G_{k,\tau}\subset G_{i,\gamma}$ (or
vice versa). More importantly, $F_{i,\gamma}$ is homeomorphic to
$H(\infty)-\delta H(\infty)$ for some ramified model $H(\infty)$. Let
$L_{\infty}=\cap_{i=0}^{\infty}L_{i}\subset X_{\infty}$, so that $X_{\infty
}-L_{\infty}$ is the disjoint union of all of the $F_{i,\gamma}$'s.\medskip

\textbf{Claim.} $L_{\infty}$ \emph{is cellular in }$M$.

\textit{Proof.} $L_{\infty}$ is cell-like because it is the intersection of a
nested sequence of cell-like sets, namely $L_{\infty}=\cap_{i=0}^{\infty
}(L_{i}\cap X_{\infty})$. Each set $L_{i}\cap X_{\infty}$ is cell-like because
it is a retract of $X_{\infty}$, for as noted above, $X_{\infty}%
-\operatorname*{int}H_{i,\gamma}$ is a strong deformation retract of
$X_{\infty}$ for each $m-2$ handle $H_{i,\gamma}$ of $N_{i}$. Granted that
$L_{\infty}$ is cell-like, then it is cellular because it is an intersection
of manifolds having $(m-3)$-spines. This establishes the Claim.\medskip

Since $L_{\infty}$ is cellular, it can be shrink to a point, producing
$X_{\#}=X_{\infty}/L_{\infty}$. This space $X_{\#}$ has particularly nice
structure, being a countable null wedge of subsets of $M$ each homeomorphic to
$H(\infty)/\delta H(\infty)$. This latter set is a cone on a cantor set's
worth of $m-3$ cells (as in Bing's original $m=3$ case). To obtain $X_{\ast}$
from $X_{\#}$, one can use the shrinkability of the spun Bing collection of
$(m-3)$-cells (see the Appendix) to shrink each copy of $H(\infty)/\delta
H(\infty)$ to a cone on a cantor set, by shrinking to a cantor set each of the
parallel copies of $\delta H(\infty)\times\{\operatorname*{point}\}$ that make
up $H(\infty)/\delta H(\infty)\approx\delta H(\infty)\times\lbrack0,1]/\delta
H(\infty)\times0=c\delta H(\infty)$.\newpage

\centerline{\bf Appendix.  Shrinking the spun Bing decomposition.}\bigskip
\end{proof}

\indent The $k$-times spun Bing decomposition of euclidean $(3+k)$-space into
points and tame $(k+1)$-discs, described earlier and again below, is at the
very heart of the work of Parts I and II (and Cannon's subsequent work, too).
The purpose of this Appendix is to establish:

\begin{theorem}
\label{spun-bing}The spun Bing decomposition is shrinkable.
\end{theorem}

\indent The ramified versions of the spun Bing decomposition, which appeared
earlier in Part II, are also shrinkable. This will be evident from the
construction below, which also works in the ramified context, with identical
motions and equal efficiency.

The $k$-times spun Bing decomposition has been around for some time (as I
found out after needing and proving the above theorem in 1974). L. Lininger,
in order to exhibit some nonstandard topological involutions of spheres,
showed in \cite{Li} that the once-spun version was shrinkable, and suggested
that the higher dimensional versions were also. Neuzil's thesis \cite{Neu}
handled a class of decompositions of 4-space that included the once-spun
version. Recently, R. Daverman \cite[Cor 11.7]{Da3} showed that the
shrinkability of the general $k$-times-spun, versions, $k\geq2$, followed from
a general mismatch theorem. This was adapted to the $k=1$ case in \cite{CD}.

\indent The proof given here is elementary, in the sense of Bing's original
proof, for that is what it is modeled on. In fact, it is almost fair to say
that one simply takes Bing's original proof and spins it.

\begin{proof}
[Proof of Theorem \ref{spun-bing}]The proof is divided into two natural steps.
The first step describes a certain decomposition of $\mathbb{R}^{3}%
\times\lbrack0,\infty)$, gotten by unfolding, or unraveling, the original Bing
decomposition of $\mathbb{R}^{3}\times0$, and then shows that this
decomposition of $\mathbb{R}^{3}\times\lbrack0,\infty)$ is shrinkable in a
certain level-preserving manner. The second step shows how to spin these
motions to shrink the spun Bing decomposition.\medskip

\textbf{Step 1. }Consider the original Bing decomposition of $\mathbb{R}^{3}$
(from \cite{Bi1}), as shown at the bottom of Figure IIA-1. It is defined using
two disjoint, linked embeddings $\chi_{-},\chi_{+}:S^{1}\times D^{2}%
\rightarrow S^{1}\times\operatorname*{int}D^{2}$, whose images are denoted
$S_{-}^{1}\times D^{2}$ and $S_{+}^{1}\times D^{2}$. By taking iterates of
these embeddings, e.g. $\chi_{+}(\chi_{-}(S^{1}\times D^{2}))$, one obtains
the deeper stages of the defining neighborhood sequence of solid tori. At the
$p^{\text{th}}$ stage, there are $2^{p}$ linked solid tori, denoted $\{S_{\mu
}^{1}\times D^{2}\mid\mu\in\{-,+\}^{p}\}$, embedded in the original solid
torus $S^{1}\times D^{2}$. Letting $T_{p}$ denote the union of the $2^{p}$
$p^{\text{th}}$ stage components, then the nontrivial elements of the Bing
decomposition of $\mathbb{R}^{3}$ are the components of the intersection
$\cap_{p=1}^{\infty}T_{p}$, which can be arranged to be $C^{\infty}$-smooth
arcs, as suggested in Figure IIA-1.

\setcounter{figure}{0} \renewcommand{\thefigure}{\Roman{part}A-\arabic{figure}}

\begin{figure}[th]
\centerline{
\includegraphics{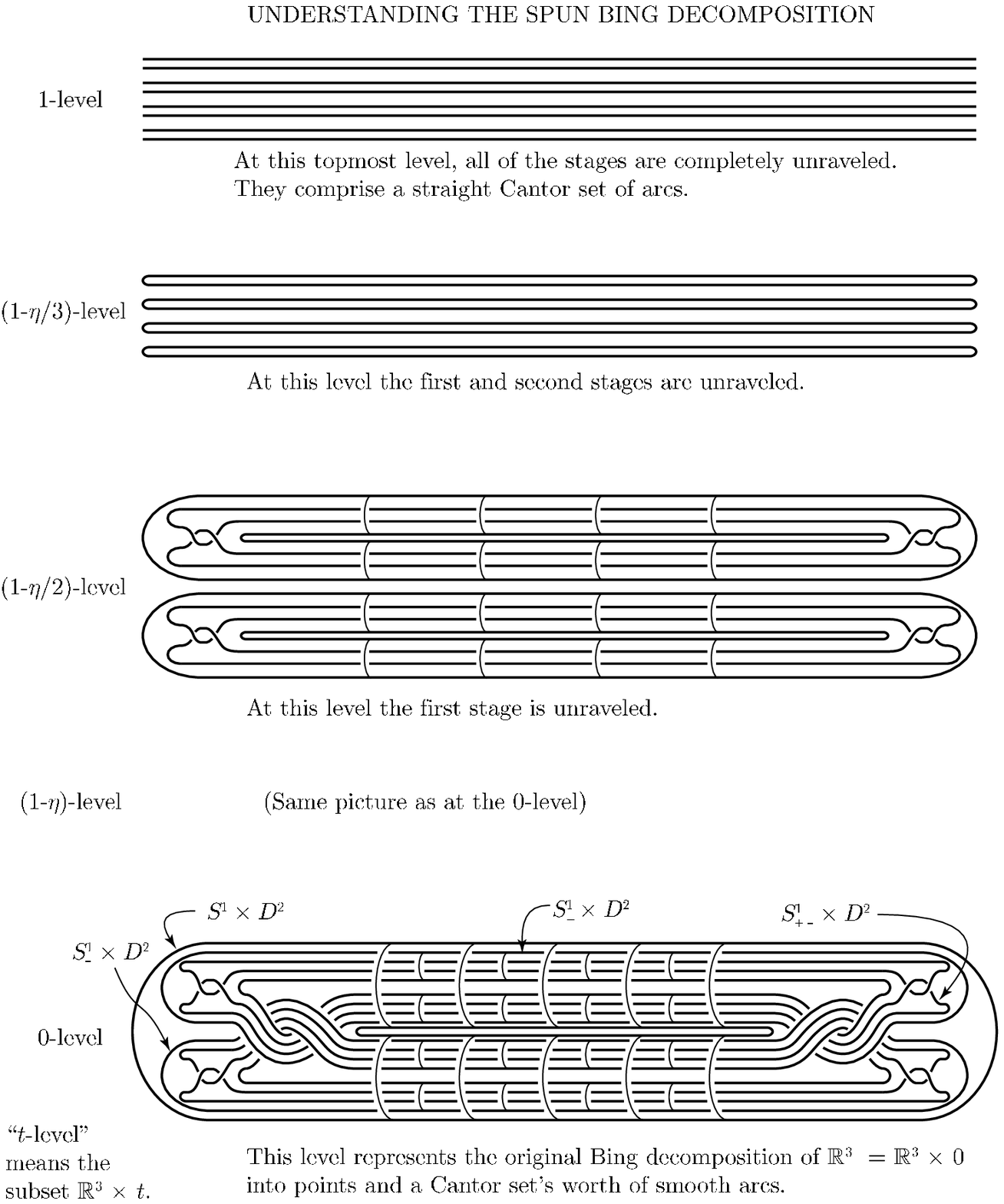}
}\caption{Unraveling the original Bing decomposition in $\mathbb{R}^{3}%
\times[0,\infty)$}%
\end{figure}

We wish to extend this decomposition of $\mathbb{R}^{3}=\mathbb{R}^{3}\times0$
to one of $\mathbb{R}^{3}\times\lbrack0,\infty)$, with all decomposition
elements lying in levels $\{\mathbb{R}^{3}\times t\}$, by unravelling it, as
suggested by Figure IIA-1. Fix some small $\eta>0$. Then, as one works up
through the levels, from $t=0$ to $t=1$, the following activity takes place:
from $t=0$ to $t=1-\eta$, nothing happens, i.e., the decomposition of each
intermediate level is just a translate of the decomposition of the 0-level;
from $t=1-\eta$ to $t=1-\eta/2$, the first stage is unraveled, i.e.,
$S_{+}^{1}\times D^{2}$ and $S_{-}^{1}\times D^{2}$ are isotoped until they
look disjoint (i.e., are separated by a hyperplane), as pictured; from
$t=1-\eta/2$ to $t=1-\eta/3$, the second stage is unraveled, etc. What one
obtains is an unraveling of the original Bing decomposition of $\mathbb{R}%
^{3}\times0$ to a standard decomposition of $\mathbb{R}^{3}\times1$ consisting
of a straight cantor set's worth of arcs. (It is interesting to ponder the
pseudoisotopy of $\mathbb{R}^{3}$ that this process yields in the limit, but
that is another matter.)

There is a natural \textquotedblleft defining neighborhood
sequence\textquotedblright\ for this decomposition of $\mathbb{R}^{3}%
\times\lbrack0,\infty)$, which we describe at this point. For each $p\geq0$
and each $\mu\in\{-,+\}^{p}$, let $C_{\mu}^{3}$ denote the convex hull of the
$p^{\text{th}}$ stage component $S_{\mu}^{1}\times D^{2}$ \textbf{after} it
has been unraveled, i.e., in any level $\mathbb{R}^{3}\times t$ for
$1-\eta/(p+1)\leq t\leq1$. Let $N_{p}^{4}$ denote the union of $\cup\{C_{\mu
}^{3}\mid|\mu|=p\}\times\lbrack1-\eta/(p+1),1+\eta/(p+1)]$ together with the
unraveling images of $\cup\{S_{\mu}^{1}\times D^{2}\mid|\mu|=p\}$ in
$\mathbb{R}^{3}\times\lbrack0,1-\eta/(p+1)]$. Then each $N_{p}^{4}$ is
homeomorphic to $B^{2}\times D^{2}$, with $N_{p}^{4}\cap\mathbb{R}^{3}\times0$
corresponding to $\partial B^{2}\times D^{2}$, and furthermore the above
described decomposition of $\mathbb{R}^{3}\times\lbrack0,\infty)$ has as its
nontrivial elements the components of $(\cap_{p=0}^{\infty}N_{p}^{4}%
)\cap\mathbb{R}^{3}\times t$, for $t$ varying from 0 to 1.

It turns out that this decomposition of $\mathbb{R}^{3}\times\lbrack0,\infty)$
is shrinkable \emph{in a level-preserving manner}. This calls for a mild
extension of Bing's original shrinking argument.

As Bing showed, the key to shrinking his decomposition is to regard
$S^{1}\times D^{2}$ as being long and thin, as drawn in Figure II A-2, and to
show that if one cuts $S^{1}\times D^{2}$ into any number of chambers by means
of transverse $2$-planes in $\mathbb{R}^{3}$, then the decomposition elements
can be isotoped in $S^{1}\times D^{2}$ so that when done each lies inside a
single chamber. To save words, we make a definition for this. Suppose
$P_{1}^{2},\ldots,P_{q}^{2}$ are $q$ parallel $2$-planes in $\mathbb{R}^{3}$,
transverse to the \textquotedblleft long\textquotedblright\ axis of
$S^{1}\times D^{2}$ as shown in Figure II A-2, such the two \textquotedblleft
bends\textquotedblright\ of $S^{1}\times D^{2}$ lie to the outer sides of the
endmost planes $P_{1}$ and $P_{q}$. These planes cut $\mathbb{R}^{3}$ into
$q+1$ \emph{chambers} (and they cut $S^{1}\times D^{2}$ into $2q$ chambers).
We will say that the $p^{th}$ stage of the Bing decomposition is
\emph{essentially }$\emph{q}$\emph{-chamberable} if, for any $\delta>0$, there
is an ambient isotopy of $\mathbb{R}^{3}$, fixed outside of $S^{1}\times
D^{2}$, such that the image of each of the $2^{p}$ $p^{\text{th}}$ stage solid
tori under the final homeomorphism lies in the $\delta$-neighborhood of one of
these chambers. Generalizing this in a natural manner, we say that the
$p^{th}$ stage of the developed Bing decomposition of $\mathbb{R}^{3}%
\times\lbrack0,\infty)$ (as described above) is \emph{essentially
q-chamberable} if, for any $\delta>0$, there is a level-preserving ambient
isotopy of $\mathbb{R}^{3}\times\lbrack0,\infty)$, fixed outside of
\[
N_{0}^{4}=(S^{1}\times D^{2}\times\lbrack0,1-\eta])\cup(C^{3}\times
\lbrack1-\eta,1+\eta])\text{,}%
\]
where $C^{3}\approx B^{3}$ is the convex hull of $S^{1}\times D^{2}$ in
$\mathbb{R}^{3}$, such that the image of each of the $2^{p}$ components of the
$p^{\text{th}}$ stage neighborhood $N_{p}^{4}$ under the final homeomorphism
lies in the $\delta$-neighborhood of one of the $q+1$ chambers of
$\mathbb{R}^{3}\times\lbrack0,\infty)$ determined by the 3-planes $P_{i}%
^{2}\times\lbrack0,\infty),1\leq i\leq q$.

\begin{figure}[th]
\centerline{
\includegraphics{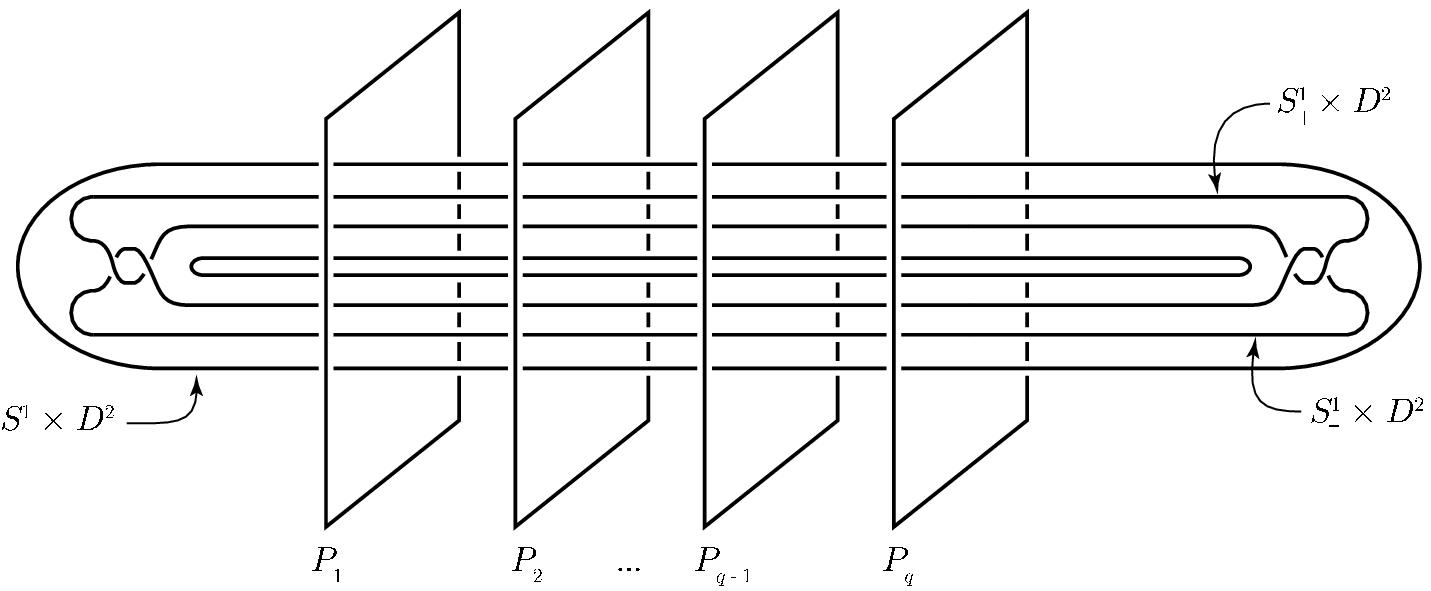}
}\end{figure}

Bing proved that for any $q\geq0$, the $q^{th}$ stage of his decomposition is
essentially $q$-chamberable. The following proposition is a simple extension
of this.\medskip

\textbf{Proposition 1.} \emph{For any }$q\geq0$\emph{, the }$q^{th}$\emph{
stage of the developed Bing decomposition of }$R^{3}\times\lbrack0,\infty
)$\emph{ is essentially q-chamberable.}

\textit{Proof of Proposition 1.} By induction on $q$. See Figure IIA-3. (also,
Bing's pictures in \cite{Bi1} are helpful here). Let the $2$-planes $P_{1}%
^{2},\ldots,P_{q}^{2}$ be given. The first step is to apply the induction
hypothesis separately inside of the two components $N_{-}^{4}$ and $N_{+}^{4}$
of $N_{1}^{4}$, to reposition in them the components of the $q^{\text{th}}$
stage $N_{q}^{4}$ which lie there, as follows. For $N_{-}^{4}$, using the
initial $q-1$ $2$-planes $P_{1}^{2},\ldots,P_{q-1}^{2}$, the induction
hypothesis provides a level-preserving isotopy of $\mathbb{R}^{3}\times
\lbrack0,\infty)$ supported in $N_{-}^{4}$ such that for each of the $2^{q-1}$
components of $N_{q}^{4}$ lying in $N_{-}^{4}$, its image under the final
homeomorphism lies in the $\delta$-neighborhood ($\delta$ arbitrarily small)
of one of the $q$ chambers of $\mathbb{R}^{3}\times\lbrack0,\infty)$
determined by $P_{1}^{2}\times\lbrack0,\infty),\ldots,P_{q-1}^{2}\times
\lbrack0,\infty)$. For $N_{+}^{4}$, one finds a similar isotopy, \textbf{using
instead} the final $q-1$ q-planes $P_{2}^{2},\ldots,P_{q}^{2}$. Let these
isotopies be applied to $\mathbb{R}^{3}\times\lbrack0,\infty)$.

\begin{figure}[th]
\centerline{
\includegraphics{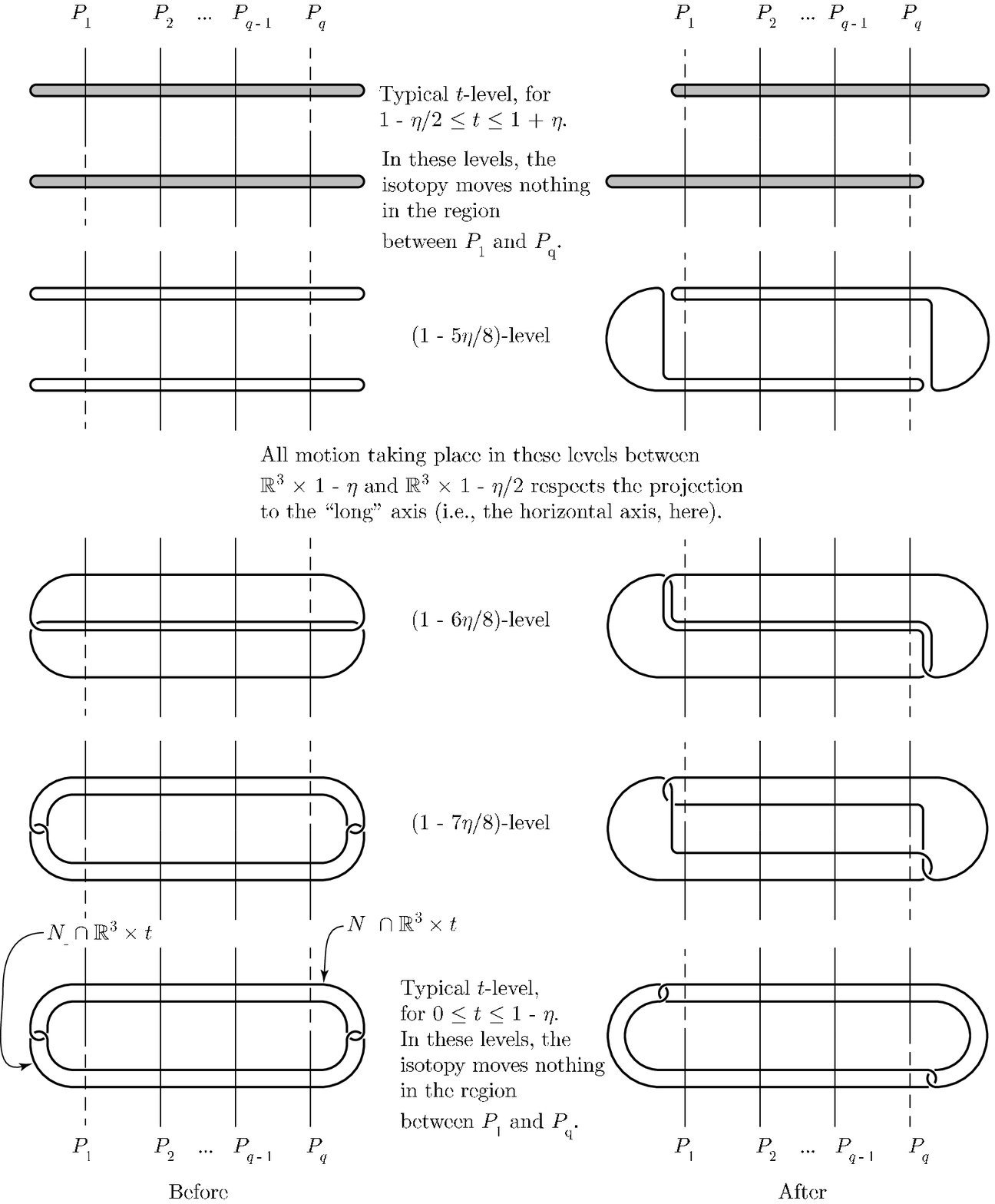}
}\caption{The left side shows how $N_{-}^{4}\cap N_{+}^{4}$ is situated in
$N_{0}^{4}$ before the isotopy is applied. The right side shows how the image
of $N_{-}^{4}\cap N_{+}^{4}$ is situated in $N_{0}^{4}$ after the isotopy is
applied, at the end of Step 1.}%
\end{figure}

The Proposition is completed by isotoping $N_{-}^{4}\cup N_{+}^{4}$ in
$N_{0}^{4}$ in a certain careful manner, as shown in Figure IIA-3.

It is best to describe this level-preserving isotopy in three pieces. Below
the $(1-\eta)$-level, one does in each $t$-level the classic Bing move, whose
motion here takes place only in the two outermost chambers (i.e. the two
chambers whose frontiers are $P_{1}^{2}\times\lbrack0,\infty)$ and $P_{q}%
^{2}\times\lbrack0,\infty))$. In the left most chamber one pulls the loop of
$N_{+}^{4}$ very tight to make it poke into this chamber no more than distance
$\delta$. In the rightmost chamber one similarly pulls the loop of $N_{-}^{4}$
very tight. Of course, these pulling motions stretch the other two end-loops
longer, but that does not matter. Above the $(1-\eta/2)$-level, where
$S_{-}^{2}\times D^{2}$ and $S_{+}^{1}\times D^{2}$ are in their unraveled
positions, the motion again takes place only in the two outermost chambers,
and again the same pulling-tight rule is applied, making small the left
end-loop of $N_{+}^{4}$ and the right end-loop of $N_{-}^{4}$. The nontrivial
part of the isotopy is what takes place in between these two levels. This is
described in Figure IIA-3, which shows the positioning of $N_{-}^{4}\cup
N_{+}^{4}$ in $N_{0}^{4}$ before and after the isotopy is applied. That is,
the left side of Figure IIA-3 shows the original unraveling motion of
$S_{-}^{1}\times D^{2}\cup S_{+}^{1}\times D^{2}$ in the various levels
$\mathbb{R}^{3}\times t,0\leq t\leq1+\eta$, and the right side of Figure IIA-3
shows the new unraveling motion, which has been gotten by applying the
level-preserving isotopy of $N_{0}^{4}$ to the original unraveling motion. The
important thing is that between the 3-planes $P_{1}^{2}\times\lbrack0,\infty)$
and $P_{q}^{2}\times\lbrack0,\infty)$ this new unraveling motion
\emph{respects projection to the \textquotedblleft long\textquotedblright%
\ axis} (i.e., the left-right axis, in the pictures). Hence the chambering
which was already achieved inside of $N_{-}^{4}$ and $N_{+}^{4}$, by the
earlier motions which took place there is not lost. It can be verified that
after applying this pictorially-described isotopy, the desired chambering has
been achieved, completing the proof.\medskip

Strictly speaking, the only part of Proposition 1 that is used in Step 2 is
the motion of $N_{-}^{4}\cup N_{+}^{4}$ in $N_{0}^{4}$ described above. But it
is useful to understand the entire Proposition, for it is a special case of
what follows.\medskip

\textbf{Step 2.} In this step, the motions of Proposition 1 are applied to
shrink the $k$-times-spun Bing decomposition. In a nutshell, the point is that
the unfolded Bing decomposition described in Step 1 is exactly one-half of the
once-spun Bing decomposition (after amalgamating arcs in the former to become
half-discs in the latter), and furthermore, the $k$-times-spun Bing
decomposition, $k\geq2$, can be gotten by $(k-1)$-spinning the unfolded Bing
decomposition (with its elements amalgamated). This entire step can best be
understood by concentrating on the $k=1$ and $k=2$ cases.

To define the $k$-times-spun Bing decomposition of $\mathbb{R}^{3+k}$ into
points and tame $(1+k)$-cells, one starts with half of the original Bing
decomposition, lying in $\mathbb{R}^{2}\times\lbrack0,\infty)$ as shown in
Figure IIA-4, with stages denoted by $B^{1}\times D^{2}$ and $\{B_{\mu}%
^{1}\times D^{2}\mid\mu\in\{-,+\}^{p}\}$. Then one imagines $\mathbb{R}^{3+k}$
as being obtained from $\mathbb{R}^{2}\times\lbrack0,\infty)$ by spinning
$\mathbb{R}^{2}\times\lbrack0,\infty)$ through a $k$-sphere's worth of
directions, keeping $\mathbb{R}^{2}\times0$ fixed. During this spinning, the
arcs of the half-Bing-decomposition sweep out $(1+k)$-cells in $\mathbb{R}%
^{3+k}$. (Equivalently, one can define this decomposition of $\mathbb{R}%
^{3+k}$ to be the restriction to the boundary $(3+k)$-sphere of the product
decomposition of $B^{1}\times D^{2}\times I^{1+k}$ into $(2+k)$-cells. But the
first model seems better suited to this proof.) The $p^{\text{th}}$ stage
defining neighborhoods of the $k$-times-spun Bing decomposition are the
thickened $(1+k)$-spheres $\{S_{\mu}^{1+k}\times D^{2}\mid\mu\in\{-,+\}^{p}%
\}$, which are the images of the cylinders $\{B_{\mu}^{1}\times D^{2}\mid
\mu\in\{-,+\}^{p}\}$ under the spinning.

\begin{figure}[th]
\centerline{
\includegraphics{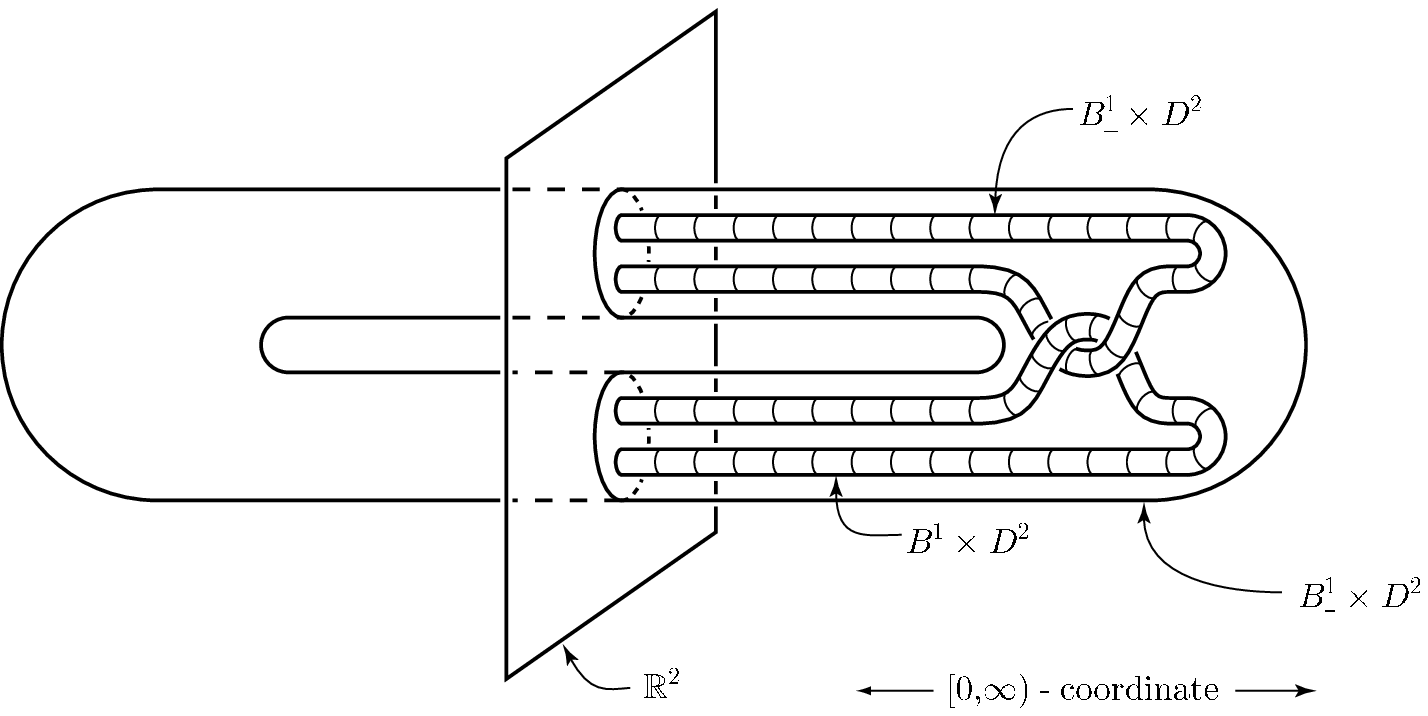}
}\caption{Half of the original Bing decomposition of $\mathbb{R}^{3}$.}%
\end{figure}

In trying to shrink this decomposition, the important coordinate is the
$\mathbb{R}^{1+k}$ coordinate in $\mathbb{R}^{2}\times\mathbb{R}%
^{1+k}=\mathbb{R}^{3+k}$, i.e., the coordinate which is perpendicular to the
fixed plane $\mathbb{R}^{2}$. In the original Bing case, i.e. when $k=0$, this
is the long-and-thin coordinate direction. For higher $k$ it might be thought
of as the flat-and-thin coordinate direction. The point is, if one projects
the $k$-times-spun Bing decomposition onto this $\mathbb{R}^{1+k}$ coordinate,
then each $(1+k)$-disc is projected homeomorphically onto a ball $I^{1+k}$ in
$\mathbb{R}^{1+k}$ (at least if the decomposition is correctly positioned,
e.g. as in Figure IIA-1, bottom frame). Also, the core $(1+k)$-sphere of each
component $S_{\mu}^{1+k}\times D^{2}$ of any stage of the defining
neighborhood sequence is flattened by this projection to be a $(1+k)$-disc.

Now we shift from round thinking to square thinking, in order to talk about
chambers. Regarding $I^{1+k}$ as $[-1,1]^{1+k}$ (which explains the reason for
the letter $I$, we can choose a set of $k$-planes in $\mathbb{R}^{1+k}$, each
being perpendicular to one of the $1+k$ coordinate axes and each passing
through $\operatorname*{int}I^{1+k}$, so that these $k$-planes cut
$\mathbb{R}^{1+k}$ into \emph{chambers}, whose intersections with $I^{1+k}$
are small. For notational purposes, let this collection of $k$-planes be
denoted $P^{k}(\vec{q})$, where the $(1+k)$-tuple of non-negative integers
$\vec{q}=(q_{1},q_{2},\ldots,q_{1+k})$ lists the number of $k$-planes
perpendicular to each axis. Let $\Vert\vec{q}\Vert=\sum_{i=1}^{1+k}q_{i}$.

Crossing these $k$-planes and chambers with $\mathbb{R}^{2}$, we get a
partitioning of $\mathbb{R}^{3+k}$, which we will use to measure smallness
with regard to the spun Bing decomposition. But there is one technical point
which has to be discussed, namely how these $(2+k)$-planes are allowed to pass
through the \textquotedblleft bends\textquotedblright\ in the $S_{\mu}%
^{1+k}\times D^{2}$'s near $\partial I^{1+k}$ (as they must do, when $k\geq
1)$. The idea here is to impose as much regularity as possible. One way to say
this accurately is, given the collection $P^{k}(\vec{q})$ of $k$-planes, and
letting $\epsilon$ be the minimum distance from any of the $k$-planes to
either of the two $k$-faces of $I^{1+k}$ parallel to it, then the
decomposition is made so taut, and the bends in the defining neighborhoods
$\{S_{\mu}^{1+k}\times D^{2}\}$ are made so sharp, and so near $\partial
I^{1+k}$, that for any of the $(1+k)$-choose-$\ell$ coordinate sub-spaces
$\mathbb{R}^{\ell}$ of $\mathbb{R}^{1+k}$ $(1\leq\ell\leq1+k)$, the projection
$\pi:\mathbb{R}^{3+k}\rightarrow R^{\ell}$, when restricted \textbf{over} the
subset $[-1+\epsilon/2,1-\epsilon/2]^{\ell}\subset\mathbb{R}^{\ell}$, is a
product map on all members of the defining neighborhood sequence (i.e., any
bending, or pinching, in the $\mathbb{R}^{\ell}$ coordinate direction occurs
outside of $\pi^{-1}([-1+\epsilon/2,1-\epsilon/2]^{\ell})$. Hence this
restricted part of the $k$-times-spun decomposition (and its defining
neighborhood sequence) looks like the product of the $(k-\ell)$-times-spun
Bing decomposition (and its defining neighborhood sequence) with the cube
$[-1+\epsilon/2,1-\epsilon/2]^{\ell}$. One consequence of this regularity is
that the two basic spun embeddings $\underline{x},\underline{x}_{+}%
:S^{1+k}\times D^{2}\rightarrow S^{1+k}\times D^{2}$ which define the
$k$-times-spun Bing decomposition (cf. start of Step 1) can be assumed to
respect the chambering of $\mathbb{R}^{3+k}$ by the codimension one
hyperplanes $\mathbb{R}^{2}\times P^{k}(\vec{q})$, that is, for each such
chamber $C,\chi_{\underline{+}}(C\cap S^{1+k}\times D^{2})\subset C$. Hence
the various finite compositions of $\chi_{-}$ and $\chi_{+}$ have this
property. This will be important in Proposition 2 below.

Generalizing the earlier definition, we say the $p^{\text{th}}$ stage of the
$k$-times-spun Bing decomposition is \emph{essentially }$q$\emph{-chamberable}
if for any collection $P^{k}(\overset{\rightarrow}{q})$ of $k$-planes as
described above, with $\Vert\vec{q}\Vert=q$, and for any $\delta>0$, there is
an ambient isotopy of $\mathbb{R}^{3+k}$, fixed outside of the $0^{\text{th}}$
stage $S^{1+k}\times D^{2}$, such that the image of each of the $2^{p}$
$p^{\text{th}}$ stage components $S_{\mu}^{1+k}\times D^{2}$ under the final
homeomorphism lies in the $\delta$-neighborhood of one of the chambers of
$\mathbb{R}^{3+k}$ determined by $P^{k}(\vec{q})$.\medskip

\textbf{Proposition 2.} \emph{For any }$q\geq0$\emph{, the }$q^{th}$\emph{
stage of the }$k$\emph{-times-spun Bing decomposition of }$R^{3+k}$\emph{ is
essentially }$q$\emph{-chamberable.}

\emph{Proof of Proposition 2. } By induction on $q$. Let $q_{i}$ be any
nonzero component of $\vec{q}=(q_{1},\ldots,q_{1+k})$. We want to regard the
$k$-times-spun Bing decomposition as being gotten by $(k-1)$-spinning the
unfolded Bing decomposition. This $(k-1)$-spinning is to be done in
$\mathbb{R}^{3+k}=\mathbb{R}^{2}\times\mathbb{R}^{1+k}$, by taking the base
$\mathbb{R}^{3}\times0$ of $\mathbb{R}^{3}\times\lbrack0,\infty)$ and
indentifying it with $\mathbb{R}^{2}\times\mathbb{R}^{1}(e_{i})$, where
$e_{i}$ is the $i^{\text{th}}$ coordinate direction in $\mathbb{R}^{1+k}$ ($i$
as in $q_{i}$, above), and then by spinning $\mathbb{R}^{3}\times
\lbrack0,\infty)$ through the $(k-1)$-sphere's worth of directions in
$\mathbb{R}^{k}(e_{1},\ldots,\widehat{e}_{i},\ldots,e_{1+k})$. The important
thing is: this spinning can be done so that

\begin{itemize}
\item[a)] the $(k-1)$-spin of the defining neighborhood sequence of the
unfolded Bing decomposition (i.e., the $N_{p}^{4}$'s in Step 1) is
corresponded to the defining neighborhood sequence of the $k$-times-spun Bing
decomposition (i.e., to the unions $\cup\{S_{\mu}^{1+k}\times D^{2}\mid\mu
\in\{-,+\}^{p}\}$, for each $p\geq0$; hence in particular the decompositions
will coincide);

\item[b)] the $q_{i}$ $(2+k)$-planes of $\mathbb{R}^{2}\times P^{k}(\vec{q})$
which are perpendicular to $\mathbb{R}^{1}(e_{i})$ become perpendicular to the
$(-\infty,\infty)$-coordinate of $\mathbb{R}^{3}=\mathbb{R}^{2}\times
(-\infty,\infty)$ in the $(k-1)$-spin structure, and

\item[c)] the remaining $q-q_{i}\ (2+k)$-planes of $\mathbb{R}^{2}\times
P^{k}(\vec{q})$ become parallel to the $\mathbb{R}^{3}$ coordinate in the
$(k-1)$-spin structure.
\end{itemize}

Granted this, then Proposition 2 follows using the argument from Proposition
1, as we now explain.

First, one applies the induction hypothesis of Proposition 2 separately to the
two first stage components $S_{-}^{1+k}\times D^{2}$ and $S_{+}^{1+k}\times
D^{2}$. For notation here, let $P_{i,1}^{k}$ and $P_{i,q_{i}}^{k}$ denote the
first and the last (i.e. the extremes) among the subcollection $P_{i}^{k}$ of
$q_{i}$ $k$-planes in $P^{k}(\vec{q})$ which are perpendicular to
$\mathbb{R}^{1}(e_{i})$ (if $q_{i}=1$, then these two $k$-planes coincide).
Let $P_{-}=P^{k}(\vec{q})-\{P_{i,q_{i}}^{k}\}$ and $P_{+}=P^{k}(\vec{q}%
{)}-\{P_{i,1}^{k}\}$. Then applying the induction hypothesis to those of the
$q^{\text{th}}$ stage components $\{S_{u}^{1+k}\mid u\in\{-,+\}^{q}\}$ which
lie inside of $S_{-}^{1+k}\times D^{2}$, using the collection $P_{-}$ of $q-1$
$k$-planes, one obtains an isotopy of $S_{-}^{1+k}\times D^{2}\;rel\;\partial$
such that the image of each of these components under the final homeomorphism
lies in the $\delta$-neighborhood ($\delta$ arbitrarily small) of one of the
chambers of $\mathbb{R}^{3+k}$ determined by $P_{-}$. Similarly, one obtains
an isotopy of $S_{+}^{1+k}\times D^{2}\;rel\;\partial$ such that the image of
each of the remaining $q^{\text{th}}$ stage components under the final
homeomorphism lies in the $\delta$-neighborhood of one of the chambers of
$\mathbb{R}^{3+k}$ determined by $P_{+}$. After applying these two isotopies
in $\mathbb{R}^{3+k}$, one completes the motion of Proposition 2 by applying
the corresponding next motion of Proposition 1, spun $k-1$ times in the
structure described above to become a motion of $\mathbb{R}^{3+k}$. The
$q_{i}$ $2$-planes that are used for the motion of Proposition 1 are the
intersections with $\mathbb{R}^{2}\times\mathbb{R}^{1}(e_{i})$ of the $q_{i}$
perpendicular planes in the collection $\mathbb{R}^{2}\times P_{i}^{k}$. The
fact that the motion from Proposition 1 was level-preserving, combined with
(c) above, ensures that the $\delta$-control imposed by the initial isotopies
supported in $S_{-}^{1+k}\times D^{2}\cup S_{+}^{1+k}\times D^{2}$ is not
diminished at all. Using the same analysis as in Proposition 1, this proof is
seen to be completed. \bigskip
\end{proof}

\centerline{\bf The Interlude}

In the (seemingly long) period between my proofs of Theorem 2 (basically
January 1975) and the nearly definitive Triple Suspension Theorem (Part IV;
October 1976), a certain fact became known, which seemed to make it more
likely that the Multiple Suspension Conjecture was true. For historical
interest, it seems worth recalling,\bigskip

\textsc{Observation. }\emph{Suppose }$H_{1}^{3},H_{2}^{3}$\emph{ are homology
3-spheres. Fix }$k\geq2$\emph{. Then }$\Sigma^{k}(H_{1}^{3}\#H_{2}^{3})\approx
S^{k+3}\Leftrightarrow\Sigma^{k}H_{i}\approx S^{k+3}$\emph{ for }%
$i=1,2$\emph{.\bigskip}

If the Multiple Suspension Conjecture were to be false, it seemed natural that
it would fail precisely for Rokhlin invariant 1 homology 3-spheres, because of
the significant relation they were known to have with regard to $PL$ versus
non-$PL$ triangulations of manifolds (cf. Prologue). But a consequence of the
above observation is the \bigskip

\textsc{Corollary. }\emph{The following conjecture is false: }$\Sigma^{k}%
H^{3}\approx S^{k+3}\Leftrightarrow\mu(H^{3})=0$\emph{, where }$\mu$\emph{ is
the Rokhlin invariant in }$Z/2$\emph{.}

The corollary follows because $\mu$ is additive under connected sum.\bigskip

The implication $\Leftarrow$ of the Observation is implicit in \cite[proof of
Prop 4]{Gl2}, \cite[Step 5]{Gl3}, and is explicitly noted in \cite{ES}. The
implication $\Rightarrow$ (which I realized in February 1975, and others did
independently around then is a consequence of \v{S}tanko's Approximation
Theorem in the codimension one setting (\cite{St3}; its fault (see \cite{AC})
is immaterial here). It says that, granted ${\mathring{c}}(H_{1}^{3}%
\#H_{2}^{3})\times\mathbb{R}^{k-1}$ is a manifold, then the natural closed
subset ${\mathring{c}}(H_{i}^{3}-\operatorname*{int}B^{3})\times
\mathbb{R}^{k-1}(i=1$ or 2) can be re-embedded to have 1-LC complement, hence
(by \cite{Da2} or \cite{Ce}) be collared, hence ${\mathring{c}}H_{i}^{3}%
\times\mathbb{R}^{k-1}$ is a manifold.

\part{The double suspension of any homology $3$--sphere is the image of
$S^{5}$ under a cell-like map}

\indent When working on the double suspension question for an arbitrary
homology 3-sphere $H^{3}$, it is natural first to ask whether it can be
reduced to a cell-like decomposition space problem, i.e., whether $\Sigma
^{2}H^{3}$ is the cell-like image of some manifold. The aim of this section is
to prove that this is so.

\begin{theorem}
\label{existence-of-cell-like-map}\textit{For any homology $n$-sphere $H^{n}$,
there is a cell-like map $f:S^{n+2}\rightarrow\Sigma^{2}H^{n}$ from the
$(n+2)$-sphere onto the double suspension of $H^{n}$.}
\end{theorem}

This was proved independently by J. Cannon in \cite{Can1}.

Recall that the Theorem is readily established in the cases where $H^{n}$ is
known to bound a contractible $(n+1)$-manifold, hence whenever $n\geq4$ (see
the remark following the Proposition, Part II). So the new content of the
theorem is the $n=3$ case (but the proof works for any $n$).

The construction which follows is the second I formulated (in February 1976);
the first (September 1975) is sketched at the end of this part. This second
construction has the advantage that the prerequisites are simpler (they being
only the half-open $h$-cobordism theorem and the local contractibility of the
homeomorphism group of a topological manifold), and also that the non-trivial
point-inverses of $f$ are easier to understand.

The following construction grew out of my attempt to understand a fundamental
manifold in PL-TOP manifold theory: the compact topological manifold $P^{5}$,
homotopy-equivalent to $S^{1}$, such that $\partial P^{5}\approx H^{3}\times
S^{1}$, where $H^{3}$ is any given homology 3-sphere. The importance of such a
manifold $P^{5}$ is that, when the Rokhlin invariant of $H^{3}$ is nonzero,
then $P^{5}$ provides a counterexample to the PL triangulation conjecture for
topological manifolds (cf. Prologue, Section I). (The goal of the double
suspension problem is to show that $P^{5}$ is in fact $cH^{3}\times S^{2}$.
The goal of the Theorem above amounts to constructing a cell-like map $p$ from
$P^{5}$ onto $cH^{3}\times S^{1}$ such that $p$ is a homeomorphism on
$\partial P^{5}$. The above-implied uniqueness of $P^{5}$ is a consequence of
the 6-dimensional topological s-cobordism theorem applied to a cobordism
constructed to join rel $\partial$ two candidates $P_{0}$ and $P_{1}$. The
cobordism, which is homeomorphic to $B^{5}\times S^{1}$, exists because
$P_{0}\cup_{\partial}P_{1}$ is homotopically equivalent to $S^{4}\times S^{1}%
$, hence is homeomorphic to $S^{4}\times S^{1}$ (applying the topological
version of \cite[Thm 1.2]{Sha1} \cite[Thm 1.1]{Sha2}; or see \cite[V, App
C]{KS2}, keeping in mind that $\mathcal{S}^{TOP}(B^{4}\times S^{1})=0$ implies
the above-mentioned fact). However, the uniqueness of $P^{5}$ is not germane here.)

The existence of the manifold $P^{5}$ is a consequence of the work of
Kirby-Siebenmann. In \cite[Section 5]{Si2} Siebenmann presented a construction
which produced $P^{5}\times S^{1}$, and $P^{5}$ itself was an implicit
consequence of this and a splitting theorem (cf. \cite[Remark 5.4]{Si2}). This
was made explicit in \cite{Mat1} and \cite{GS1}. In \cite{Sch}, Scharlemann
presented a simplified construction of $P^{5}$, which was distilled from some
earlier arguments of Kirby and Siebenmann. In doing so he exhibited a
structure on $P^{5}$ that was strong enough to let him make some significant
assertions about topological transversality at dimension 4 (cf. Prologue,
Section I). My construction below of the cell-like map $f:S^{5}\rightarrow
\Sigma^{2}H^{3}$ is a direct extension of the work of Siebenmann and Scharlemann.

One can take the following point of view about $P^{5}$. Observing that
$P^{5}\cup_{\partial}H^{3}\times B^{2}$ is a homotopy 5-sphere, hence is
homeomorphic to $S^{5}$, the problem in constructing $P^{5}$ amounts to
constructing an embedding of $H^{3}\times B^{2}$ into $S^{5}$ so that the
complement has the homotopy type of $S^{1}$. (Such a construction was known
for \textquotedblleft homology type of $S^{1}$\textquotedblright\ in place of
"homotopy type of $S^{1}$, for it was known that any oriented 3-manifold
smoothly embeds in $S^{5}$ with product neighborhood \cite[Cor. 4]{Hi1}.
Recall that the desired embedding $H^{3}\times B^{2}\rightarrow S^{5}$ above
can be made smooth, if and only if the Rokhlin invariant of $H^{3}$ is zero
(see Prologue, Section IV).

The construction of $P^{5}$ in \cite{Sch} can be briefly described as follows
(details are filled in below). Given a homology 3-sphere $H^{3}$, do 1-and
$2$-surgeries to $H^{3}\times T^{2}$ ($T^{2}=2$-torus), with the attaching
maps of the surgeries being confined to some tube $H^{3}\times D^{2}$, to
produce a manifold, denoted $\chi(H^{3}\times T^{2})$, which is
homotopy-equivalent to $S^{3}\times T^{2}$. By the local contractibility of
the homeomorphism group of a manifold, some finite cover of $\chi(H^{3}\times
T^{2})$ is homeomorphic to $S^{3}\times T^{2}$ (see the Proposition below; the
minimum information needed in this paragraph is that the universal cover
$\widetilde{\chi}(H^{3}\times T^{2})$ is homeomorphic to $S^{3}\times
\mathbb{R}^{2})$. Hence $\widetilde{\chi}(H^{3}\times T^{2})\approx
S^{3}\times\mathbb{R}^{2})$ can be compactified in the natural way by adding a
circle at infinity to produce $S^{5}=S^{3}\times\mathbb{R}^{2}\cup S^{1}$.
Letting $B^{2}\subset\operatorname*{int}(T^{2}-D^{2})$ be some $2$-cell, this
produces the desired embedding of $H^{3}\times B^{2}$ into $S^{5}$.

The proof below of the Theorem amounts to showing that embedded in the above
construction there is a natural cell-like map from $S^{5}$ onto $\Sigma
^{2}H^{3}$.

We proceed to the details of the proof. For convenience and completeness, the
preceding part of the argument will be repeated, with full details. We
continue to restrict attention to the $n = 3$ case, although the proof works
for any $n$, without alteration.

\begin{proof}
[Proof of Theorem \ref{existence-of-cell-like-map}]Let $H^{3}$ be any homology
3-sphere. Let $G^{3}=H^{3}-\operatorname*{int}B^{3}$, where $B^{3}$ is some
collared 3-cell in $H^{3}$. Consider the homology 5-cell $G^{3}\times I^{2}$.

From \cite{Ke} (cf. Proposition, Part II) there is a contractible 5-manifold
$M^{5}$ such that $\partial M^{5}=\partial(G^{3}\times I^{2})$ ($M^{5}$ is in
fact unique, by the 6-dimensional $h$-cobordism theorem rel $\partial$).
Inasmuch as $M^{5}$ is the fundamental building block of this proof, it may be
worthwhile pointing out how easily it can be constructed. Kervaire quickly
produced $M^{5}$ by doing 1- and $2$-surgeries in $\operatorname*{int}%
(G^{3}\times I^{2})$, the 1-surgeries to kill the fundamental group, and the
(equal number of) $2$-surgeries to make $H_{2}$ again 0. Another way of
exhibiting $M^{5}$ comes from the observation that $M^{5}\cup_{\partial}%
G^{3}\times I^{2}$ is a homotopy 5-sphere, hence is $S^{5}$. Thus $M^{5}$ can
be regarded as the complement $S^{5}-\alpha(G^{3}\times I^{2})$, where
$\alpha:G^{3}\times I^{2}\rightarrow S^{5}$ is any (PL say) embedding. Any
such $\alpha$ works, because $S^{5}-\alpha(G^{3}\times I^{2})$ is not only
acyclic (by duality), but is simply connected, since $\alpha(G^{3}\times
I^{2})$ has a $2$-dimensional spine. Such an embedding $\alpha$ is most easily
produced by first embedding a $2$-dimensional spine of $G^{3}$ into $S^{5}$,
and then taking a regular neighborhood of the image. That this neighborhood is
homeomorphic to $G^{3}\times I^{2}$ can be verified by comparing the 0,1 and
$2$-handle structure of the regular neighborhood to that of $G^{3}\times
I^{2}$. (Note that the embedding $\alpha$ is in fact unique, by the above
mentioned uniqueness of $M\;rel\;\partial$).

Returning to the proof of the Theorem, one has the natural map $I^{2}%
\rightarrow T^{2}$ which identifies the opposite sides of the square $I^{2}$
to produce the $2$-torus $T^{2}$. Doing these identifications to each
1-dimensional square $g\times\partial I^{2}$ in $G^{3}\times\partial
I^{2}\subset M^{5}$, one produces from $M^{5}$ a manifold-with-boundary
$N^{5}$, with $\partial N^{5}\approx\partial B^{3}\times T^{2}$, and a map
$\phi:N^{5}\rightarrow B^{3}\times T^{2}$ which is a homeomorphism on the
boundary, and is a homotopy equivalence $rel\;\partial$. Now comes the
basic\medskip

\textbf{Proposition (Kirby-Siebenmann)}. \emph{Some finite cover }%
$\widehat{\phi}:\widehat{N}^{5}\rightarrow B^{3}\times\widehat{T}^{2}\approx
B^{3}\times T^{2}$\emph{ is homotopic }$rel\;\partial$\emph{ to a
homeomorphism.\medskip}

Actually, it follows from the work of Kirby-Siebenmann that $\phi$ itself is
homotopic $rel\;\partial$ to a homeomorphism \cite[V, Thm. C.2]{KS2}
\cite{HW}, but that is unnecessarily strong for the purposes at hand.

\emph{Proof of Proposition }(From \cite[Lemma 4.1]{Si2}). Regard $N^{5}$ as a
5-dimensional h-cobordism-with-boundary from $N_{0}^{4}\equiv B^{2}%
\times0\times T^{2}$ to $N_{1}^{4}\equiv B^{2}\times1\times T^{2}$, the
product-boundary being $\partial B^{2}\times\lbrack0,1]\times T^{2}=\partial
N_{0}\times\lbrack0,1]=\partial N_{1}\times\lbrack0,1]$, where we are
regarding $B^{3}$ as $B^{2}\times\lbrack0,1]$. For notational convenience
below, assume without loss that $\phi^{-1}(\partial B^{3}\times T^{2}%
)=\partial N^{5}$. By the half-open 5-dimensional h-cobordism theorem
(\cite{Co}; recall that it uses only engulfing), the restricted map
$\phi|:N^{5}-N_{1}^{4}\rightarrow B^{2}\times\lbrack0,1)\times T^{2}$ is
homotopic $rel\partial$ to a homeomorphism, say $\phi_{0}:N^{5}-N_{1}%
^{4}\overset{\approx}{\longrightarrow}B^{2}\times\lbrack0,1)\times T^{2}$, and
similarly the restricted map $\phi|:N^{5}-N_{0}^{4}\rightarrow B^{2}%
\times(0,1]\times T^{2}$ is homotopic $rel\;\partial$ to a homeomorphism
$\phi_{1}:N^{5}-N_{0}^{4}\overset{\approx}{\longrightarrow}B^{2}%
\times(0,1]\times T^{2}$. The goal now is to make $\phi_{0}$ agree with
$\phi_{1}$ over say $B^{2}\times\lbrack\frac{1}{2},1)\times T^{2}$, by using
meshing (\`{a} la Cernavskii) and taking finite covers and then applying local
contractibility of the homeomorphism group of a manifold. Amplifying this, it
is an easy matter using meshing in the radial direction of $B^{2}$, and also
in the $(0,1)$-coordinate, to arrange that $\phi_{1}\phi_{0}^{-1}|B^{2}%
\times(0,1)\times T^{2}$ respects the $B^{2}\times(0,1)$ coordinate to an
arbitrarily fine degree. Now, let $\widehat{\phi}_{0},\widehat{\phi}_{1}$ be
large finite covers of $\phi_{0},\phi_{1}$ so that the homeomorphism
$\widehat{\phi}_{1}\widehat{\phi}_{0}^{-1}:B^{2}\times(0,1)\times
T^{2}\rightarrow B^{2}\times(0,1)\times T^{2}$ arbitrarily closely respects
the $T^{2}$ coordinate, in addition to the $B^{2}\times(0,1)$ coordinate. That
is, $\widehat{\phi}_{1}\widehat{\phi}_{0}^{-1}$ is arbitrarily close to the
identity. By local contractibility, there is a homeomorphism $\rho:B^{2}%
\times\lbrack0,1)\times T^{2}\rightarrow B^{2}\times\lbrack0,1)\times T^{2}$
such that $\rho=identity$ on $B^{2}\times\lbrack0,\frac{1}{4}]\times T^{2}%
\cup\partial B^{2}\times\lbrack0,1)\times T^{2}$ and $\rho=\widehat{\phi}%
_{1}\widehat{\phi}_{0}^{-1}$ over $\phi B^{2}\times\lbrack\frac{1}{2},1)\times
T^{2}$. Replacing $\widehat{\phi}_{0}$ by $\rho\widehat{\phi}_{0}$, one sees
that the homeomorphisms $\rho\widehat{\phi}_{0}$ and $\widehat{\phi}_{1}$
agree over $B^{2}\times\lbrack\frac{1}{2},1)\times T^{2}$, and so together
they define the desired homeomorphism of $\widehat{N}^{5}$ to $B^{2}%
\times\lbrack0,1]\times T^{2}$.\medskip

Let $\psi:\widehat{N}^{5}\rightarrow B^{3}\times T^{2}$ be the homeomorphism
of finite covers promised by the Proposition. Let $\widetilde{\psi}%
:\widetilde{N}^{5}\rightarrow B^{3}\times\mathbb{R}^{2}$ be an induced
homeomorphism of universal covers. Let $\gamma:\mathbb{R}^{2}\rightarrow
\operatorname*{int}D^{2}$ be a radial homeomorphism, say $\gamma(x)=x/(1+\Vert
x\Vert)$, where $D^{2}$ denotes the unit disc in $\mathbb{R}^{2}$ (the letter
$D^{2}$ is used here, rather than $B^{2}$, for notational distinction and
clarity. Henceforth one should think of the $B^{2}$'s and $B^{3}$'s as coming
from the first coordinate of $B^{3}\times\mathbb{R}^{2}$, and $D^{2}$ as
coming from the second coordinate.) Identifying $B^{5}$ with the joint
$B^{3}\ast S^{1}$, so that $B^{5}-S^{1}=B^{3}\ast S^{1}-S^{1}=B^{3}%
\times\overset{\circ}{c}S^{1}=B^{3}\times\operatorname*{int}D^{2}$, where
$S^{1}=\partial D^{2}$ and $\overset{\circ}{c}S^{1}$ denotes the open cone on
$S^{1}$, then the embedding $id_{B^{3}}\times\gamma:B^{3}\times\mathbb{R}%
^{2}\rightarrow B^{3}\times\operatorname*{int}D^{2}\hookrightarrow B^{3}\ast
S^{1}=B^{5}$ establishes a compactification of $B^{3}\times\mathbb{R}^{2}$ by
a circle to a space homeomorphic to $B^{5}$. Preceding this embedding with the
homeomorphism $\widetilde{\psi}:\widetilde{N}^{5}\rightarrow B^{3}%
\times\mathbb{R}^{2}$ establishes such a compactification of $\widetilde
{N}^{5}$; this will be denoted $\widetilde{N}^{5}\cup S^{1}\approx B^{5}$.

At this point one has produced the manifold $P^{5}$ which disproves the $PL$
triangulation conjecture for topological manifolds (cf. earlier remarks, fifth
paragraph before Proposition.)

The task now is to construct a cell-like surjection $f_{0}:\widetilde{N}%
^{5}\cup S^{1}\rightarrow G^{3}\ast S^{1}(\approx\Sigma^{2}G^{3})$ such that
$f_{0}$ carries $\partial\widetilde{N}^{5}\cup S^{1}(\approx S^{4})$
homeomorphically onto $\partial G^{3}\ast S^{1}(=S^{2}\ast S^{1})$. Thus each
point-inverse of $f_{0}$ will intersect $\partial\widetilde{N}^{5}\cup S^{1}$
in at most a single point. Given such a cell-like map $f_{0}$, one produces
the desired cell-like map $f:S^{5}\rightarrow H^{3}\ast S^{1}(\approx
\Sigma^{2}H^{3})$ simply by gluing 5-balls onto the source and target of
$f_{0}$, and extending $f_{0}$ over these 5-balls by coning. In symbols,
\[
f:S^{5}\approx\widetilde{N}^{5}\cup S^{1}\cup_{\partial}B^{5}\overset
{f_{0}\cup\text{homeo.}}{\longrightarrow}G^{3}\ast S^{1}\cup_{\partial}%
B^{3}\ast S^{1}=(G^{3}\cup_{\partial}B^{3})\ast S^{1}=H^{3}\ast S^{1}.
\]
To construct $f_{0}$, we analyze in more detail the structure of the universal
cover $\widetilde{N}^{5}$ of $N^{5}$ and its compactification $\widetilde
{N}^{5}\cup S^{1}$. For each $(i,j)\in\mathbb{Z}\oplus\mathbb{Z}$, let
$I^{2}(i,j)$ be the square $[i-\frac{1}{4},i+\frac{1}{4}]\times\lbrack
j-\frac{1}{4},j+\frac{1}{4}]$ in $\mathbb{R}^{2}$; it is centered at the point
$(i,j)$ and has side-length $\frac{1}{2}$. Let $G_{1}^{3}=G^{3}-\partial
G^{3}\times\lbrack0,1)\subset\operatorname*{int}G^{3}$, where $\partial
G^{3}\times\lbrack0,2)$ is some boundary collar for $\partial G^{3}$ in
$G^{3}$. In $G^{3}\times\mathbb{R}^{2}$, take each 5-dimensional block
$G_{1}^{3}\times I^{2}(i,j)$, for each pair $(i,j)$ \textbf{except} $(0,0)$,
and replace it, fixing boundary, with a copy of the contractible manifold.
$M^{5}$, which will be denoted $M^{5}(i,j)$. It turns out that this produces
$\widetilde{N}^{5}$. The reasoning here is in two steps. First, it is clear
that $\widetilde{N}^{5}$ can be identified with the space gotten from
$G^{3}\times\mathbb{R}^{2}$ by replacing \textbf{all} of the blocks $G_{1}%
^{3}\times I^{2}(i,j)$ with copies of $M^{5}$. Second, leaving $G_{1}%
^{3}\times I^{2}(0,0)$ unreplaced (the reason for which will become clear)
does not affect the homeomorphism type of the resultant union, because for
example the space
\[
F^{5}\equiv\left(  G^{3}\times\lbrack-\frac{1}{2},\frac{1}{2}]\times
\lbrack-\frac{1}{2},\infty)-\cup_{j=1}^{\infty}G_{1}^{3}\times I^{2}%
(0,j)\right)  \cup\bigcup_{j=1}^{\infty}M^{5}(0,j)
\]
is homeomorphic fixing boundary to the space $(F^{5}-G_{1}^{3}\times
I^{2}(0,0))\cup M^{5}(0,0)$, by a simple sliding motion in the $[-\frac{1}%
{2},\infty)$-coordinate direction which takes each $M(0,j)$ onto $M(0,j-1)$,
for $j\geq1$.

We now construct a simple but important cell-like surjection $\rho
:D^{2}\rightarrow D^{2}$ (see Figure III-2). This map $\rho$ will be the
identity on $\partial D^{2}$, and the nontrivial point-inverses of $\rho$ will
comprise a countable null collection of $2$-discs, each intersecting $\partial
D^{2}$ in a single point (\emph{null} means that for any $\epsilon>0$, there
are only finitely many members of the collection having diameter $\geq
\epsilon)$. For each relatively prime pair $(p,q)\in\mathbb{Z}\oplus
\mathbb{Z}-(0,0)$ (i.e. $gcd\{p,q\}=1$, so this includes all pairs where $p$
or $q$ is $\pm1$), let $I_{\ast}(p,q)=\bigcup_{\ell\geq1}I(\ell p,\ell q)$,
which is a union of squares in $\mathbb{R}^{2}$ converging to infinity in the
direction $\theta$, where $\cos\theta=p/\Vert(p,q)\Vert$ and $\sin
\theta=q/\Vert(p,q)\Vert$. For each $(p,q)$, we wish to join together the
components of $I_{\ast}(p,q)$ by using bands to connect adjacent squares in
the sequence, so that the resultant union, to be denoted $I_{\#}(p,q)$, is
contractible (See Figure III-1; this operation is amplified in the next
paragraph). In fact, $I_{\#}(p,q)$ will be homeomorphic to a $2$-disc minus a
boundary point. This connecting operation is to be done so that if one defines
$D(p,q)=\gamma(I_{\#}(p,q))\cup(p,q)/\Vert(p,q)\Vert\subset D^{2}$, where
$\gamma:\mathbb{R}^{2}\overset{\approx}{\longrightarrow}\operatorname*{int}%
D^{2}$ is as above and the point $(p,q)/\Vert(p,q)\Vert$ lies on $\partial
D^{2}$, then the collection $\{D(p,q)\}$ is a disjoint null collection of
$2$-discs in $D^{2}$, each intersecting $\partial D^{2}$ in a single point.
Given the collection $\{D(p,q)\}$, then $\rho$ can be taken to be any map
$\rho:D^{2}\rightarrow D^{2}$, fixed on $\partial D^{2}$, such that the
nontrivial point-inverses of $\rho$ are precisely the sets $\{D(p,q)\}$.

\setcounter{figure}{0} \renewcommand{\thefigure}{\Roman{part}-\arabic{figure}}

\begin{figure}[th]
\centerline{
\includegraphics{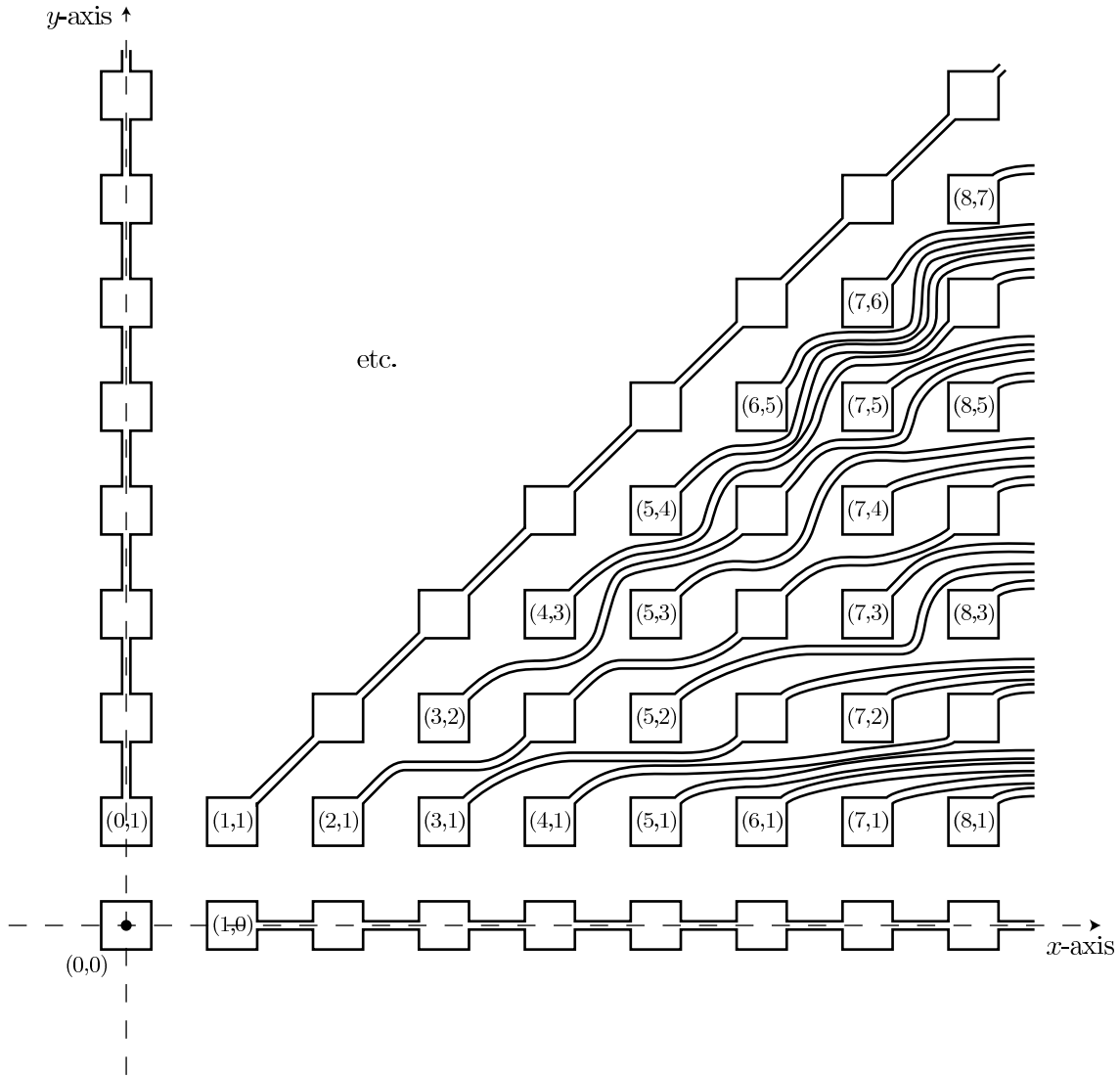}
}\caption{Connecting the squares of $I_{*}(p,q)$ to make $I_{\#}(p,q)$, in
$\mathbb{R}^{2}$. Each labeled pair $(p,q)$ is relatively prime.}%
\end{figure}

\begin{figure}[th]
\centerline{
\includegraphics{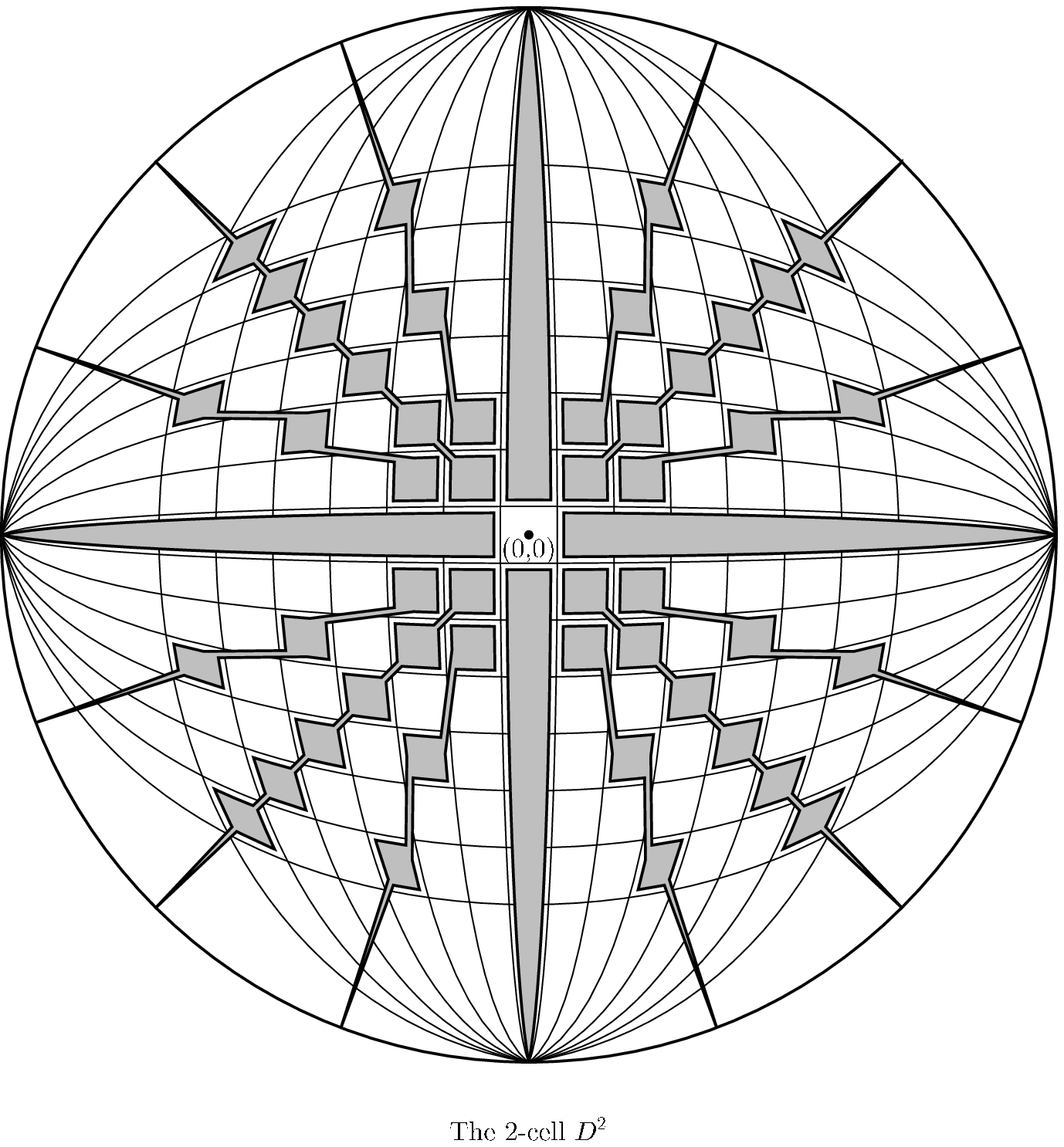}
}\caption{The nontrivial point-inverses of the cell-like surjection
$\rho=D^{2}\to D^{2}$. Each nontrivial point-inverse is homeomorphic to a
2-cell, and intersects $\partial D^{2}$ in a single point.}%
\end{figure}

One way to do the above connecting operation with precision is, given
$I_{\ast}(p,q)$, to adjoin to it the straight ray $\mathbb{R}^{1}%
(p,q)(\approx\lbrack0,\infty))$ which starts at the point $(p,q)$ and passes
through all of the points $(\ell p,\ell q)$, $\ell\geq1$. These rays
$\{\mathbb{R}(p,q)\}$ are all disjoint, but a given $\mathbb{R}(p,q)$ may
unfortunately pass through other $I(p^{\prime},q^{\prime})$'s. To get the
desired disjointness here, one can argue as follows. Let $\pi:\mathbb{R}%
^{2}\rightarrow\mathbb{R}^{2}$ be a map, bounded close to $id_{\mathbb{R}^{2}%
}$ (say by making it the identity on the grid $(\mathbb{Z}+\frac{1}{2}%
)\times\mathbb{R}^{1}\cup\mathbb{R}^{1}\times(\mathbb{Z}+\frac{1}{2}%
)\subset\mathbb{R}^{2})$, such that the only nontrivial point-inverses of
$\pi$ are the sets $\pi^{-1}((i,j))=I^{2}(i,j),(i,j)\in\mathbb{Z}%
\oplus\mathbb{Z}-(0,0)$, and such that for each relatively prime pair $(p,q)$,
$\pi$ leaves invariant (not fixed) the set $\mathbb{R}(p,q)\cap\bigcup
_{\ell\geq1}[\ell p-\frac{1}{2},\ell p+\frac{1}{2}]\times\lbrack\ell
q-\frac{1}{2},\ell q+\frac{1}{2}].$ Then each preimage set $\pi^{-1}%
(\mathbb{R}(p,q))$ looks like $I_{\ast}(p,q)$ with its components joined
together by arcs. So one can take as $I_{\#}(p,q)$ a small thickening of
$\pi^{-1}(\mathbb{R}(p,q))$.

Given $\rho$ as described above, let $\rho_{t}:D^{2}\rightarrow D^{2},0\leq
t\leq1$, be a pseudoisotopy of $\rho_{0}=$ identity to $\rho_{1}=\rho$
(\emph{pseudoisotopy} means here that each $\rho_{t},t<1$, is a homeomorphism).

In order to define the desired cell-like map $f_{0}:\widetilde{N}\cup
S^{1}\rightarrow G^{3}\ast S^{1}$, we first define a certain cell-like map
$g_{0}:G^{3}\ast S^{1}\rightarrow G^{3}\ast S^{1}$. For each $x\in G^{3}$,
define $g_{0}|x\ast S^{1}$ to be the map $\rho_{t(x)}:x\ast S^{1}\rightarrow
x\ast S^{1}$, where $x\ast S^{1}$ is being identified with $D^{2}=0\ast S^{1}$
in the obvious manner, and where $t(x)\in\lbrack0,1]$ is defined using the
previously chosen collar $\partial G^{3}\times\lbrack0,2)$, by
\[
t(x)=\left\{
\begin{array}
[c]{lll}%
t & \text{if} & x\in\partial G^{3}\times t,\;\;t\;\epsilon\;[0,1),\;\text{or}%
\\
1 & \text{if} & x\in G^{3}-\partial G^{3}\times\lbrack0,1)=G_{1}^{3}.
\end{array}
\right.
\]
Thus $g_{0}|\partial G^{3}\ast S^{1}=identity$, and the nontrivial
point-inverses of $g_{0}$ comprise a countable null collection of cell-like
subsets of $G^{3}\ast S^{1}$, each intersecting $S^{1}$ in a single point.
This is because the nontrivial point-inverses of $g_{0}$ are the sets%
\[
G_{1}^{3}\times\gamma(I_{\#}^{2}(p,q))\cup(p,q)/\Vert(p,q)\Vert\subset
G^{3}\times\operatorname*{int}D^{2}\cup S^{1}=G^{3}\ast S^{1}\text{,}%
\]
where $(p,q)$ ranges over all relatively prime pairs in $\mathbb{Z}%
\oplus\mathbb{Z}-(0,0)$, and each of these set is homeomorphic to $c(G_{1}%
^{3}\times I)$.

Now, according to the discussion earlier, we can regard $\widetilde{N}^{5}\cup
S^{1}$ as being obtained from $G^{3}\ast S^{1}$ by removing each block
$G_{1}^{3}\times\gamma(I^{2}(i,j)),\ (i,j)\in\mathbb{Z}\oplus\mathbb{Z}%
-(0,0)$, and replacing it with a copy $M(i,j)$ of $M^{5}$. Since $g_{0}$ sends
each block $G_{1}^{3}\times\gamma(I^{2}(i,j))$ to a single point in
$S^{1}\subset G^{3}\ast S^{1}$, then $g_{0}$ gives rise to a well-defined map
\textbf{after} these replacement operations are done in the \textbf{source} of
$g_{0}$. This new map is the desired $f_{0}:\widetilde{N}^{5}\cup
S^{1}\rightarrow G^{3}\ast S^{1}$.

If this definition is to be made more formally, let $\chi:G^{3}\ast
S^{1}\rightarrow\widetilde{N}^{5}\cup S^{1}$ be a map which restricts on each
block $G_{1}^{3}\times\gamma(I^{2}(i,j))$ to a degree 1 identity-on-boundary
map onto $M(i,j)$, and which is the identity elsewhere, and then note that
$g_{0}\chi^{-1}$ well-defines the desired map $f_{0}$.

The nontrivial point-inverses of $f_{0}$ are each homeomorphic to the
one-point-compactification $F^{5}\cup\infty$ of
\[
F^{5}\equiv(G_{1}^{3}\times\lbrack-\frac{1}{2},\frac{1}{2}]\times\lbrack
\frac{1}{2},\infty)-\cup_{j=1}^{\infty}G_{1}^{3}\times I^{2}(0,j))\cup
\bigcup_{j=1}^{\infty}M^{5}(0,j).
\]
The compactum $F^{5}\cup\infty$ is cell-like because it is contractible, which
in turn follows from the fact that the space consisting of two copies of
$M^{5}$ glued together along $G^{3}\times I\subset\partial M^{5}$, strongly
deformation retracts to either single copy of $M^{5}$, since each space is
contractible (this fact is to be contrasted to the fact that $M^{5}$ itself
does not strongly deformation retract to $G^{3}\times I)$.

This completes the construction of $f_{0}$ and hence $f$. It turns out, then,
that $f$ has a countable null collection of nontrivial point-inverses, each a
contractible ANR. One can seek to improve the point-inverses of $f$, for
example by taking spines to lower their dimension, keeping in mind that the
ultimate goal is to make $f$ a homeomorphism. That leads to the work in Part IV.
\end{proof}

The following paragraphs describe the original proof I formulated that the
double suspension conjecture for any homology 3-sphere is reducible to a
cell-like decomposition problem for $S^{5}$. It rests on some clever
4-dimensional analysis of A. Casson. I found this proof useful only
psychologically, because I was never able to make the decomposition of $S^{5}$
nice enough so that I could work with it. Interestingly, one message of
Cannon's work in \cite{Can2} is that given such a decomposition as I produced
(or he, in \cite{Can2}), it can always be changed into a standard
decomposition which one can work with. This proof is presented here in part to
advertise Cannon's unpublished work, which focuses attention on some
interesting questions in 4-dimensional topology.

\begin{proof}
[Original Proof of Theorem \ref{existence-of-cell-like-map} ]Suppose, then,
that $H^{3}$ is any homology 3-sphere. It was pointed out in the Interlude
following Part II how the implication $\Sigma^{2}(H^{3}\#H^{3})\approx
S^{5}\Rightarrow\Sigma^{2}H^{3}\approx S^{5}$ followed from a certain
splitting construction. Hence it suffices here to work with a homology
3-sphere of Rokhlin invariant 0, namely $H^{3}\#H^{3}$, instead of an
arbitrary $H^{3}$. The following argument shows:\medskip

\noindent\emph{Given a homology 3-sphere }$H^{3}$\emph{ of Rokhlin invariant
}$0$\emph{, there is a cell-like map }$p:P^{5}\rightarrow cH^{3}\times S^{1}%
$\emph{ from some }$5$\emph{-manifold-with-boundary }$P^{5}$\emph{ onto
}$cH^{3}\times S^{1}$\emph{ such that the only nontrivial point-inverses of
}$p$\emph{ lie in }$p^{-1}(c\times S^{1})\subset\;\operatorname*{int}P^{5}%
$\emph{.\medskip}

(Recall that if one wishes to produce from $p$ a cell-like map from $S^{5}$
onto $\Sigma^{2}H^{3}$, one merely glues copies of $H^{3}\times B^{2}$ onto
the boundaries of $P^{5}$ and $cH^{3}\times S^{1}$, and extends $p$ over these
sets via the identity.)

If $H^{3}$ has Rokhlin invariant 0, it is a standard consequence that
$H^{3}=\partial M^{4}$, where $M^{4}$ is some simply-connected parallelizable
$PL$ manifold such that $M^{4}\cup c\partial M^{4}$ has the homotopy type of
$\#_{k}S^{2}\times S^{2}$, the connected sum of $k$ copies of $S^{2}\times
S^{2}$, for some $k$. One would like to be able to do 4-dimensional
(simply-connected) surgery to $\operatorname*{int}M^{4}$, to convert $M^{4}$
to a contractible 4-manifold $N^{4}$ (which could even be topological for the
purposes at hand). For then one could let $P^{5}=N^{4}\times S^{1}$, and
$p=n\times id_{S^{1}}$, where $n:N^{4}\rightarrow cH^{3}$ is gotten by
collapsing to a point a spine of $N^{4}$. Unfortunately, such surgery is not
known to be possible. But, Andrew Casson has shown that such surgery is
possible to some extent, and using his work one can construct the desired
manifold $P^{5}$ and map $p:P^{5}\rightarrow cH^{3}\times S^{1}$.

Casson has shown the following:\footnote{This is unpublished, although some
notes from Casson's lectures, taken by C. Gordon and R. Kirby, have been in
circulation since Casson's work in 1974.} There exist disjoint open sets
$U_{1},\cdots,U_{k}$ in $\operatorname*{int}M^{4}$, bounded away from
$\partial M^{4}$, such that

\begin{itemize}
\item[(i)] each $U_{i}$ is proper-homotopy-equivalent to $S^{2}\times
S^{2}-\ast$ ($\ast=\operatorname*{point}$), in the following special manner:
Each $U_{i}$ is homeomorphic to $S^{2}\times S^{2}-C_{i}$, where $C_{i}$ is
some compact subset of $S^{2}\times S^{2}$ which is cell-like and satisfies
the cellularity criterion of McMillan \cite{McM1} (which is that the end of
$S^{2}\times S^{2}-C_{i}$ is 1-connected; from this one can show that
$S^{2}\times S^{2}-C_{i}$ is proper-homotopy-equivalent to $S^{2}\times
S^{2}-\ast)$, and

\item[(ii)] the homomorphism $\oplus_{i=1}^{k}H_{2}(U_{i})\rightarrow
H_{2}(M)$ is an intersection-form-preserving isomorphism.
\end{itemize}

As a consequence of (i) and (ii), it follows that $M^{4}-\bigcup_{i=1}%
^{k}U_{i}$ is \v{C}ech homotopically 3-connected.

(One of Casson's fundamental questions is: \emph{Are the ends of the }$U_{i}%
$\emph{'s diffeomorphic (or even homeomorphic) to }$S^{3}\times R^{1}$\emph{?}
If so, one can do surgery on $M^{4}$ by replacing each $U_{i}$ by an open
4-cell, thereby producing the desired contractible manifold $N^{4}$ mentioned
above. Although unnecessary for our purposes, it may be worth recalling what
Casson's typical open set $U_{i}$ looks like. Regard $S^{2}\times S^{2}$ as
consisting of a 0-handle $B^{4}$, two $2$-handles attached to a pair of
linking solid tori in $\partial B^{4}$, and a 4-handle. Let $W_{1},W_{2}$ be
two Whitehead continua in $\partial B^{4}$, each the familiar intersection of
a nest of solid tori, constructed in the two given solid tori in $\partial
B^{4}$ (\cite[\S 4]{Wh1}; see also \cite[\S 11]{Bi4}). Thus $W_{1}$ and
$W_{2}$ link geometrically in $\partial B^{4}$ Let $C\approx\Sigma W_{1}%
\vee\Sigma W_{2}$ be gotten by coning $W_{1}\cup W_{2}$ to the origin in
$B^{4}$, and then coning $W_{1}$ and $W_{2}$ separately to the center points
of the respective $2$-handles in which they lie. Then $U=S^{2}\times S^{2}-C$
is the model open set of Casson. His general open set is gotten by replacing
$W_{1}$ and $W_{2}$ above by two Whitehead-like continua which are produced by
suitably ramifying the original Whitehead construction.)

Let $U_{i}$ be any open 4-manifold as described in (i) above.\medskip

\textbf{Proposition.} $U_{i}\times S^{1}$\emph{ is homeomorphic to }%
$(S^{2}\times S^{2}-\ast)\times S^{1}$\emph{, in such a manner that the
}$S^{1}$\emph{ coordinate is preserved at }$\infty$\emph{. That is, there is a
homeomorphism of pairs }$h:(U_{i}\bigcup\infty,\infty)\times S^{1}%
\overset{\rightarrow}{\approx}(S^{2}\times S^{2},\ast)\times S^{1}$\emph{ such
that }$h|\infty\times S^{1}=$\emph{ identity, where }$U_{i}\bigcup\infty
$\emph{ denotes the one-point compactification of }$U_{i}$.\medskip

\emph{Proof of Proposition. }This is a cell-like decomposition problem.
Regarding $U_{i}$ as the subset $S^{2}\times S^{2}-C_{i}$ of $S^{2}\times
S^{2}$, the assertion is that there is a homeomorphism $h:(S^{2}\times
S^{2}/C_{i})\times S^{1}\rightarrow S^{2}\times S^{2}\times S^{1}$ such that
$h$ carries $\{C_{i}\}\times S^{1}$ onto $\ast\times S^{1}$ via the identity.
To construct $h$, it suffices as usual to construct a map $f:S^{2}\times
S^{2}\times S^{1}\rightarrow S^{2}\times S^{2}\times S^{1}$ such that $f$ is a
homeomorphism over the complement of $\ast\times S^{1}$, and $f^{-1}%
(\ast\times t)=C_{i}\times t$ for each $t\in S^{1}$. For then one can define
$h=f\pi^{-1}$, where $\pi:S^{2}\times S^{2}\times S^{1}\rightarrow(S^{2}\times
S^{2}/C_{i})\times S^{1}$ is the natural quotient map. The map $f$ is
constructed by shrinking the decomposition $\{C_{i}\times t\mid t\in S^{1}\}$
of $S^{2}\times S^{2}\times S^{1}$ keeping $\ast\times S^{1}$ fixed, where
without loss $\ast\in C_{i}$. This shrinking can be done using a dual skeleton
engulfing argument in the manner of \cite[Thm 1]{EG}.\medskip

Given the Proposition, one produces the manifold $P^{5}$ by removing from
$M^{4}\times S^{1}$ each open set $U_{i}\times S^{1}$, and sewing in its place
$\operatorname*{int}B^{4}\times S^{1}$. To do this precisely, let
$\alpha:\operatorname*{int}B^{4}-0\rightarrow S^{2}\times S^{2}-\ast$ be an
embedding such that $S^{2}\times S^{2}-(\ast\cup$ image $(\alpha))$ is compact
(i.e., image $(\alpha)$ is a deleted neighborhood of $\ast$ in $S^{2}\times
S^{2}$) and such that the $\alpha$-image of the $\partial B^{4}$-end of
$\operatorname*{int}B^{4}-0$ lies toward $\ast$ . Let $\beta=\alpha\times
id_{S^{1}}:(\operatorname*{int}B^{4}-0)\times S^{1}\rightarrow(S^{2}\times
S^{2}-\ast)\times S^{1}$. Then we can define $P^{5}=M^{4}\times S^{1}%
-\bigcup_{i=1}^{k}(U_{i}\times S^{1}-h_{i}^{-1}$ (image $\beta$))
$\bigcup_{i=1}^{k}(\operatorname*{int}B^{4}\times S^{1})_{i}$ where, for each
$i$, the open subset $(\operatorname*{int}B^{4}-0)\times S^{1}$ of the
$i^{\text{th}}$ copy of $\operatorname*{int}B^{4}\times S^{1}$ is identified
to $h_{i}^{-1}$ (image $\beta$) by the homeomorphism $h_{i}^{-1}\beta$.

This manifold $P^{5}$ is known to topologists, being the unique manifold which
is homotopy equivalent, fixing boundary, to $cH^{3}\times S^{1}$ (cf.
introductory remarks of Part III; here however $P^{5}$ is not so interesting,
because the Rokhlin invariant of $H^{3}$ is 0). What is useful about the above
description of $P^{5}$ is that it provides a layering of $P^{5}$ into
cell-like subsets. Namely, for each $\theta\in S^{1}$, let $Q_{\theta}%
=(M^{4}-\bigcup_{i=1}^{k}(U_{i}\times\theta))\cup\bigcup_{i=1}^{k}%
(\operatorname*{int}B^{4}\times\theta)_{i}\subset P^{5}$ (although $Q_{\theta
}$ is compact, it is not necessarily an ANR, because its two parts may not
match up very smoothly where they come together). Let $H^{3}\times
\lbrack0,1)\hookrightarrow M^{4}-\bigcup_{i=1}^{k}U_{i}$ be an open boundary
collar for $\partial M^{4}$ in $M^{4}$. Define $p:P^{5}\rightarrow
cH^{3}\times S^{1}$ to be the \textquotedblleft
level-preserving\textquotedblright\ map gotten by sending the collar
$H^{3}\times\lbrack0,1)\times S^{1}$ homeomorphically onto $(cH^{3}-c)\times
S^{1}$, and by sending each compact set $Q_{\theta}-H^{3}\times\lbrack
0,1)\times\theta$ to the point $c\times\theta,\theta\in S^{1}$.\medskip

\textbf{Proposition}. $p$\emph{ is cell-like.}\textit{\medskip}

It follows from properties (i) and (ii) above, together with the construction
of $P^{5}$, that for any open interval $(\theta_{0},\theta_{1})$ in $S^{1}$,
the manifold $p^{-1}(\theta_{0},\theta_{1}))$ is contractible. Hence each
point-inverse of $p$ has arbitrarily small contractible neighborhoods.

This completes the description of my original map $p:P^{5}\rightarrow
cH^{3}\times S^{1}$.
\end{proof}

\part{The triple suspension of any homology $3$--sphere is $S^{6}$}

The purpose of this part is to show how I proved the following

\begin{theorem}
[Triple Suspension]\textit{For any homology 3-sphere $H^{3}$, the triple
suspension $\Sigma^{3}H^{3}$ is homeomorphic to $S^{6}$.} \indent Recall that
Cannon subsequently has improved this by showing that $\Sigma^{2}H^{3}\approx
S^{5}$ \cite{Can1}.
\end{theorem}

My proof (done in October 1976) made use of my earlier suspension work,
specifically the result that $\Sigma^{2}sp(H^{3})\approx S^{6}$, where
$sp(H^{3})$ is the homology $4$-sphere gotten by spinning $H^{3}$ (defined
below). Aside from this fact, this part does not use anything from Parts I and
II. However, it does use very strongly the construction of Part III.

Given a compact metric space $H$ containing a collared $n$-cell $B^{n}$ (in
general, $H$ will be a homology $n$-sphere), the $k$-\textit{spin} of
$H(k\geq0)$ is the space $sp^{k}H\equiv\partial((H-\operatorname*{int}%
B^{n})\times I^{k+1})$, i.e., $sp^{k}H=(H-\operatorname*{int}B^{n}%
)\times\partial I^{k+1})\cup\partial B^{n}\times I^{k+1}\subset
(H-\operatorname*{int}B^{n})\times I^{k+1}$ (note that the 0-spin of $H$ is
the double of $H-\operatorname*{int}B^{n}$ along $\partial B^{n}$).

After proving the above theorem, I realized that the result could be
interpreted as a converse to work of Glaser \cite{Gl2}-\cite{Gl4} in 1969-70,
by casting it in the following manner.

\begin{theorem}
\textit{Given any compact metric space $H$ containing a collared $n$-cell, and
given any $k\geq2$, then
\[
\Sigma^{k}H\approx S^{n+k}\Leftrightarrow\Sigma^{2}(sp^{k-2}(H))\approx
S^{n+k}.
\]
}
\end{theorem}

The work of Glaser amounted to establishing the implication $\Rightarrow$ (see
\cite{Si4}, where this argument, extracted from Glaser's work, is presented in
a single page). My work below contains a proof for the implication
$\Leftarrow$ (the $k=2$ case is discussed in the Interlude following Part II).
The real case of interest here is when $k=3$ and $H$ is a homology 3-sphere,
in which case we know from Part II that the right hand conclusion is true. The
goal of the remainder of Part IV is to prove this special case, i.e., to show
that if $\Sigma^{2}sp(H^{3})\approx S^{6}$, then $\Sigma^{3}H^{3}\approx
S^{6}$.

The idea of the proof is to use the hypothesis to understand better a certain
cell-like map $\overline{f}:S^{6}\rightarrow\Sigma^{3}H^{3}$, which is
constructed in completely analogous fashion to the cell-like map
$f:S^{5}\rightarrow\Sigma^{2}H^{3}$ from Part III. Before expanding on this,
we emphasize that $\overline{f}$ is \textbf{not} the suspension of $f$; that
would make $\overline{f}$ have the undesirable feature of having uncountably
many nontrivial point-inverses. Instead, $\overline{f}$ is to have a countable
null sequence of nontrivial point-inverses, just as $f$ did.

The rule for constructing $\overline{f}$ is everywhere in Part III to increase
the dimension from 5 to 6, by increasing the dimension of the second
coordinate from 2 to 3. Thus, $G^{3}\times I^{2}$ becomes $G^{3}\times I^{3}$;
$M^{5}$ becomes $M^{6}$, which is the (unique) contractible manifold such that
$\partial M^{6}=\partial(G^{3}\times I^{3})$; and $N^{5}$ becomes $N^{6}$,
with its rel-boundary homotopy equivalence $\overline{\phi}:N^{6}\rightarrow
B^{3}\times T^{3}$. (The bar symbol, e.g. $\overline{f},\overline{\phi
},\overline{\rho},\overline{g}_{0}$, is used over maps which are the
one-dimension-higher analogues of maps from Part III.) The compactification
$\widetilde{N}^{5}\cup S^{1}$ becomes $\widetilde{N}^{6}\cup S^{2}$, using the
radial homeomorphism $\overline{\gamma}(x)=x/(1+\Vert x\Vert):\mathbb{R}%
^{3}\overset{\approx}{\rightarrow}\operatorname*{int}D^{3}$, where $D^{3}$ is
the unit ball in $\mathbb{R}^{3}$. In the construction of the map
$\overline{\rho}:D^{3}\rightarrow D^{3}$, which is the analogue of $\rho
:D^{2}\rightarrow D^{2}$, one works in $\mathbb{R}^{3}$ instead of
$\mathbb{R}^{2}$, using the 3-dimensional cubes $I^{3}(i,j,k)\equiv\lbrack
i-\frac{1}{4},i+\frac{1}{4}]\times\lbrack j-\frac{1}{4},j+\frac{1}{4}%
]\times\lbrack k-\frac{1}{4},k+\frac{1}{4}]$, for $(i,j,k)\in\mathbb{Z}%
\oplus\mathbb{Z}\oplus\mathbb{Z}-(0,0,0)$. In particular, Figure III-1 becomes
a 3-dimensional picture, to be compactified by adding a $2$-sphere at $\infty
$. Each nontrivial point-inverse of $\overline{\rho}$ is a 3-cell
$\overline{\gamma}(I_{\#}^{3}(p,q,r))\cup(p,q,r)/\Vert(p,q,r)\Vert\subset
D^{3}$, where $(p,q,r)\in\mathbb{Z}\oplus\mathbb{Z}\oplus\mathbb{Z}-(0,0,0)$
is a relatively prime triple, i.e., a triple such that $gcd\{p,q,r\}=1$. These
changes and adaptations are all routine, and should require only a few moments
thought to absorb.

The idea of the following proof is this. Given the hypothesis that $\Sigma
^{2}sp(H^{3})\newline(\equiv\Sigma^{2}\partial(G^{3}\times I^{2}))$ is
homeomorphic to $S^{6}$, or more germanely that $c(\partial(G^{3}\times
I^{2}))\times\mathbb{R}^{1}$ is a manifold, it will be shown that then the
point-inverses of $\overline{f}:S^{6}\rightarrow\Sigma^{3}H^{3}$ can be
improved to be cellular arcs, which are embedded in $S^{6}$ in a certain
regular manner. Then it will be shown that this decomposition of $S^{6}$ into
points and arcs is shrinkable, and consequently $\overline{f}$ is approximable
by homeomorphisms.

The precise point of view that we adopt for this procedure is the following.
The map $\overline{f}:S^{6}\rightarrow\Sigma^{3}H^{3}$ will be factored into
two cell-like maps
\[
\overline{f}=\beta\alpha:S^{6}\overset{\alpha}{\longrightarrow}Q^{6}%
\overset{\beta}{\longrightarrow}\Sigma^{3}H^{3},
\]
and each of these maps will be shown to be approximable by homeomorphisms.
This approximating will be trivial for $\alpha$, but will require some work
for $\beta$ (whose nontrivial point-inverses are arcs).

A key item for the construction below is the following description of $M^{6}$.
It is inspired by Glaser's constructions in \cite{Gl2}-\cite{Gl4}.

Define
\[
M_{\ast}^{6}=G^{3}\times I^{2}\times\lbrack0,1]/\{G_{1}^{3}\times I_{1}%
^{2}\times t\mid\frac{1}{3}\leq t\leq\frac{2}{3}\}\text{,}%
\]
that is, $M_{\ast}^{6}$ is the quotient space gotten from $G^{3}\times
I^{2}\times\lbrack0,1]$ by identifying to points each of the subsets
$G_{1}^{3}\times I_{1}^{2}\times t,\frac{1}{3}\leq t\leq\frac{2}{3}$, where
$G_{1}^{3}\subset\operatorname*{int}G^{3}$ and $I_{1}^{2}\subset
\operatorname*{int}I^{2}$ denote smaller copies of $G^{3}$ and $I^{2}$,
obtained as usual by taking complements of open collars. It is in the
following claim that the hypothesis of the theorem is used. (compare
\cite[Prop. 1]{Gl2}).\medskip

\textbf{Claim.} $M_{\ast}^{6}$\emph{ is a contractible manifold.\medskip}

\textit{Proof.} $M_{\ast}^{6}$ is contractible because it is gotten by
identifying to an arc the spine of a manifold. To see that $M_{\ast}^{6}$ is a
manifold, one uses the hypothesis of the theorem, which amounts to assuming
that $c(\partial(G^{3}\times I^{2}))\times\mathbb{R}^{1}$ is a manifold. From
this one sees that $M_{\ast}^{6}$ is a manifold along the interior of the
arc-spine. This leaves the two endpoints of the arc-spine to analyze. Each
endpoint has a deleted neighborhood homeomorphic to $L^{5}\times\mathbb{R}%
^{1}$, where $L^{5}=G^{3}\times I^{2}\cup_{\partial}c(\partial(G^{3}\times
I^{2}))$. Since $L^{5}$ is homotopically a 5-sphere and $L^{5}\times
\mathbb{R}^{1}$ is a manifold, it follows that $L^{5}\times\mathbb{R}^{1}$ is
homeomorphic to $S^{5}\times\mathbb{R}^{1}$ (see \cite[App. 1]{Si6}). Thus
$M_{\ast}^{6}$ is a manifold, henceforth denoted $M^{6}$.

By analogy with Part III, we regard the space $(\widetilde{N}^{6}\cup
S^{2})\cup_{\partial}B^{6}\ (\approx S^{6})$ and the map $\overline
{f}:(\widetilde{N}^{6}\cup S^{2})\cup_{\partial}B^{6}\rightarrow H^{3}\ast
S^{2}\ (\approx\Sigma^{3}H^{3})$ as being obtained by doing modifications to
the source of a certain cell-like map
\[
\overline{g}_{0}\cup homeo:G^{3}\ast S^{2}\cup_{\partial}B^{6}\rightarrow
G^{3}\ast S^{2}\cup_{\partial}B^{3}\ast S^{2}=H^{3}\ast S^{2}\text{.}%
\]
The map $\overline{g}_{0}:G^{3}\ast S^{2}\rightarrow G^{3}\ast S^{2}$ is
defined, as earlier, by letting $\overline{g}_{0}|x\ast S^{2}\rightarrow x\ast
S^{2}$, for $x\in G^{3}$, where $\{\overline{\rho}_{t}\}$ is a pseudoisotopy
of $\overline{\rho}_{1}=\overline{\rho}:D^{3}\rightarrow D^{3}$ to
$\overline{\rho}_{0}\equiv id_{D^{3}}$, and $D^{3}=x\ast S^{2}$. The
nontrivial point-inverses of $\overline{g}_{0}$ are the sets $\{G_{1}%
^{3}\times\overline{\gamma}(I_{\#}^{3}(p,q,r))\cup(p,q,r)/\Vert(p,q,r)\Vert
\}$, where $(p,q,r)\in\mathbb{Z}\oplus\mathbb{Z}\oplus\mathbb{Z}-(0,0,0)$ is a
relatively prime triple, as above. In order to be able to uniformly describe
these nontrivial point-inverses of $\overline{g}_{0}$, we choose for each
triple $(p,q,r)$ a neighborhood of $\overline{\gamma}(I_{\#}^{3}(p,q,r))$ in
$\operatorname*{int}D^{3}$, to be uniformly denoted $I^{2}\times
\lbrack0,\infty)$, such that all of these neighborhoods are disjoint, and such
that the spine $I_{1}^{2}\times\lbrack1,\infty)$ of $I^{2}\times
\lbrack0,\infty)$ coincides with the set $\overline{\gamma}(I_{\#}%
^{3}(p,q,r))$. From now on we will refer to the nontrivial point-inverses of
$\overline{g}_{0}$ as the sets $\{G_{1}^{3}\times I_{1}^{2}\times
\lbrack1,\infty)\cup\infty\}$, which have disjoint pinched neighborhoods
$\{G^{3}\times I^{2}\times\lbrack0,\infty)\cup\infty\}$ in $G^{3}\ast S^{2}$,
and it will be understood that the indices $\{(p,q,r)\}$ are being suppressed.

Recall that the space $\widetilde{N}^{6}\cup S^{2}\ (\approx B^{6})$ is
regarded as being obtained from $G^{3}\ast S^{2}$ by replacing each of the
blocks $G_{1}^{3}\times\overline{\gamma}(I^{3}(i,j,k)),\ (i,j,k)\in
\mathbb{Z}\oplus\mathbb{Z}\oplus\mathbb{Z}-(0,0,0)$, by a copy $M^{6}(i,j,k)$
of $M^{6}$. Using the above explicit description of $M^{6}$ as $M_{\ast}^{6}$,
we can equivalently describe $\widetilde{N}^{6}\cup S^{2}$ as the quotient
space obtained from $G^{3}\ast S^{2}$ by the acyclic quotient map
$\overline{\chi}:G^{3}\ast S^{2}\rightarrow\widetilde{N}^{6}\cup S^{2}$, where
$\overline{\chi}$ makes the following quotient-identifications in each
nontrivial point-inverse of $\overline{g}_{0}$:
\[
G^{3}\times I^{2}\times\lbrack0,\infty)\overset{\left.  \overline{\chi
}\right\vert }{\longrightarrow}G^{3}\times I^{2}\times I^{2}\times
\lbrack0,\infty)/\{G_{1}^{3}\times I_{1}^{2}\times t\mid\ell+\frac{1}{3}\leq
t\leq\ell+\frac{2}{3}\text{, }\ell\geq1\}
\]
(cf. above convention on notation). See Figure IV-1. Here, one is thinking of
the copies $M^{6}(i,j,k)$ of $M^{6}$ in $\widetilde{N}^{6}\cup S^{2}$ as being
say the various blocks $\overline{\chi}(G^{3}\times I^{2}\times\lbrack
\ell,\ell+1])$, for $\ell\geq1$ (actually, these blocks are enlarged-by-collar
versions of the genuine $M^{6}(i,j,k)$'s, but that is not important). Using
this description of $\widetilde{N}^{6}\cup S^{2}$, then the map $\overline
{f}_{0}:\widetilde{N}^{6}\cup S^{2}\rightarrow G^{3}\ast S^{2}$ is the map
which identifies to a point each of the sets $\overline{\chi}(G_{1}^{3}\times
I_{1}^{2}\times\lbrack1,\infty)\cup\infty)$. Note that this typical nontrivial
point-inverse of $\overline{f}_{0}$ looks like an infinite sequence of
(singlefold) suspensions of $G_{1}^{3}\times I_{1}^{2}$ (except that the first
\textquotedblleft suspension\textquotedblright, $G_{1}^{3}\times I_{1}%
^{2}\times\lbrack1,\frac{4}{3}]/G_{1}^{3}\times I_{1}^{2}\times\frac{4}{3}$,
is really a cone), strung together by arcs joining adjacent members of the
sequence, and then compactified by the point at $\infty$. See Figure IV-1,
second from the top.

There is now an apparent factoring of $\overline{f}_{0}$ into two cell-like
maps,
\[
\overline{f}_{0}=\beta_{0}\alpha_{0}:\widetilde{N}^{6}\cup S^{2}%
\overset{\alpha_{0}}{\longrightarrow}Q_{0}^{6}\overset{\beta_{0}%
}{\longrightarrow}G^{3}\ast S^{2},
\]
which is described by factoring into two maps each of the trivial restricted
maps
\[
\overline{f}_{0}|:\overline{\chi}(G_{1}^{3}\times I_{1}^{2}\times
\lbrack1,\infty)\cup\infty)\rightarrow\infty\in S^{2}%
\]
(see Figure IV-1).

\begin{figure}[th]
\setcounter{figure}{0} \renewcommand{\thefigure}{\Roman{part}-\arabic
{figure}} \centerline{
\includegraphics{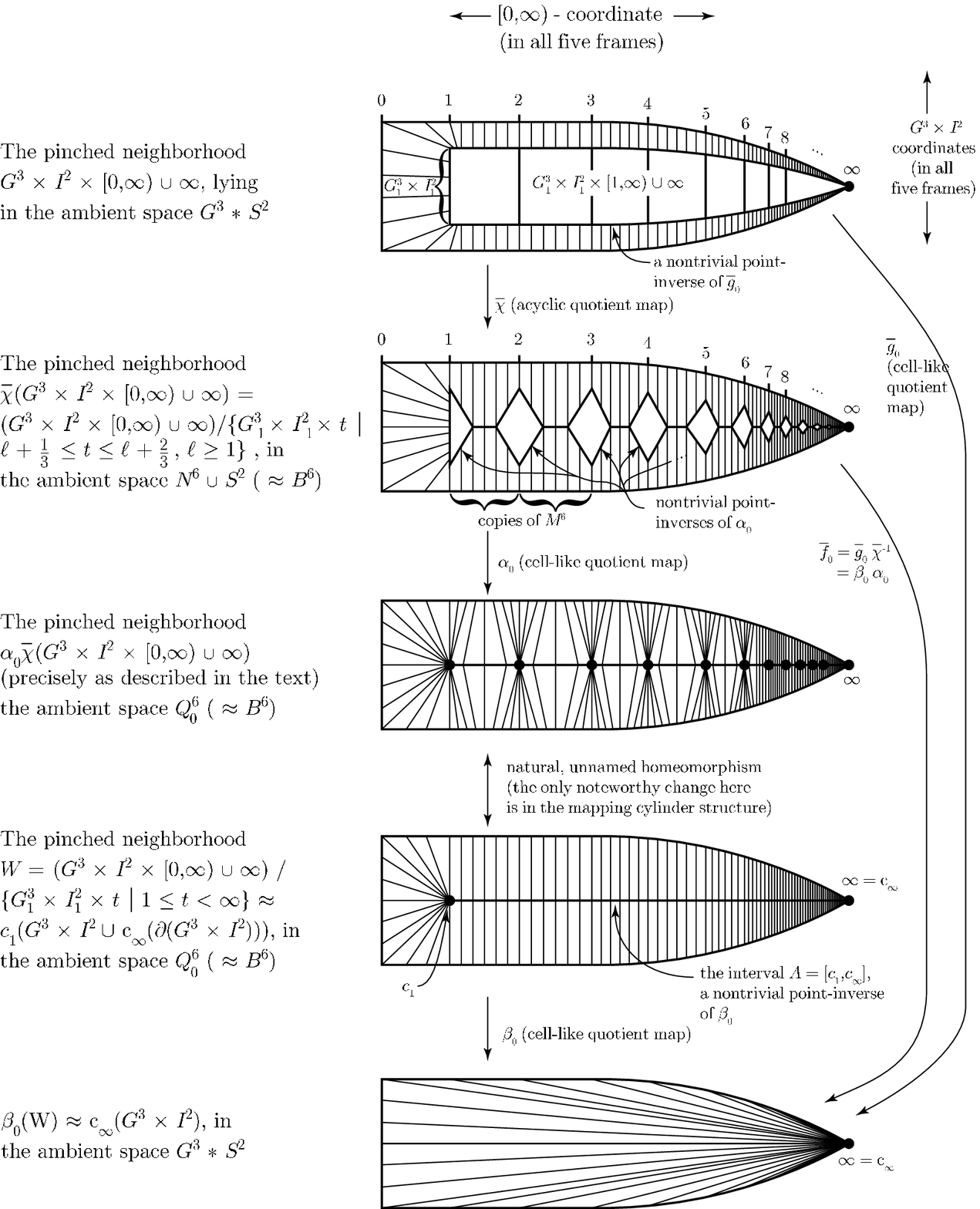}
}\caption{The nontrivial point inverses of the various maps $\bar\chi,
\alpha_{0}, \beta_{0}, \bar{g}_{0}$, and $\bar{f}_{0}$, together with their
pinched neighborhoods}%
\end{figure}

The first map $\alpha_{0}|\overline{\chi}(G_{1}^{3}\times I_{1}^{2}%
\times\lbrack1,\infty)\cup\infty)$ identifies to a point each of the countably
many suspension (or cone) subsets in each of the sets $\overline{\chi}%
(G_{1}^{3}\times I_{1}^{2}\times\lbrack1,\infty)\cup\infty)$, thereby
producing from $\widetilde{N}^{6}\cup S^{2}$ a quotient space $Q_{0}^{6}$. The
second map $\beta_{0}$ completes the remaining identifications by identifying
each set $\alpha_{0}(\overline{\chi}(G_{1}^{3}\times I_{1}^{2}\times
\lbrack1,\infty)\cup\infty))$ to a point. In symbols, this becomes
\begin{multline*}
\overline{\chi}(G_{1}^{3}\times I_{1}^{2}\times\lbrack1,\infty)\cup\infty)=\\
(G_{1}^{3}\times I_{1}^{2}\times\lbrack1,\infty)\cup\infty)/\{G_{1}^{3}\times
I_{1}^{2}\times t\mid\ell+\frac{1}{3}\leq t\leq\ell+\frac{2}{3}\text{, }%
\ell\geq1\}\\
\overset{\left.  \alpha_{0}\right\vert }{\longrightarrow}(G_{1}^{3}\times
I_{1}^{2}\times\lbrack1,\infty)\cup\infty)/(\{G_{1}^{3}\times I_{1}^{2}\times
t\mid\ell+\frac{1}{3}\leq t\leq\ell+\frac{2}{3}\text{, }\ell\geq1\}\cup\\
\{G_{1}^{3}\times I_{1}^{2}\times\lbrack1,\frac{4}{3}]\}\cup\{G_{1}^{3}\times
I_{1}^{2}\times\lbrack\ell+\frac{2}{3},\ell+\frac{4}{3}]\mid\ell
\geq1\})\medskip\\
=[1,\infty]/\{[1,\frac{4}{3}],[\frac{5}{3},\frac{7}{3}],[\frac{8}{3},\frac
{10}{3}],\ldots\}\approx arc\overset{\left.  \beta_{0}\right\vert
}{\longrightarrow}\infty\in S^{2}.
\end{multline*}
Let
\begin{align*}
\alpha &  \equiv\alpha_{0}\cup_{\sigma}id_{B^{6}}:(\widetilde{N}^{6}\cup
S^{2})\cup_{\sigma}B^{6}\ (\approx S^{6})\rightarrow Q^{6}\\
&  \equiv Q_{0}^{6}\cup_{\sigma}B^{6}%
\end{align*}
and let
\[
\beta\equiv\beta_{0}\cup_{\partial}id_{B^{3}\ast S^{2}}:Q^{6}\rightarrow
G^{3}\ast S^{2}\cup_{\partial}B^{3}\ast S^{2}=H^{3}\ast S^{2}.
\]
Then $\overline{f}=\beta\alpha$.\medskip

\textbf{Claim.}\emph{ }$\alpha$\emph{ is approximable by
homeomorphisms.\medskip}

\textit{Proof.} This follows because $\alpha_{0}|\widetilde{N}^{6}$ is
arbitrarily majorant-closely approximable by homeomorphisms, which in turn
follows because $\alpha_{0}|\widetilde{N}^{6}$ satisfies the Bing Shrinking
Criterion, namely, $\alpha_{0}|\widetilde{N}^{6}$ has a countable discrete
collection of nontrivial point-inverses, each of which is cellular. This is
because, recalling the present description of $\widetilde{N}^{6}$ as a
quotient of $G^{3}\ast S^{2}-S^{2}\approx G^{3}\times\mathbb{R}^{2}$, the
nontrivial point-inverses of $\alpha_{0}|\widetilde{N}^{6}$ are each
homeomorphic to $c(G_{1}^{3}\times I_{1}^{2})$ or to $\Sigma(G_{1}^{3}\times
I_{1}^{2})$, each of which has a deleted neighborhood homeomorphic to
$(G^{3}\times I^{2}\cup c(\partial(G^{3}\times I^{2})))\times\mathbb{R}^{1}$
or to $\Sigma(\partial(G^{3}\times I^{2}))\times\mathbb{R}^{1}$, which are
each homeomorphic to $S^{5}\times\mathbb{R}^{1}$ (cf. earlier remarks on
$L^{5}\times\mathbb{R}^{1})$.

We now turn our attention to $\beta:Q^{6}\rightarrow H^{3}\ast S^{2}$, where
now we know from the preceding Claim that $Q^{6}$ is a manifold, homeomorphic
to $S^{6}$. The nontrivial point inverses of $\beta$ consist of a countable
null sequence of arcs. Inspection reveals that we can regard any one of these
arcs, $A$ say, as having a compact pinched neighborhood $W$ in $Q$, described
by the following model:
\begin{align*}
A  &  \equiv\{G_{1}^{3}\times I_{1}^{2}\times t\mid1\leq t<\infty\}\cup
\infty\ (\approx\lbrack1,\infty])\subset W\\
&  \equiv(G^{3}\times I^{2}\times\lbrack0,\infty)/\{G_{1}^{3}\times I_{1}%
^{2}\times t\mid1\leq t<\infty\}\cup\infty.
\end{align*}
There is a useful cone structure on $W$, displayed by writing
\[
W=c_{1}(G^{3}\times I^{2}\cup c_{\infty}(\partial(G^{3}\times I^{2}))),
\]
where $G^{3}\times I^{2}$ here is to be regarded as
\[
G^{3}\times I^{2}\times0\cup\partial(G^{3}\times I^{2})\times\lbrack0,1]
\]
in the first description, where $c_{\infty}$ denotes coning to the point
$\infty$, and where $c_{1}$ denotes coning to a point labelled $1$. In this
latter description, the arc $A$ becomes the interval $A=[c_{1},c_{\infty
}]\subset W$. Note that $W$ has a natural manifold interior,
$\operatorname*{int}W={\mathring{c}}_{1}(G^{3}\times I^{2}\cup c_{\infty
}(\partial(G^{3}\times I^{2})))$, but the \textquotedblleft
boundary\textquotedblright\ (i.e., the base of the cone) is not necessarily a
manifold. From the definitions, $A\cap\operatorname*{int}W=A-c_{\infty}$, and
also $A-(c_{1}\cup c_{\infty})$ has a product neighborhood in
$\operatorname*{int}W$ (but this neighborhood, after deleting the core, need
not be simply connected, as it has the homotopy type of $\partial(G^{3}\times
I^{2}))$.

Using this model, we can describe $Q^{6}$ and $\beta:Q^{6}\rightarrow
H^{3}\ast S^{2}$ this way. To obtain $Q^{6}$, start with $G^{3}\ast S^{2}%
\cup_{\partial}B^{6}$, and replace each of the countably many nontrivial
$g_{0}$-point-inverse sets, i.e. each of the pinched neighborhoods
$G^{3}\times I^{2}\times\lbrack0,\infty)\cup\infty$, with a copy of $W$ in the
obvious manner. The map $\beta$ is the quotient map which identifies to a
point each arc $A$ in each copy of $W$ in $Q^{6}$. We will denote this
countable null collection of arcs and their pinched neighborhoods in $Q^{6}$
by $\{A_{\ell}\}$ and $\{W_{\ell}\}$.

To complete the proof, it remains to establish:\medskip

\textbf{Proposition.}\emph{ }$\beta:Q^{6}\rightarrow H^{3}\ast S^{2}$\emph{ is
approximable by homeomorphisms. That is, the above decomposition of }$Q^{6}%
$\emph{ into arcs }$\{A_{\ell}\}$\emph{ and points is shrinkable. \medskip}

The inspiration for this proposition and its proof comes from a theorem of
Bing \cite[Thm. 3]{Bi2}, who showed that a decomposition of a manifold whose
nontrivial elements comprise a countable (not necessarily null) collection of
flat arcs is a shrinkable decomposition. (Bing did not require his collection
to be null, because he gave an ingenious but simple argument about limits
\cite[p. 368]{Bi2}, which in effect let him assume that his collection was
null. Since our collection above is null to begin with, we will not need this
part of his argument.)

The heart of Bing's argument is a certain technical proposition \cite[Lemma
4]{Bi2} dealing with shrinking individual arcs, which we state here in a form
applicable to our situation.

\begin{lemma}
[Shrinking an individual arc]\label{shrinking-an-individual-arc}\emph{Suppose
}$A_{0}$\emph{ is one of the arcs in the collection }$\{A_{\ell}\}$\emph{, and
suppose }$\epsilon>0$\emph{ is given. Then there exists a neighborhood }%
$U$\emph{ of }$A_{0}$\emph{ in }$N_{\epsilon}(A_{0})$\emph{ (the }$\epsilon
$\emph{-neighborhood of }$A_{0}$\emph{ in }$Q^{6}$\emph{) and there exists a
homeomorphism }$h:Q^{6}\rightarrow Q^{6}$\emph{, supported in }$U$\emph{, such
that for any arc }$A_{\ell}$\emph{, if }$A_{\ell}\cap U\neq\emptyset$\emph{,
then }$\operatorname*{diam}h(A_{\ell})<\epsilon$\emph{.}
\end{lemma}

Given this lemma, it is an easy matter to establish that the Bing Shrinking
Criterion holds for the decomposition $\{A_{\ell}\}$, for one only has to
apply the lemma to those finitely many members of $\{A_{\ell}\}$ whose
diameters are bigger than some preassigned $\epsilon>0$.

To prove his original version of the Lemma, Bing used the flatness of his arc
$A_{0}$. It turns out that our nonflat arcs have enough regularity so that all
the requisite motions of the lemma can be performed. In order to describe
these motions, we make a definition. An embedded arc $A$ (parametrized by
$[1,\infty]$, here) in a manifold-without-boundary $Q$ has the \emph{almost
covering retraction} property (toward the $1$-end, in the definition here)
provided that for any $a\in\lbrack1,\infty]$ and any $\epsilon>0$, there is a
homeomorphism $h:Q\rightarrow Q$, supported in $N_{\epsilon}([a,\infty])$,
such that the restriction $h|[1,\infty]$ is $\epsilon$-close to the
retraction-inclusion map $[1,\infty]\rightarrow\lbrack1,a]\hookrightarrow Q$.
(This concept has been used profitably by several people, among them
Bryant-Seebeck \cite{BS}, Cantrell-Lacher \cite{CL}, Bryant \cite{Bry} and
Price-Seebeck \cite{PS}.) The simplest example of an arc with this property is
a flat arc, or more generally, any subarc of an open arc having a product
neighborhood. Also, we have:

\begin{lemma}
\label{almost-covering-retraction-property}\textit{Suppose $A_{0}$ is one of
the arcs in the collection $\{A_{\ell}\}$. Then $A_{0}$ has the almost
covering retraction property, toward the 1-end}.
\end{lemma}

\textbf{Note.} If $\pi_{1}(G^{3})\neq1$, then $A_{0}$ cannot have the almost
covering retraction property toward the $\infty$-end, as a straightforward
fundamental group analysis shows.

\begin{proof}
[Proof of Lemma \ref{shrinking-an-individual-arc} from Lemma
\ref{almost-covering-retraction-property}]The following argument is adapted
from \cite[Lemma 4]{Bi2}. Pictures are indispensable (but are not provided
here). Given $A_{0}\approx\lbrack1,\infty]$ and given $\epsilon>0$, choose a
partition $1=a_{1}<a_{2}<\ldots<a_{n}<a_{n+1}=\infty$ of $[1,\infty]$ so fine
that each subinterval $[a_{1},a_{i+1}]$ in $Q$ has diameter $<\epsilon/6$. Let
$\delta\in(0,\frac{\varepsilon}{12})$ be so small that $2\delta$ is less than
the minimum distance between any two nonintersecting subintervals in this
partitioning, and such that for any arc $A_{\ell}$ aside from $A_{0}$, if
$A_{\ell}$ intersects $N_{\delta}(A_{0})$, then $\operatorname*{diam}A_{\ell
}<\epsilon/2$. Choose, in the order $U_{n},h_{n},U_{n-1},h_{n-1},\ldots
,U_{1},h_{1}$, a sequence of open subsets $\{U_{i}\}$ of $Q$ and
homeomorphisms $\{h_{i}\}$ of $Q$, satisfying the following conditions
(letting $U_{n+1}=\varnothing$ and $h_{n+1}=identity$). Each $U_{i}$ is a
neighborhood of $[a_{i},\infty]$ such that

\begin{enumerate}
\item $U_{i}-U_{i+1}\subset N_{\delta}([a_{i},a_{i+1}])$ (and hence, applying
induction, $\cup_{j=i}^{n}U_{j}\subset N_{\delta}([a_{i},\infty]))$,

\item (cf. (5) below for $h_{i+1})h_{i+1}(U_{i})\subset N_{\delta}%
([a_{i},a_{i+1}])$ (and hence by (1), $h_{i+1}(U_{i})\cap\bigcup_{j=i+2}%
^{n}U_{j}=\varnothing)$, and

\item no arc $A_{\ell}$ intersects both $U_{i}$ and $(\cup_{j=i+2}^{n}%
U_{j})-U_{i+1}$ (note this latter set misses $A_{0}$).
\end{enumerate}

Each $h_{i}$ is chosen, using the almost covering retraction property for the
retraction $[1,\infty]\rightarrow\lbrack1,a_{i}]$, so that

\begin{enumerate}
\item[(4)] $h_{i}$ is supported in $U_{i}$, and

\item[(5)] $h_{i}([a_{i-1},\infty])\subset N_{\delta}([a_{i-1},a_{i}])$.
\end{enumerate}

The desired homeomorphism of Lemma \ref{shrinking-an-individual-arc} is
$h=h_{n}h_{n-1}\ldots h_{2}h_{1}$, which is supported in $U=\bigcup_{i=1}%
^{n}U_{i}\subset N_{\delta}(A_{0})$. To verify the conclusion, suppose
$A_{\ell}$ is such that $A_{\ell}\cap U\neq\varnothing$. Let $i$ be the least
index such that $U_{i}\cap A_{\ell}\neq\varnothing$. Hence $A_{\ell}$ misses
the supports of $h_{1},\ldots,h_{i-1}$ and so $h(A_{\ell})=h_{n}\ldots
h_{i}(A_{\ell})$. By (3), $A_{\ell}$ misses $(\bigcup_{j=i+2}^{n}%
U_{j})-U_{i+1}$, and so $A_{\ell}\cap U\subset U_{i}\cup U_{i+1}$. Now,
$h_{i+1}h_{i}(U_{i}\cup U_{i-1})\subset U_{i}\cup U_{i+1}$. By (2),
$h_{i+2}(U_{i+1})\subset N_{\delta}([a_{i+1},a_{i+2}])$, and by (1) and (4),
$U_{i}-U_{i+1}$ is not moved by $h_{i+2}$. Consequently $h_{i+2}(U_{i}\cup
U_{i+1})\subset N_{\delta}([a_{i},a_{i+2}])$, and $h_{i+2}(U_{i}\cup U_{i+1})$
misses $\bigcup_{j=i+3}^{n}U_{j}$. So $h(A_{\ell}\cap U)\subset h_{n}\ldots
h_{i+2}(U_{i}\cup U_{i+1})=h_{i+2}(U_{i}\cup U_{i+1})\subset N_{\delta}%
([a_{i},a_{i+2}])$. So $\operatorname*{diam}h(A_{\ell}\cap U)<\epsilon/2$. Now
by choice of $\delta$, $\operatorname*{diam}(A_{\ell}-U)<\epsilon/2$. Hence
Lemma \ref{shrinking-an-individual-arc} is established from Lemma
\ref{almost-covering-retraction-property}
\end{proof}

\begin{proof}
[Proof of Lemma \ref{almost-covering-retraction-property}]Because of the
product and cone structure of $\operatorname*{int}W_{0}$, which in particular
ensures that any subarc [1,a] of $A_{0},\ a<\infty$, has the almost covering
retraction property toward the 1-end, it is clear that Lemma
\ref{almost-covering-retraction-property} can be deduced from the\medskip

\textbf{Claim.} \emph{Given any }$\delta>0$\emph{, there is a homeomorphism
}$h$\emph{ of }$Q$\emph{, supported in }$N_{\delta}(\infty)$\emph{, such that
}$h(A_{0})\subset\operatorname*{int}W_{0}$\emph{ where }$\infty$\emph{ denotes
the }$\infty$\emph{-endpoint of }$A_{0}$\emph{.\medskip}

This claim is established by engulfing. Let $W_{1}\subset(\operatorname*{int}%
W_{0})\cup\infty$ be a copy of $W_{0}$ gotten by squeezing $W_{0}$ inwards a
bit, keeping $\infty$ fixed, using an interior collar of $W_{0}$ which is
pinched at $\infty$. Then $\operatorname*{int}W_{0}\cup\operatorname*{ext}%
W_{1}=Q^{6}-\infty$. Given $\delta>0$, the idea is to use engulfing to produce
two homeomorphisms $h_{0},h_{1}:Q^{6}\rightarrow Q^{6}$, each supported in
$N_{\delta}(\infty)$, such that $h_{0}(\operatorname*{int}W_{0})\cup
h_{1}(\operatorname*{ext}W_{1})=Q^{6}$. Then $h=h_{0}^{-1}h_{1}$ is the
desired homeomorphism of the Claim.

In order to construct these homeomorphisms using dual-skeleton engulfing
arguments, it suffices as usual to construct certain engulfing homotopies
(compare the following engulfing argument to that, in say, \cite{Se}). Letting
$U$ denote either $\operatorname*{int}W_{0}$ or $\operatorname*{ext}W_{1}$,
one wants establish that for any neighborhood $V_{0}$ of the point $\infty$ in
$Q^{6}$, there is a smaller neighborhood $V_{1}$ of $\infty$ and a homotopy
deformation of $U\cup V_{1}$ into $U$ such that all of the motion takes place
in $V_{0}$. Since $\operatorname*{cl}U$ is an ANR, it suffices to show that
for any $\delta>0$, there is a map $\operatorname*{cl}U\rightarrow U$ which is
$\delta$-close to the identity. The existence of this map is clear for
$U=\operatorname*{int}W_{0}$, because $W_{0}$ has an interior collar. For
$U=\operatorname*{ext}W_{1}$, one must analyze the situation more carefully.
Probably the quickest argument is to pass to the target space under the
cell-like map $\beta:Q^{6}\rightarrow H^{3}\ast S^{2}$. Using the fact that
$U$ is saturated with respect to the point-inverses of $\beta$, and that the
nontrivial point-inverses in $U$ have diameter tending to 0 as they approach
$\operatorname*{fr}U$, and using the map lifting property of cell-like maps,
(see for example \cite[Lemmas 2.3, 3.4]{Lac1}, \cite[Thm 1]{Ko} or \cite[Lemma
3]{Hav}), it suffices to show in the target $H^{3}\ast S^{2}$ that for any
$\delta>0$ there is a map $\beta(\operatorname*{cl}U)\rightarrow\beta(U)$
which is $\delta$-close to the identity. The construction of this map is
fairly clear, using the uniform structure of $H^{3}\ast S^{2}$ and the
explicit description of $\beta(\operatorname*{cl}U)$ available from
definitions, and using the interior collar on $B^{3}\ast S^{2}\subset
H^{3}\ast S^{2}$ to make sets disjoint from $G^{3}\ast S^{2}$ by pushing them
into $\operatorname*{int}(B^{3}\ast S^{2})$. This completes the proof of the
Theorem.\medskip
\end{proof}

The above proof does not apply in dimension $5$ for the simple reason that
there it would require the hypothesis that $\Sigma^{2}(H^{3}\#H^{3})\approx
S^{5}$, which is essentially what one is trying to prove. Still, the above
arguments in Lemmas \ref{shrinking-an-individual-arc} and
\ref{almost-covering-retraction-property} seem at first glance almost to work
in dimension $5$ to shrink the nontrivial point-inverses of the original
cell-like map $f:S^{5}\rightarrow\Sigma^{2}H^{3}$ constructed in Part III. But
the difficulty is that the nontrivial point-inverses of $f$, unlike arcs, have
thickness as well as length. Hence the almost covering retraction principle
fails to do the job completely. It would be interesting if such an argument
could be made to work, for that would provide a brief proof of the Double
Suspension Theorem.

\end{document}